# A formal Lie correspondence

by Vincent Bagayoko

IMJ-PRG

*Email:* bagayoko@imj-prg.fr

**Abstract**

We establish an equivalence between categories of "formally nilpotent" Lie algebras and exponential groups in characteristic zero. It extends the equivalences of Mal'cev, Lazard, Quillen and Warfield, and applies to groups under composition of generalized formal series or automorphisms of algebras of generalized formal series. We obtain first-order transfer results from finite dimensional nilpotent objects to formally nilpotent ones. We give applications to solving equations over groups, to the theory of nilpotent exponential groups as per Miasnikov-Remeslennikov, and to definability problems in certain groups of formal series.

# Introduction

**Lie theory in the formal realm** Elements of Lie theory frequently appear when studying algebras and groups of generalized formal series [20, 26, 31, 22, 24] and maps between them that preserve the formal structure. This includes Écalle's work on transseries [16] and alien derivations [17] in relation to Dulac's problem [14], various works on the problem of normalization [38, 37, 36] or identification of periodic points [29] in dynamical systems. More recently [7], a fragment of the formal Lie correspondence was used to obtain an equivalence between certain derivations and automorphisms on formal series, subsuming earlier instances (e.g. [15, Lemma 8.1], [35, Section 1.3]). Lastly, a more explicit interplay between Lie algebras of derivations on formal series and groups under compositions of formal series was used [4] to extend the results on normalization of [37] to more general fields of formal series. Our overarching goal is to study the model theory of certain groups of formal series and germs definable in o-minimal structures (see [6]) requires an understanding of their intricate algebraic structure. We see Lie theory as providing the natural language for such a project.

No framework exists yet for stating formal versions of the classical theorems of Lie theory that would apply to all algebras and groups of generalized formal series. In this paper, we provide one that includes these, as well as classical nilpotent or topologically nilpotent objects. It covers what can be interpreted as the purely formal content of Lie theory, and applies to objects that can be construed as formally nilpotent. These turn out to be Lie algebras equipped with infinite sums [18, 7] and groups equipped with infinite linearly ordered products [3].

**Past correspondences** A Lie correspondence is an equivalence of categories between a category of Lie algebras and a category of groups, with structure preserving homomorphisms as arrows. Beyond the original Lie correspondence [28, 10] between connected simply connected Lie groups and finite dimensional Lie algebras over $\mathbb{C}$ or $\mathbb{R}$, a number of correspondences have been established, in particular for nilpotent-like objects. Using analytic methods, Mal'cev obtained [32, Theorem 9] a correspondence between uniquely divisible nilpotent groups and nilpotent rational Lie algebras. Lazard gave [27, Theorems 4.2 and 4.3] a more explicit correspondence based on valuation theoretic methods, for a class of groups slightly more general than that of lower central complete groups (see Definition 1.9). For the latter, the correspondence was proved independently by Quillen [40, Appendix 3, Theorem 3.3] using Hopf algebras. Stewart gave [46, Theorem 2.4.2] a matrix-based treatment of the correspondence that applies to locally nilpotent uniquely divisible groups and locally nilpotent rational Lie algebras.

These correspondences present a few shortcomings in our perspective. The first one is that they apply to rational Lie algebras, as opposed to Lie algebras over more general fields. This is due in part to the fact that an appropriate notion of "groups with exponents in a ring", now called exponential groups, only emerged [33] after refining previous notions [30, 8], and in part to the fact that the algebraic structure of exponential groups is more involved than that of pure groups. This is all the more problematic as rings of exponents in groups of formal series may be interpretable in the pure group language (see Proposition 6.10), and should therefore be added as a primitive from the standpoint of elimination theory. In the case of nilpotent objects, this issue was solved by Warfield [47, Theorem 12.11].

A second limitation is the breadth of the categories of objects involved. Indeed the main examples of groups of formal series that we are interested in are perfect (see Example 5.28). Even for examples that are residually nilpotent, the topological notion of infinite products coming from the lower central filtration is coarser (i.e. has fewer multipliable families) than our notion of multipliability, and thus cannot always be applied (see Remarks 2.14 and 4.15). Lastly, categories of nilpotent objects are not cocomplete, while the categories of locally nilpotent ones are neither complete nor cocomplete.

Our correspondence extends the above equivalences while solving these issues. It is orthogonal (see Remark 5.27) to the Milnor-Moore correspondence [34] between formal groups and arbitrary Lie algebras.

**Outline of the paper**  We introduce in Section 1 our conventions for algebras and exponential groups, all over a fixed commutative and unital $\mathbb{Q}$-algebra $k$. Section 2 defines summability modules, that are minor generalizations of strong spaces [18] and summability spaces [7]. Summability modules are modules $V$ over $k$ equipped with an array $(\Sigma_I)_{I \in \mathbf{Set}}$ of partial linear maps $\Sigma_I : V^I \longrightarrow V$ seen as summation operators, that share properties of the array of finite summation operators $\Sigma_I^{\min} : V^{(I)} \longrightarrow V$ sending a finitely supported map $v : I \longrightarrow V$ to the finite sum $\sum_{i \in \mathrm{supp}\, v} v(i)$ (see Remark 2.3 for set theoretic concerns). Typical non-trivial examples of summability modules are modules of Hahn series with formal summability [23, 43, 24], or rings of power series $k[[x]]$ with the valuation theoretic notion of sum. Linear maps between summability modules that commute with the summation operators are said strongly linear.

In Section 3, we define summability algebras, as summability modules with a strongly bilinear internal law. We consider infinitary versions $k\langle\!\langle I \rangle\!\rangle^{\mathrm{na}}$ of the free non-associative algebra over a set of formal variables $X_i, i \in I$, as well as their associative algebra and Lie algebra versions $k\langle\!\langle I \rangle\!\rangle$ and $\mathrm{Lie}\langle\!\langle I \rangle\!\rangle$ respectively. The Lie algebra objects of the equivalence form the category $\mathbf{\Sigma Lie}^{\mathrm{ev}}$ of summability Lie algebras with evaluations. These are summability Lie algebras $(L, \Sigma)$ equipped with strongly linear evaluation morphisms $\mathrm{ev}_a^{\mathrm{Lie}} : \mathrm{Lie}\langle\!\langle I \rangle\!\rangle \longrightarrow L$ with $\mathrm{ev}_a^{\mathrm{Lie}}(X_i) = a(i)$ for all $i \in I$, for each $a \in \mathrm{dom}\, \Sigma_I$. The definition extends to summability algebras. We show that locally nilpotent algebras (Proposition 3.25) and residually nilpotent algebras that are complete for the lower central filtration (Example 3.26) have evaluations. We also show that algebras of Noetherian series (Corollary 3.31) and certain algebras under composition of strongly linear maps on modules of Noetherian series have evaluations (Proposition 3.30). Lastly, we prove:

**Theorem 1 (completeness).** [Theorem 3.44] *The category $\mathbf{\Sigma Lie}^{\mathrm{ev}}$ of summability Lie algebras with evaluations is complete and cocomplete.*

Section 4 concerns multipliability exponential groups. They are exponential groups $\mathcal{G}$ together with an array $(\Pi_{(I,<)})_{(I,<)}$, indexed by all linearly ordered sets, of partial functions $\Pi_{(I,<)}\colon \mathcal{G}^I \longrightarrow \mathcal{G}$ that share some properties of the finitary ordered product operators $\Pi_{(I,<)}^{\min}\colon \mathcal{G}^{(I)} \longrightarrow \mathcal{G}$. The main examples of such structures are groups $1+\mathfrak{m}$ where $\mathfrak{m}$ is an associative algebra with evaluations, as well as subgroups of such structures, including groups of formal series under composition (see Section 4.2). The formal exponential map $\exp$ of [7] yields a multipliability exponential group $\mathrm{Gr}\langle\!\langle I\rangle\!\rangle := \exp(\mathrm{Lie}\langle\!\langle I\rangle\!\rangle)$. The category of group objects of the formal Lie equivalence is the category $\mathbf{\Pi Gr}^{\mathrm{ev}}$ of multipliability exponential groups with evaluations $\mathrm{ev}_f^{\mathrm{Gr}}\colon \mathrm{Gr}\langle\!\langle I\rangle\!\rangle \longrightarrow \mathcal{G}$, with $\mathrm{ev}_f^{\mathrm{Gr}}(\exp(X_i)) = f(i)$ for all $i \in I$ whenever $f \in \mathrm{dom}\,\Pi_{(I,<)}$. We have $1+\mathfrak{m} \in \mathbf{\Pi Gr}^{\mathrm{ev}}$ if $\mathfrak{m} \in \mathbf{\Sigma Lie}^{\mathrm{ev}}$ (Proposition 4.25).

The formal Lie correspondence is proved in Section 5:

**Theorem 2 (formal Lie correspondence).** [Theorem 5.16] *There is an isomorphism of categories between $\mathbf{\Sigma Lie}^{\mathrm{ev}}$ and $\mathbf{\Pi Gr}^{\mathrm{ev}}$.*

The proof essentially relies on the single-object equivalence between $\mathrm{Gr}\langle\!\langle I\rangle\!\rangle$ and $\mathrm{Lie}\langle\!\langle I\rangle\!\rangle$. It is based on the Baker-Campbell-Hausdorff formula (see [44, Section I.7]):

$$\log(\exp(X)\cdot\exp(Y)) = X+Y+\frac{1}{2}[X,Y]+\cdots,$$

which is a series in Lie brackets $[P,Q]=PQ-QP$, and on its inverse formulas [32, 27, 46]:

$$\begin{aligned}
\exp(X+Y) &= \exp(X)\exp(Y)\,[\![\exp(X),\exp(Y)]\!]^{1/2}\cdots \\
\exp([X,Y]) &= [\![\exp(X),\exp(Y)]\!]\,[\![\exp(X),[\![\exp(X),\exp Y]\!]\,]\!]^{1/2}\cdots
\end{aligned}$$

which are infinite ordered products in rational powers of commutators $[\![S,T]\!]=S^{-1}T^{-1}ST$. We obtain the general equivalence by evaluating into objects of $\mathbf{\Sigma Lie}^{\mathrm{ev}}$ and $\mathbf{\Pi Gr}^{\mathrm{ev}}$. We show that the formal Lie correspondence extends the previously known ones:

**Theorem 3 (specializations).** [Theorem 5.26] *The formal Lie correspondence specializes to isomorphisms between the categories with the following objects*

*a)* *nilpotent Lie algebras and exponential groups,*

*([47, Theorem 12.11])*

*b)* *locally nilpotent Lie algebras and exponential groups,*

*([46, Theorem 2.4.2] when $k=\mathbb{Q}$)*

*c)* *lower central complete Lie algebras and exponential groups.*

*([40, Appendix 3, Theorem 3.3] when $k=\mathbb{Q}$)*

In Section 6, we give applications of the correspondence. We work in a first-order, purely functional language $\mathcal{L}_k$ of so-called mixed structures, which are objects of $\mathbf{\Sigma Lie}^{\mathrm{ev}}$ or $\mathbf{\Pi Gr}^{\mathrm{ev}}$ equipped with the primitives of both languages of Lie algebras and exponential groups. Terms in $\mathcal{L}_k$, although they act as Lie polynomials or group words in nilpotent mixed structures, involve infinitary sums and products in general (see Section 6.1). We first show that all mixed structures share some first-order properties that are common to all nilpotent mixed structures:

**Theorem 4.** [Corollary 6.6] *Suppose that $k$ is a field of characteristic $0$. Let $\varphi \in \mathcal{L}$ be a sentence of the form $\forall \overline{x}\exists!\overline{y}(\theta(\overline{x},\overline{y}))$ where $\theta(\overline{x},\overline{y})$ is a positive boolean combination of atomic formulas. If $\varphi$ holds in all finite dimensional nilpotent mixed structures, then $\varphi$ holds in all mixed structures.*

Using a generalisation [2, Theorem 7.10] of a result of Smel'kin [45], we obtain:

**Theorem 5.** [Corollary 6.8] *Suppose that $k$ is a field of characteristic $0$. Let $\mathcal{G} \in \mathbf{\Pi Gr}^{\mathrm{ev}}$, let $n > 0$, let $g_1, \ldots, g_n \in \mathcal{G}$ and $\lambda_1, \ldots, \lambda_n \in k$ with $\lambda_1 + \cdots + \lambda_n \neq 0$. There is a unique $f \in \mathcal{G}$ with $g_1 f^{\lambda_1} \cdots g_n f^{\lambda_n} = 1$.*

We show (Proposition 6.9) that nilpotent exponential groups in our sense satisfy the Hall-Petresco identities included in Warfield's definition of nilpotent exponential groups. Together with Example 1.7, this entails that actual Lie groups over $k \in \{\mathbb{R}, \mathbb{C}\}$ that are nilpotent are $k$-nilpotent groups as per Warfield. We show (Proposition 6.10) that the ordered exponential group of finitely nested hyperseries of [5] is intepretable in its pure group reduct. We then establish the central formula in mixed structures (Proposition 6.11), and show (Proposition 6.12) that multipliability exponential groups with evaluations that are $n$-Engel for $n > 0$ are nilpotent.

The paper requires no knowledge in model theory beyond some familiarity with first-order languages and interpretations.

# Table of contents

# 1  Preliminaries

Throughout the paper, we fix an associative, commutative and unital $\mathbb{Q}$-algebra $k$. In certain applications, it will be taken to be a field of characeristic 0. We write **Set** for the category of sets with functions as arrows, and $\mathrm{id}_S$ for the identity arrow $S \longrightarrow S$ on an $S \in \mathbf{Set}$. We write **On** for the class of all ordinals. We sometimes identify a number $m \in \mathbb{N}$ with the set $\{0, \ldots, m-1\}$.

## 1.1  Orderings

Our orderings will always be strict: an ordering is an anti-reflexive and transitive relation. A linear ordering is an ordering which is total. An ordered set is said Noetherian if all its non-empty subsets have a finite, non-zero number of minimal elements. Noetherian linearly ordered sets are exactly well-ordered sets.

    The category of all linearly ordered sets with nondecreasing maps as arrows is denoted by **Los**. The forgetful functor $\mathbf{Los} \longrightarrow \mathbf{Set}$ is denoted $\underline{\cdot} : I \longmapsto \underline{I}$. For $I = (\underline{I}, <) \in \mathbf{Los}$ and $i \in \underline{I}$, we write

$$i^- := (\{j \in \underline{I} : j < i\}, <) \quad \text{and}$$
$$i^+ := (\{j \in \underline{I} : i < j\}, <),$$

We also write $I^* = (\underline{I}, <^*)$ where $<^*$ is the *reverse ordering*

$$\forall i, j \in \underline{I}, (i <^* j \iff j < i).$$

We also sometimes denote $<^*$ by $>$. Let $J = (\underline{J}, <), I = (\underline{I}, <_I) \in \mathbf{Los}$. Consider a linearly ordered set $I$ with ordering $<_I$ and a family $(J_i)_{i \in \underline{I}}$ of linearly ordered subsets $(\underline{J_i}, <)$ of $J$. We write

$$J = \coprod_I (J_i)_{i \in \underline{I}}$$

if we have

$$\underline{J} = \bigcup_{i \in I} \underline{J_i} \quad \text{and} \quad \forall i, i' \in I, (i <_I i' \implies (\forall (j, j') \in \underline{J_i} \times \underline{J_{i'}}, (j < j'))).$$

## 1.2 Monoids and groups

A *magma* is a set $M$ together with a binary law $\cdot : M \times M \longrightarrow M$. A (non-associative) *monoid* is a magma $(M, \cdot)$ together with an element $1 \in M$, such that $1 \cdot m = m \cdot 1 = m$ for all $m \in M$.

Let $(\mathcal{G}, \cdot, 1)$ be a group. Given $f, g \in \mathcal{G}$, we will denote the commutator of $(f, g)$ by

$$[\![f, g]\!] := f^{-1} g^{-1} f g. \qquad (1.1)$$

Note that $(\mathcal{G}, [\![\cdot, \cdot]\!], 1)$ is a non-associative monoid. The *opposite group* $\mathcal{G}^{\mathrm{op}}$ of $\mathcal{G}$ is $\mathcal{G}^{\mathrm{op}} = (\mathcal{G}, \cdot^{\mathrm{op}}, 1)$ where $f \cdot^{\mathrm{op}} g = g \cdot f$ for all $f, g \in \mathcal{G}$. If $X$ is a set, then we write $\mathcal{G}^X$ for the group, under pointwise product, of functions $I \longrightarrow \mathcal{G}$. The *support* of an $f \in \mathcal{G}^X$ is the set

$$\operatorname{supp} f := \{x \in X : f(x) \neq 1\}.$$

We denote by $\mathcal{G}^{(X)}$ the subset of elements $f \in \mathcal{G}^X$ whose support is finite, which is a normal subgroup of $\mathcal{G}^X$.

**Remark 1.1.** We use double brackets for commutators in order to distinguish this from the frequently occuring Lie brackets in Lie algebras. The standard convention (1.1) for the commutator (as opposed to $[\![f, g]\!] = f g f^{-1} g^{-1}$) is that adopted by Mal'cev and Stewart [32, 46] but not by Lazard [27].

## 1.3 Algebras

An *algebra* (over $k$) is a $k$-module $(A, +, 0, (a \mapsto \lambda a)_{\lambda \in k})$ together with a map

$$\cdot : A \times A \longrightarrow A$$

which is $k$-bilinear. The first-order theory of algebras in the language

$$\mathcal{L}_{k\text{-alg}} = \langle +, 0, (\lambda \, . \,)_{\lambda \in k}, \cdot \rangle$$

is denoted Alg. We understand ideals of algebras to be two-sided.

We now fix an algebra $(A, +, 0, (a \mapsto \lambda a)_{\lambda \in k}, \cdot)$. We say that $A$ is *abelian* if $a \cdot b = b \cdot a$ for all $a, b \in A$. We write Ab for the first-order theory of abelian algebras. We say that $A$ is *associative* iif $a \cdot (b \cdot c) = (a \cdot b) \cdot c$ for all $a, b, c \in A$. We write Ass for the first-order theory of associative algebras. We say that $A$ is a *Lie algebra* if

$$\begin{aligned} a \cdot a &= 0 \text{ and} \\ a \cdot (b \cdot c) + b \cdot (c \cdot a) + c \cdot (a \cdot b) &= 0 \end{aligned}$$

for all $a, b, c \in A$. We write Lie for the first-order theory of Lie algebras. We have $a \cdot b = -b \cdot a$ for all $a, b \in A$, so $A$ is abelian if and only if $\cdot$ is identically zero. If $A$ associative algebra, then we write $[A]$ for the algebra $(A, +, 0, (a \mapsto \lambda a)_{\lambda \in k}, [\cdot, \cdot])$ where $[\cdot, \cdot]$ is the Lie bracket

$$\forall s, t \in A, [s, t] = s \cdot t - t \cdot s,$$

which is a Lie algebra. We say that $A$ is *unital* if there is a $1 \in A \setminus \{0\}$ such that $(A \setminus \{0\}, \cdot, 1)$ is a monoid.

Let $A$ be an algebra. We set $A_1 = A$, and for each $n > 0$, we define $A_{n+1}$ as the submodule of $A$ generated by elements of the form $a \cdot b$ and $b \cdot a$ where $a \in A_k$ and $b \in A_{n+1-k}$ for each $k \in \{1, \ldots, n\}$. Each $A_n$ is an ideal, and $(A_n)_{n>0}$ is called the *lower central series* of $A$. For $m > 0$, we say that $A$ is $m$-*nilpotent* if $A_{m+1} = \{0\}$, that it is nilpotent if it is $n$-nilpotent for a certain natural number $n$. We say that $A$ is *residually nilpotent* if $\bigcap_{n \in \mathbb{N}} A_n = \{0\}$. Lastly, we say that $A$ is *locally nilpotent* if every finitely generated subalgebra of $A$ is nilpotent. The *lower central completion* $\tilde{A}$ of a residually nilpotent algebra $A$ is the inverse limit $\tilde{A} = \varprojlim (A/A_n)_{n>0}$.

**Definition 1.2.** *A residually nilpotent algebra $A$ is said* l**ower central complete** *if the natural embedding $A \longrightarrow \tilde{A}$ is surjective.*

**Remark 1.3.** *If $A$ is a Lie algebra, then for $n > 0$, one can equivalently define $A_{n+1}$, to be the submodule spanned by elements $a \cdot b$ only, for $a \in A$ and $b \in A_n$.*

## 1.4 Exponential groups

**Definition 1.4.** [33, Definition 2] *An* **exponential group** *(over $k$), or $k$-group, is a group $(\mathcal{G}, \cdot, 1)$ together with a* **power map**, *i.e. a function*

$$\begin{aligned} \cdot^{\cdot} : k \times \mathcal{G} &\longrightarrow \mathcal{G} \\ (\lambda, f) &\longmapsto f^{\lambda} \end{aligned}$$

*satisfying the following axioms for all $f, g \in \mathcal{G}$ and $\lambda, \mu \in k$.*

- **EG1.** $f^1 = f$.
- **EG2.** $f^{\lambda+\mu} = f^{\lambda} f^{\mu}$.
- **EG3.** $(f^{\lambda})^{\mu} = f^{(\lambda\mu)}$.
- **EG4.** *if $fg = gf$, then $(fg)^{\lambda} = f^{\lambda} g^{\lambda}$.*
- **EG5.** $(fgf^{-1})^{\lambda} = fg^{\lambda}f^{-1}$.

We write $\mathcal{L}_{k\text{-gr}} = \langle 1, \cdot, (\cdot^{\lambda})_{\lambda \in k} \rangle$ for the first-order language of exponential groups, i.e. the first-order language of groups expanded by a unary function symbol $\cdot^{\lambda}$ for each $\lambda \in k$. Exponential groups are construed $\mathcal{L}_{k\text{-gr}}$-structures in the expected way.

**Remark 1.5.** *For $k = \mathbb{Q}$, exponential groups are exactly uniquely divisible groups, i.e. groups in which any element $g$ has a unique $n$-th root $g^{1/n}$ for each $n \in \mathbb{N} \setminus \{0\}$. The unique power map is given by $(m/n, g) \mapsto (g^{1/n})^m$, for all $m \in \mathbb{Z}$ and $n \in \mathbb{N} \setminus \{0\}$.*

**Remark 1.6.** *The notion of exponential group is defined for any unital associative ring $k$, and generalizes to the case of possibly non-abelian group the notion of right $k$-module.*

**Example 1.7.** *Let $k \in \{\mathbb{R}, \mathbb{C}\}$ and let $G$ be a nilpotent real or complex (connected and simply connected) Lie group with lie algebra $\mathfrak{g}$. Then the exponential map $\exp: \mathfrak{g} \longrightarrow G$ is bijective. The operation $(a, b) \mapsto \log(\exp(a)\exp(b))$ on $\mathfrak{g}$ is [9, Chapter 3, § 9 Proposition 13] the Baker-Campbell-Hausdorff operation (3.6). Therefore the group $G$, together with the power map $(\lambda, g) \mapsto \exp(\lambda \log(g))$, is isomorphic to the exponential group $\mathrm{Gr}(\mathfrak{g})$ given by the formal Lie correspodence (Theorem 5.16). In particular, it is an exponential group over $k$.*

Let $\mathcal{G}, \mathcal{H}$ be exponential groups. A homormorphism of exponential groups $\varphi : \mathcal{G} \longrightarrow \mathcal{H}$ is a group homomorphism such that $\varphi(f^{\lambda}) = \varphi(f)^{\lambda}$ for all $f \in \mathcal{G}$ and all $\lambda \in k$, i.e. an $\mathcal{L}_{k\text{-gr}}$-homomorphism. An exponential subgroup of $\mathcal{G}$ is a subgroup of $\mathcal{G}$ that is closed under the operations $g \mapsto g^{\lambda}$ for all $\lambda \in k$, i.e. an $\mathcal{L}_{k\text{-gr}}$-substructure of $\mathcal{G}$. An ideal of an exponential group $\mathcal{G}$ is the kernel of a homomorphism of exponential groups with domain $\mathcal{G}$. Equivalently [33, Proposition 5], it is a normal exponential subgroup of $\mathcal{G}$ such that for all $f, g \in \mathcal{G}$ and $\lambda \in k$, we have

$$[\![f, g]\!] \in \mathcal{H} \Longrightarrow f^{-\lambda} g^{-\lambda} (fg)^{\lambda} \in \mathcal{H}. \tag{1.2}$$

The category **Gr** of exponential groups over $k$ with homomorphisms of exponential groups has [33, Theorem 4] free objects $\mathrm{Free}(I)$ for each set $I$.

The lower central series $(\mathcal{G}_n)_{n\in\mathbb{N}}$ of an exponential group $\mathcal{G}$ is defined inductively as follows. Set $\mathcal{G}_1 = \mathcal{G}$. For $n > 0$, define $\mathcal{G}_{n+1}$ to be the smallest ideal of $\mathcal{G}$ containing the set $[\![\mathcal{G}, \mathcal{G}_n]\!] := \{[\![f, g]\!] : (f, g) \in \mathcal{G} \times \mathcal{G}_n\}$. An exponential group is said $m$-nilpotent for $m > 0$ if $\mathcal{G}_{m+1} = \{1\}$, nilpotent if it is $m$-nilpotent for a certain $m > 0$, and *locally nilpotent* if all its finitely generated exponential subgroups are nilpotent. It is said *residually nilpotent* if $\bigcap_{n>0} \mathcal{G}_n = \{0\}$. The *lower central completion* of a residually nilpotent exponential group $\mathcal{G}$ is the inverse limit $\tilde{\mathcal{G}} := \varprojlim (\mathcal{G}/\mathcal{G}_n)_{n>0}$.

**Lemma 1.8.** *Let $\mathcal{G}$ be a residually nilpotent exponential group. Then there is a unique power map on $\tilde{\mathcal{G}}$ for which it is the inverse limit in* **Gr** *of the system $(\mathcal{G}/\mathcal{G}_n)_{n>0}$.*

**Proof.** We have a power map

$$k \times \tilde{\mathcal{G}} \longrightarrow \tilde{\mathcal{G}}$$
$$(\lambda, (f_n \mathcal{G}_n)_{n>0}) \longmapsto (f_n^\lambda \mathcal{G}_n)_{n>0}.$$

The axioms **EG1**, **EG2**, **EG3** and **EG5** are preserved by inverse limits as they are purely equational universal sentences. The axiom **EG4** however requires some care. Let $\lambda \in k$ and let $\tilde{f} = (f_n \mathcal{G}_n)_{n>0}$ and $\tilde{g} = (g_n \mathcal{G}_n)_{n>0}$ commute in $\tilde{\mathcal{G}}$, So $[\![f_n, g_n]\!] \in \mathcal{G}_n$ for each $n > 0$. Let $n > 0$. As $\mathcal{G}_n$ is an ideal, it is closed under the operation of (1.2), so $[\![f_n, g_n]\!] \in \mathcal{G}_n$ entails that $f_n^\lambda g_n^\lambda (f_n g_n)^{-\lambda} \in \mathcal{G}_n$. But this means that $\tilde{f}^\lambda \tilde{g}^\lambda = (\tilde{f}\tilde{g})^\lambda$, i.e. **EG4** holds. It is clear that this is the only power map on $\tilde{\mathcal{G}}$ for which it is the desired inverse limit. □

**Definition 1.9.** *A residually nilpotent exponential group is said* **lower central complete** *if the natural embedding $\mathcal{G} \longrightarrow \tilde{\mathcal{G}}$ is surjective.*

**Remark 1.10.** In [1], Amaglobeli and Remeslennikov proposed two additional definitions of $m$-nilpotency and showed that all three coincide for $m \in \{1, 2\}$. Our notion of nilpotency is what they call lower $m$-nilpotency, i.e. the strongest notion of nilpotency.

## 1.5 Free monoids

Let $I$ be a set with $\varnothing \notin I$, thought of as an alphabet, and let ( and ) be elements outside of $I$. The set of associative words on $I$ is the set

$$I^\star := \bigcup_{n \in \mathbb{N}} I^n,$$

where $\varnothing \in I^0$ denotes the empty word. This is a monoid under concatenation $(u, v) \mapsto uv$. On $(I \sqcup \{(,)\})^\star$, we have a parenthesized concatenation $(\cdot, \cdot)$ given by $(u, v) \mapsto (uv)$ for all $u, v \neq \varnothing$ and sending $(u, \varnothing)$ and $(\varnothing, u)$ to $u$ for all $u$. The set of non-associative words on $I$ is the smallest subset of $(I \sqcup \{(,)\})^\star$ containing $\{\varnothing\} \sqcup I$ and closed under $(\cdot, \cdot)$. We denote this set $I^{(\star)}$. This is an inductive structure with a unique writing property, and corresponding rank

$$\mathrm{rank} : I^{(\star)} \longrightarrow \mathbb{N}$$

with $\mathrm{rank}(\varnothing) = 0$, $\mathrm{rank}(i) = 1$ for all $i \in I$ and $\mathrm{rank}((uv)) = \max(\mathrm{rank}(u), \mathrm{rank}(v)) + 1$ for all $u, v \neq \varnothing$. The structure $(I^{(\star)}, (\cdot, \cdot), \varnothing)$ is a monoid.

We also define a map content : $I^{(\star)} \longrightarrow \{J \subseteq I : J \text{ is finite}\}$ by induction on the rank, as follows. We set $\text{content}(\varnothing) := \varnothing$, $\text{content}(i) := \{i\}$ for all $i \in I$ and $\text{content}((uv)) := \text{content}(u) \cup \text{content}(v)$ for all $u, v \neq \varnothing$. So $\text{content}(w)$ is the set of letters that occur in $w$ excluding parentheses. Likewise, we define the content $\text{content}(u)$ of an associative word $u$ as the set of letters that occur in $u$, i.e. the smallest subset with $u \in \text{content}(u)^\star$.

Let $\mathcal{M} = (M, \cdot, 1)$ be a non-associative monoid. Let $f : I \longrightarrow M$ be a family. Given a non-associative word $w \in I^{(\star)}$, we define the evaluation of $w$ at $f$ to be the element $w(f) \in M$ defined by induction on the rank $\text{rank}(w)$ as

$$\begin{aligned}\varnothing(f) &:= 1 \\ i(f) &:= f(i) \quad \text{for all } i \in I, \text{ and} \\ (uv)(f) &:= u(f) \cdot v(f) \quad \text{for all } u, v \neq \varnothing.\end{aligned}$$

If $\mathcal{M} = (M, \cdot)$ is only a magma, then we define $u(f)$ as above for all $u \in I^{(\star)} \setminus \{\varnothing\}$ only. If $\mathcal{M} = (\mathcal{G}, [\![\cdot, \cdot]\!], 1)$ is the commutator monoid of a group $(\mathcal{G}, \cdot, 1)$, then we write $u[\![f]\!] := u(f)$ for all $f : I \longrightarrow \mathcal{G}$. If $\mathcal{M} = (L, [\cdot, \cdot])$ is the Lie-bracket magma of a Lie algebra $(L, [\cdot, \cdot])$, then we write $u[f] := u(f)$ for all $f : I \longrightarrow L$.

## 2 Formal infinite sums

In this section, we give a presentation of a version of strong vector spaces [18] or summability spaces [7] that are $k$-modules.

### 2.1 Summability modules

Let $(V, +, 0, (v \mapsto \lambda v)_{\lambda \in k})$ be a $k$-module. For each set $I$, let us be given a partial map $\Sigma_I : V^I \longrightarrow V$ whose domain $\text{dom}\,\Sigma_I$ is a submodule of $V^I$, and which is $k$-linear on $\text{dom}\,\Sigma_I$. Consider the following axiomatic properties for all sets $I$, $J$:

**SM1.** $V^{(I)} \subseteq \text{dom}\,\Sigma_I$ and $\Sigma_I v = \sum_{v(i) \neq 0} v(i)$ for all $v \in V^{(I)}$.

**SM2.** If $\varphi : I \longrightarrow J$ is a bijection, then for all $v \in \text{dom}\,\Sigma_J$, we have $v \circ \varphi \in \text{dom}\,\Sigma_I$ and $\Sigma_J v = \Sigma_I (v \circ \varphi)$.

**SM3.** If $I = \bigsqcup_{j \in J} I_j$ for a family of subsets $(I_j)_{j \in J}$ and $v \in \text{dom}\,\Sigma_I$, then

  **SM3a.** $v \upharpoonright I_j \in \text{dom}\,\Sigma_{I_j}$ for each $j \in J$,

  **SM3b.** $(\Sigma_{I_j}(v \upharpoonright I_j))_{j \in J} \in \text{dom}\,\Sigma_J$, and

  **SM3c.** $\Sigma_J (\Sigma_{I_j}(v \upharpoonright I_j))_{j \in J} = \Sigma_I v$.

**SM4.** If $I \cap J = \varnothing$ and $(v, w) \in \text{dom}\,\Sigma_I \times \text{dom}\,\Sigma_J$, then $v \sqcup w \in \text{dom}\,\Sigma_{I \sqcup J}$.

**SM5.** For $v \in \text{dom}\,\Sigma_I$ and a family of functions $f_i : X_i \longrightarrow k$ with finite domains $X_i$, $i \in I$, we have $(f_i(x)\,v(i))_{i \in I \wedge x \in X_i} \in \text{dom}\,\Sigma_{\{(i,x):i \in I \wedge x \in X_i\}}$.

**Definition 2.1.** *A* **summability module** *is a $k$-module together with a class of functions $\Sigma_I$ for all sets $I$, that satisfy* **SM1**-**SM5**. *If $(V, (\Sigma_I)_{I \in \mathbf{Set}})$ is a summability module, then we say that $\Sigma = (\Sigma_I)_{I \in \mathbf{Set}}$ is a* **summability structure** *on $V$. Given a set $I$, elements of $\text{dom}\,\Sigma_I$ are called* **summable families**.

Given a summability module $(V, \Sigma)$, we sometimes write $\sum_{i \in I} v(i)$ for the sum $\Sigma_I v$ of a family $v \in \text{dom}\,\Sigma_I$.

**Remark 2.2.** In any non-trivial summability module $(V, \Sigma)$, for any $v \in V \setminus \{0\}$ and any infinite set $I$, the constant family $(v)_{i \in I}$ is never summable. Indeed, if it were, fixing an $i_0 \in I$ and a bijection $\varphi: I \longrightarrow I \setminus \{i_0\}$, we would have $\sum_{i \in I} v = \left( \sum_{i \in I \setminus \{i_0\}} v \right) + v$ by **SM3** and **SM1**, where $\left( \sum_{i \in I \setminus \{i_0\}} v \right) = \sum_{i \in I} v$ by **SM2**. But this yields $v = \sum_{i \in I} v - \sum_{i \in I} v = 0$: a contradiction.

**Remark 2.3.** If $f: I \longrightarrow V$ is a family where $|\operatorname{supp} f| > |V|$, then there must be a element $v_0 \in V \setminus \{0\}$ such that the set $J := f^{-1}(\{v_0\})$ is infinite. So $f$ cannot be summable, otherwise by **SM3a** the family $f \restriction J$ would be summable, contradicting Remark 2.2. Thus the supports of summable families in $(V, \Sigma)$ have cardinality bounded by $|V|$, so in view of **SM2**, the proper class $(V, \Sigma)$ can effectively be represented by a set $(V, (\Sigma_I)_{I \subseteq V})$. We do not adopt this convention in order to avoid having to specify injections into $V$ everytime we want to consider some operations on the index sets summable families.

**Remark 2.4.** Whereas the axioms **SM1**–**SM4** allow some analytical notions of summability (e.g. absolute summability of complex-valued sequences), the axiom **SM5** forbids this.

**Remark 2.5.** In [7], the scalar ring $k$ was taken to be a field of characteristic zero. However all the results in that paper apply without change in our context.

**Example 2.6.** Every module $V$ has a minimal summability structure $\Sigma^{\min}$ where $\operatorname{dom} \Sigma_I^{\min} = V^{(I)}$ for all sets $I$.

**Definition 2.7.** *Given a summability module $(V, \Sigma)$, a **summability submodule** of $(V, \Sigma)$ is a submodule $W \subseteq V$ such that the sum of each family $w: I \longrightarrow W$ that is summable in $V$ lies in $W$.*

A summability submodule $W$ of $(V, \Sigma)$ has a natural structure $\Sigma^W$ of summability module where for each set $I$, the function $\Sigma_I^W$ is the restriction $\Sigma_I^W = \Sigma_I \restriction (W^I \cap \operatorname{dom} \Sigma_I)$ of $\Sigma_I$ to $W^I \cap \operatorname{dom} \Sigma_I$.

**Example 2.8.** Let $(V, \Sigma)$ be a summability module and let $W \subseteq V$ be a summability submodule. Then the quotient module $V/W$ has a natural structure $\Sigma^{/W}$ of summability module, where for all sets $I$, we have $\operatorname{dom} \Sigma_I^{/W} = \{ f + W : f \in \operatorname{dom} \Sigma_I \}$ and $\Sigma_I^{/W}(f + W) = (\Sigma_I f) + W$ for all $f \in \operatorname{dom} \Sigma_I$. That this is well-defined follows from the facts that $W$ is a summability submodule and that $\Sigma_I$ is a linear map.

Given two summability modules $(V, \Sigma^V)$ and $(W, \Sigma^W)$, a linear map $\phi: V \longrightarrow W$ is said strongly linear if for each summable family $v: I \longrightarrow V$, the family $\phi \circ v: I \longrightarrow W$ is summable with sum $\Sigma_I^W \phi \circ v = \phi(\Sigma_I^V v)$. For instance, if $W$ is a summability submodule of $(V, \Sigma)$, then the quotient map $V \longrightarrow V/W$ is strongly linear. The $k$-module $\operatorname{Lin}^+(V, W)$ of strongly linear maps $V \longrightarrow W$ is [7, Proposition 1.35] a summability module under the following summability structure: a family $(\varphi_i)_{i \in I}: I \longrightarrow \operatorname{Lin}^+(V, W)$ is summable if for all summable $v: J \longrightarrow V$, the family $(\varphi_i(v(j)))_{(i,j) \in I \times J}: I \times J \longrightarrow W$ is summable. Then $\Sigma_I(\varphi_i)_{i \in I}$ is the map sending $v \in V$ to $\Sigma_I^W(\varphi_i(v))_{i \in I}$. We write **ΣMod** for the category of summability modules, with strongly linear maps as arrows.

## 2.2 Bornologies

Given a set $X$, a *bornology* on $X$ is a filter on $X$ that contains all finite subsets, i.e. a set of subsets of $X$ containing all finite subsets of $X$, that is closed under finite unions and taking subsets. Examples include the set $\operatorname{Fin}(X)$ of all finite subsets of $X$ and the set $\mathfrak{P}(X)$ of all subsets of $X$.

If $(V_x, \Sigma^x)_{x \in X}$ is a family of summability spaces and $\mathcal{F}$ is a filter on $X$, then the subset $\prod_{x \in X}^{\mathcal{F}} V_x$ of the Cartesian product $\prod_{x \in X} V_x$ of maps $f$ whose support $\operatorname{supp} f = \{x \in X : f(x) \neq 0\}$ lies in $\mathcal{F}$, under pointwise sum and scalar product, is a $k$-module. It also has a natural summability structure $\Sigma^{\mathcal{F}}$, where for all sets $I$, the domain of $\Sigma_I^{\mathcal{F}}$ is the set of families $\varphi : I \longrightarrow \prod_{x \in X}^{\mathcal{F}} V_x$ such that:

$$\bigcup_{i \in I} \operatorname{supp} \varphi(i) \in \mathcal{F} \qquad \text{and} \qquad \text{for all } x \in X, \text{ the family } (\varphi(i)(x))_{i \in I} \text{ is in } \operatorname{dom} \Sigma_I^x.$$

The sum of such a family is the map $\Sigma_I^{\mathcal{F}} \varphi$ that sends $x \in X$ to $\Sigma_I^x (\varphi(i)(x))_{i \in I}$. See [18, Section 2.1].

**Example 2.9.** Taking $\mathcal{F}$ to be $\mathfrak{P}(X)$ (resp. $\operatorname{Fin}(X)$), the summability module $\prod_{x \in X}^{\mathcal{F}} V_x$ is the product (resp. coproduct) of the family $(V_x, \Sigma^x)_{x \in X}$ in $\boldsymbol{\Sigma}\mathbf{Mod}$.

We fix a set $X$ and a bornology $\mathcal{F}$ on $X$. We write

$$k(X; \mathcal{F}) = \prod_{x \in X}^{\mathcal{F}} k$$

where $k$ is equipped with the minimal summability structure $\Sigma^{\min}$. An $f : I \longrightarrow k(X; \mathcal{F})$ is summable if and only if $\bigcup_{i \in I} \operatorname{supp} f(i) \in \mathcal{F}$ and for all $x \in X$, the set $\{i \in I : x \in \operatorname{supp} f(i)\}$ is finite.

For $x \in X$, the set $\{x\}$ lies in $\mathcal{F}$, so its indicator function $\mathbb{1}_x : X \longrightarrow \{1, 0\} \subseteq k$ lies in $k(X, \mathcal{F})$. By definition, for all $f \in k(X; \mathcal{F})$, the family $(f(x) \mathbb{1}_x)_{x \in X}$ is summable, with

$$f = \sum_{x \in X} f(x) \mathbb{1}_x. \tag{2.1}$$

**Example 2.10.** Let $(X, <)$ be an ordered set and let $\mathcal{N}$ be the set of Noetherian subsets of $(X, <)$. Then $\mathcal{N}$ is a bornology on $X$ (see [7, Proposition 1.15]), and we call $k(X; \mathcal{N})$ a *module of Noetherian series*.

**Lemma 2.11.** *Let $(V, \Sigma)$ be a summability module. A linear map $\phi : k(X; \mathcal{F}) \longrightarrow V$ is strongly linear if and only if for all $f \in k(X; \mathcal{F})$, the family $(\phi(f(x) \mathbb{1}_x))_{x \in X}$ is summable with sum $\sum_{x \in X} \phi(f(x) \mathbb{1}_x) = \phi(f)$.*

**Proof.** The proof is the same as in [7, Proposition 1.29], where only the properties of bornologies are used. $\square$

Let $X$ be a set, let $\mathcal{F}$ be a bornology on $X$ and let $(V, \Sigma)$ be a summability module. We say that a map $\tilde{\phi} : k(X; \mathcal{F}) \longrightarrow V$ is a strongly linear extension of a map $\phi : X \longrightarrow V$ strongly linear if and only if $\phi$ is strongly linear and $\tilde{\phi}(\mathbb{1}_x) = \phi(x)$ for all $x \in X$.

**Lemma 2.12.** *A map $\phi : X \longrightarrow V$ has a strongly linear extension if and only if for all $U \in \mathcal{F}$, the family $(\phi(u))_{u \in U}$ is summable in $V$. The strongly linear extension is unique, given by $\tilde{\phi}(f) = \sum_{x \in X} f(x) \phi(x)$ for all $f \in k(X; \mathcal{F})$.*

**Proof.** If $\phi$ satisfies the conditions, then the extension $\tilde{\phi}$ is well-defined, and it is strongly linear by Lemma 2.11. The unicity follows from (2.1). Conversely if $\tilde{\phi}$ is well-defined and strongly linear, then for $U \in \mathcal{F}$, as $(\mathbb{1}_u)_{u \in U}$ is summable in $k(X; \mathcal{F})$, so is $(\tilde{\phi}(\mathbb{1}_u))_{u \in U} = (\phi(u))_{u \in U}$. $\square$

**Remark 2.13.** The summability structure on $k^X = k(X; \mathfrak{P}(X))$ is topological, in the sense that a family $(v_i)_{i \in I}$ is summable in $k^X$ if and only if the net of finite partial sums $(\sum_{j \in J} v(j))_{J \subseteq I \text{ is finite}}$ converges in the product topology on $k^X$ for the discrete topology on $k$.

**Remark 2.14.** The summability structure on $k(X;\mathcal{F})$ is also topological, given by the representation of $k(X;\mathcal{F})$ as the direct limit $k(X;\mathcal{F}) = \varinjlim (k^U)_{U \in \mathcal{F}}$ where $k^U$ as a topological space is as in Remark 2.13. Unfortunately, this topological presentation is not preserved when taking summability submodules or quotients by summability submodules (see [18, Introduction]).

## 3 Summability algebras

We now endow summability modules with strongly bilinear maps and study the resulting structure of algebras with infinite sums.

**Definition 3.1.** A **summability algebra** *is an algebra* $(A, +, 0, (a \mapsto c\,a)_{c \in k}, \cdot)$ *together with a summability structure* $\Sigma$ *on the underlying $k$-module, that satisfies the axioms*

**SA.** *For all summable families* $a : I \longrightarrow A$ *and* $b : J \longrightarrow A$, *the family* $a \cdot b := (a(i) \cdot b(j))_{(i,j) \in I \times J} : I \times J \longrightarrow A$ *is summable, with* $\Sigma_{I \times J}\, a \cdot b = (\Sigma_I a) \cdot (\Sigma_I b)$.

A **summability associative algebra** *is a summability algebra whose underlying algebra is an associative algebra*. A **summability Lie algebra** *is a summability algebra whose underlying algebra is a Lie algebra*.

**Example 3.2.** Any algebra $A$ is a summablity algebra for the minimal summability structure $\Sigma^{\min}$ of Example 2.6.

**Example 3.3.** If $(V, \Sigma)$ is a summability module, then the summability module $\operatorname{Lin}^+(V) := \operatorname{Lin}^+(V, V)$ is an associative summability algebra under composition (see [7, Proposition 1.35]).

**Example 3.4.** If $(A, \Sigma)$ is a summability associative algebra, then $([A], \Sigma)$ is a summability Lie algebra.

**Example 3.5.** If $X$ is a set, $\mathcal{F}$ is a bornology on $X$, and $(A_x, \Sigma_x)_{x \in X}$ is a family of summability algebras, then the pointwise product turns $\prod_{x \in X}^{\mathcal{F}} A_x$ into a summability algebra.

**Example 3.6.** Let $A$ be an algebra and let $\mathfrak{q} \subseteq A$ be a residually nilpotent ideal such that $A$ is complete in the (Hausdorff) $\mathfrak{q}$-adic topology on $A$. Then $A \simeq \varprojlim (A/\mathfrak{q}_n)_{n > 0}$ has a natural structure of summability algebra where $a : I \longrightarrow A$ is summable if for all $n > 0$, the set $I_n := \{i \in I : a(i) \notin \mathfrak{q}_n\}$ is finite, and then $\Sigma_I a$ is defined as the limit of the Cauchy sequence $(\Sigma_{i \in I_n} a(i))_{n > 0}$.

**Remark 3.7.** If $k_0 \subseteq k$ is a (unital) subalgebra of $k$, then the reduct of a summability algebra $A$ over $k$ as a $k_0$-algebra is a summability algebra over $k_0$.

Given a summability algebra $(A, \Sigma)$, a *summability subalgebra* of $A$ is a subalgebra $B \subseteq A$ which is a summability submodule of the underlying summability module. A *summability ideal* of $A$ is an ideal of $A$ which is a summability submodule of $A$.

**Example 3.8.** Let $(A, \Sigma)$ be a summability algebra and let $\mathfrak{m} \subseteq A$ be a summability ideal. Then $A/\mathfrak{m}$ with the summability structure of Example 2.8 is a summability algebra.

## 3.1 Algebras of Noetherian series

**Definition 3.9.** *An ordered non-associative monoid is a non-associative monoid $(M,\cdot,1)$ together with a partial ordering $<$ such that for all $x,y,z \in M$, we have*

$$x < y \Longrightarrow (z\,x < z\,y \wedge x\,z < y\,z).$$

Note that the strict ordering is compatible with products on the left and on the right. Let $(M,\cdot,1,<)$ be an ordered non-associative monoid $(M,\cdot,1,<)$. Recall from Example 2.10 that we have summabiliy module $k(M,\mathcal{N})$ for the bornology $\mathcal{N}$ on $M$ of all Noetherian subsets of $(M,<)$. The module $k(\!(M)\!):=k(M;\mathcal{N})$ is a unital summability algebra for the Cauchy product

$$\forall m \in M, (a \cdot b)(m) := \sum_{uv=m} a(u)\,b(v) \qquad (3.1)$$

for all $a,b \in k(\!(M)\!)$. Its unit is $\mathbb{1}_1$. See [7, Section 1.4] for more details (the proofs there are given in the case of associative monoids over a field, but they carry out without change for non-associative monoids over $k$). In the sequel we write $\mathbb{A} = k(\!(M)\!)$.

**Lemma 3.10.** *Let $A$ be a summability algebra. Let $\varphi: M \longrightarrow A$ be a map such that $(\varphi(n))_{n \in N}$ is summable whenever $N \subseteq M$ is Noetherian. Let $\tilde{\phi}: \mathbb{A} \longrightarrow A$ be its unique strongly linear extension as per Lemma 2.12. Then $\tilde{\phi}$ is an algebra morphism if and only if $\phi(m \cdot n) = \phi(m) \cdot \phi(n)$ for all $m,n \in M$.*

**Proof.** This follows from the definition (3.1). □

We define an ordering $\prec$ on $\mathbb{A}$ as follows: for $a,b \in \mathbb{A}$,

$$a \prec b \quad \text{if} \quad (b \neq 0 \text{ and } \forall m \in \operatorname{supp} b, \exists n \in \operatorname{supp} a, n < m). \qquad (3.2)$$

The map $m \mapsto \mathbb{1}_m : (M,\cdot,1,<) \longrightarrow (\mathbb{A},\cdot,1,\prec)$ is an embedding of ordered non-associative monoids. Given two linear maps $\phi,\psi: \mathbb{A} \longrightarrow \mathbb{A}$, we write $\phi \prec \psi$ if $\phi(a) \prec \psi(a)$ for all $a \in \mathbb{A} \setminus \{0\}$. We say that $\phi$ is *contracting* if $\phi \prec \mathrm{id}$, and that it is *tangent to the identity* if $\varphi - \mathrm{id} \prec \mathrm{id}$. We write $\mathrm{Lin}^+_\prec(\mathbb{A})$ for the set of contracting strongly linear maps $\mathbb{A} \longrightarrow \mathbb{A}$, so $\mathrm{id}_\mathbb{A} + \mathrm{Lin}^+_\prec(\mathbb{A})$ is the set of strongly linear maps that are tangent to the identity. Both $\mathrm{Lin}^+_\prec(\mathbb{A})$ and $k\,\mathrm{id}_\mathbb{A} + \mathrm{Lin}^+_\prec(\mathbb{A})$ are summability subalgebras of $\mathrm{Lin}^+(\mathbb{A})$ (see [7, Corollary 3.6]).

## 3.2 Formal power series

Fix a set $I$. The Magnus algebras of formal power series and formal non-associative power series over $I$ respectively are denoted $k\langle\!\langle I \rangle\!\rangle$ and $k\langle\!\langle I \rangle\!\rangle^{\mathrm{na}}$. They can be seen as algebras of Noetherian series $k\langle\!\langle I \rangle\!\rangle = k(\!(I^\star)\!)$ and $k\langle\!\langle I \rangle\!\rangle^{\mathrm{na}} = k(\!(I^{(\star)})\!)$ for certain orderings $<$ and $<_{\mathrm{na}}$ on $I^\star$ and $I^{(\star)}$ such that $(I^\star,<)$ and $(I^{(\star)},<_{\mathrm{na}})$ are Noetherian. Indeed fix a well-ordering $<$ on $I$. Then the lexicographic ordering on $I^\star$ is a well-ordering for which $I^\star$ is an ordered monoid. This also holds for $I^{(\star)}$ seen as a submonoid of $(I \sqcup \{(,)\})^\star$, extending the ordering $<$ to $I \sqcup \{(,)\}$ by putting $($ and $)$ at the bottom.

If $m \in \mathbb{N}$, then we write $k\langle\!\langle m \rangle\!\rangle = k\langle\!\langle \{0,\ldots,m-1\} \rangle\!\rangle$ and $k\langle\!\langle m \rangle\!\rangle^{\mathrm{na}} = k\langle\!\langle \{0,\ldots,m-1\} \rangle\!\rangle^{\mathrm{na}}$. We often write elements of $k\langle\!\langle I \rangle\!\rangle^{\mathrm{na}}$ as formal power series in non-commuting, non-associative variables $X_i, i \in I$. Indeed, setting $X_u = \mathbb{1}_u$ for all $u \in I^{(\star)}$, we have

$$X_u \cdot X_v = X_{(uv)} \qquad (3.3)$$

for all $u,v \in I^{(\star)}$. Moreover, each $P \in k\langle\!\langle I \rangle\!\rangle^{\mathrm{na}}$ is the sum

$$P = \sum_{u \in I^{(\star)}} P(u)\,X_u. \qquad (3.4)$$

The algebra $k\langle\!\langle I\rangle\!\rangle^{\mathrm{na}}$ will play the role of free object in our categories of algebras, once we define a notion of evaluation morphism from $k\langle\!\langle I\rangle\!\rangle^{\mathrm{na}}$ to certain summability algebras. We usually write $X_I$ for the (summable) family $I \longrightarrow k\langle\!\langle I\rangle\!\rangle$ that sends each $i \in I$ to $X_i$.

**Remark 3.11.** If $J \subseteq I$ is a subset, then we write $\pi_J^{\mathrm{na}}$ (resp. $\pi_J$) for the strongly linear projection $k\langle\!\langle I\rangle\!\rangle^{\mathrm{na}} \longrightarrow k\langle\!\langle J\rangle\!\rangle^{\mathrm{na}}$ (resp. $k\langle\!\langle I\rangle\!\rangle \longrightarrow k\langle\!\langle J\rangle\!\rangle$) that sends a $P$ to its restriction to $J^{(\star)}$ (resp. $J^\star$): these are strongly linear morphisms of algebras.

**Remark 3.12.** We have a valuation map $\mathrm{val} : k\langle\!\langle I\rangle\!\rangle \longrightarrow \mathbb{N} \cup \{+\infty\}$ with $\mathrm{val}(0) = +\infty$, and

$$\mathrm{val}(P) := \min \{n : \mathrm{supp}\, P \cap I^n \neq \varnothing\}.$$

for all $P \neq 0$. This induces the topology of topological vector space on $k\langle\!\langle I\rangle\!\rangle$ called the valuation topology that has the collection of sets $\{P \in k\langle\!\langle I\rangle\!\rangle : \mathrm{val}(P) \geqslant n\}$ as a basis of open neighborhoods of 0. This topology is coarser than the product topology (see Remark 2.14). The same applies to $k\langle\!\langle I\rangle\!\rangle^{\mathrm{na}}$.

## 3.3 Lyndon words

Let $I$ be a non-empty set. Given a well-ordering on $I$ with lexicographic extension $<_{\mathrm{lex}}$ on $I^\star$ we write $<_{\mathrm{gr}}$ for the ordering on $I^\star$ given by $u < v$ if $\mathrm{rank}(u) < \mathrm{rank}(v)$ or $\mathrm{rank}(u) = \mathrm{rank}(v)$ and $u <_{\mathrm{lex}} v$ (i.e. the pullback by the map $u \mapsto (\mathrm{rank}(u), u) : I^\star \longrightarrow \mathbb{N} \times I^\star$ of lexicographic ordering on $\mathbb{N} \times I^\star$). It is easy to see that $I^\star$ is a well-ordered ordered associative monoid for this ordering (under concatenation). We call $<_{\mathrm{gr}}$ the graded extension of $<$. A Lyndon word on $I$ is a non-trivial associative word $w$ on $I$ which is $<_{\mathrm{lex}}$-smaller than all its right segments, i.e. with $w <_{\mathrm{lex}} v$ whenever $w = uv$ for non-trivial words $u, v \in I^\star$. We write $\mathrm{Lydon}(I)$ for the set of all Lyndon words on $I$.

For each Lyndon word $w$, we define a non-associative word $w^{()}$ by induction on the rank of $w$ as follows. If $w$ is a letter then $w^{()} := w$. Otherwise $w$ can be decomposed [41, (5.1.1)] uniquely as a concatenation $w = uv$ where $u, v$ are Lyndon words and $v$ is longest. We then set $w^{()} = (u^{()} v^{()})$. We have the classical result:

**Proposition 3.13.** *Let $f : I \longrightarrow k\langle\!\langle I\rangle\!\rangle$ be the family with $f(i) = X_i$ for all $i \in I$. Let $w$ be a Lyndon word. Then the support of $w^{()}[f] - X_w$ lies strictly above $w$ for the ordering $<_{\mathrm{gr}}$.*

**Proof.** The classical statement [41, Theorem 5.1] gives the result for $<_{\mathrm{lex}}$, not $<_{\mathrm{gr}}$. However it is easy to see by induction on the rank of $w$ that all words in the support of $w^{()}[f]$ have the same rank, so $<_{\mathrm{lex}}$ and $<_{\mathrm{gr}}$ coincide on this set. □

**Proposition 3.14.** [41, Theorem 4.9] *Suppose that $S$ is a finite set and let $\mathrm{Lie}(S)$ be the free Lie algebra on $S$. Then $\{u[X] : u \in \mathrm{Lyndon}(S)\}$ is a basis of $k$-module of $\mathrm{Lie}(S)$.*

## 3.4 Some limits and colimits

We consider the categories $\mathbf{\Sigma Alg}$, $\mathbf{\Sigma Ass}$ and $\mathbf{\Sigma Lie}$ of summability algebras, summability associative algebras and summability Lie algebras respectively, with strongly linear algebra homomorphisms as arrows. We will justify that these categories are complete and cocomplete.

Although $\mathbf{\Sigma Alg}$ is not a variety of algebras in the classical sense, the inclusions $\mathbf{\Sigma Ass} \subset \mathbf{\Sigma Alg}$ and $\mathbf{\Sigma Lie} \subset \mathbf{\Sigma Alg}$ have the variants of the classical left adjoints described hereafter. An equational theory is a first-order theory in $\mathfrak{P}_{k\text{-alg}}$ containing only sentences $\forall \bar{x}(t(\bar{x}) = 0)$ where $t(\bar{x})$ is a term.

**Definition 3.15.** *Let $T$ be an equational theory. We say that the full subcagegory $\Sigma T$ of summability algebras $(A, \Sigma)$ with $A \vDash T$ is a **variety of summability algebras**.*

**Proposition 3.16.** *Let $\Sigma T$ be a variety of summability algebras. Then the inclusion $\Sigma T \hookrightarrow \Sigma\mathbf{Alg}$ has a left adjoint.*

**Proof.** Given an object $(A, \Sigma)$ in $\Sigma T$, let $I_T(A)$ denote the smallest summability ideal of $A$ containing all $t(\bar{a})$ where $\forall \bar{x}(t(\bar{x}) = 0) \in T$ and $\bar{a} \subseteq A$. Then $T(A) := A/I_T(A)$ is a summability algebra with $T(A) \vDash T$. If $f : A \longrightarrow B$ is an arrow in $\Sigma\mathbf{Alg}$, then $f(I_T(A)) \subseteq I_T(B)$ as $f^{-1}(I_T(B))$ is a summability ideal of $A$ that contains all $t(\bar{a})$ where $\forall \bar{x}(t(\bar{x}) = 0) \in T$ and $\bar{a} \subseteq A$. So we have a strongly linear morphism $T(f) : T(A) \longrightarrow T(B)$ and one can see that this defines a functor $\Sigma\mathbf{Alg} \longrightarrow \Sigma T$. That it is left ajoint to the inclusion follows as in the classical case. $\square$

**Lemma 3.17.** *Let $T$ be an equational theory such that any $t(\bar{x})$ with $\forall \bar{x}(t(\bar{x}) = 0) \in T$ is a multilinear term. Let $(A, \Sigma)$ be a summability algebra and let $B \subseteq A$ be a subset with $t(\bar{b}) = 0$ whenever $\forall \bar{x}(t(\bar{x}) = 0) \in T$ and $\bar{b} \subseteq B$. The smallest summability subalgebra of $A$ containing $B$ is a model of $T$.*

**Proof.** Let $S \subseteq A$ be the smallest summability algebra of $A$ containing $B$. We claim that $S$ can be obtained as the union $S = \bigcup_{\alpha \in \mathbf{On}} B_\alpha$ where $B_0$ is the subalgebra generated by $B$, and for all $\alpha > 0$,

$$B_\alpha := \left\{ \Sigma_I a : a \in \operatorname{dom} \Sigma_I \cap \left( \bigcup_{\beta < \alpha} B_\beta \right)^I \right\}.$$

Indeed, this union is clearly contained in $S$. Now if $a : I \longrightarrow \bigcup_{\alpha \in \mathbf{On}} B_\alpha$ is summable in $A$, then as $I$ is a set, there must be an ordinal $\alpha \in \mathbf{On}$ with $a(I) \subseteq B_\alpha$, whence $\Sigma_I a \in B_{\alpha+1}$. So $\bigcup_{\alpha \in \mathbf{On}} B_\alpha$ is a summability submodule of $A$ containing $B$. It is easy to see by induction on $\alpha$ and by **SA** that each $B_\alpha$ is a subalgebra of $A$. Therefore $S = \bigcup_{\alpha \in \mathbf{On}} B_\alpha$.

Let $\varphi = \forall \bar{x}(t(\bar{x})) \in T$ where $t(\bar{x})$ is a term of arity $n \in \mathbb{N}$. Recall that $t(\bar{x})$ acts a multilinear function in all algebras. We claim that the map $\bar{a} \mapsto t(\bar{a}) : A^n \longrightarrow A$ is strongly multilinear, meaning that for all summable families $a_1 : I_1 \longrightarrow A, \ldots, a_n : I_n \longrightarrow A$, the family $t(a_1, \ldots, a_n) : I_1 \times \cdots \times I_n \longrightarrow A$ given by $t(a_1, \ldots, a_n)(i_1, \ldots, i_n) := t(a_1(i_1), \ldots, a_n(i_n))$ for all $(i_1, \ldots, i_n) \in I_1 \times \cdots \times I_n$ is summable, with sum $\Sigma_{I_1 \times \cdots \times I_n} t(a_1, \ldots, a_n) = t(\Sigma_{I_1} a_1, \ldots, \Sigma_{I_n} a_n)$. Indeed this follows, by induction on the total degree of $t(\bar{x})$, from **SA**. We next claim that for all ordinals $\alpha$, we have $t(B_\alpha^n) = \{0\}$. Indeed, this is true for $\alpha = 0$ as $B \vDash T$. Let $\alpha > 0$ such that $t(B_\beta^n) = \{0\}$ for all $\beta < \alpha$ and write $X_\alpha := \bigcup_{\beta < \alpha} B_\beta$. Note that $t(X_\alpha^n) = \{0\}$. Let $(b_1, \ldots, b_n) \in B_\alpha$. So there are summable families $a_1 : I_1 \longrightarrow X_\alpha, \ldots, a_n : I_n \longrightarrow X_\alpha$ with $(\Sigma_{I_1} a_1, \ldots, \Sigma_{I_n} a_n) = (b_1, \ldots, b_n)$. Now $t(b_1, \ldots, b_n) = \Sigma_{I_1 \times \cdots \times I_n} t(a_1, \ldots, a_n) = \Sigma_{I_1 \times \cdots \times I_n} 0 = 0$ by our previous claim. We conclude by induction that $t(B_\gamma^n) = \{0\}$ for all $\gamma \in \mathbf{On}$, so $t(S^n) = \{0\}$. This shows that $S \vDash T$. $\square$

The equalizer, in the category of algebras, of two arrows in $\Sigma\mathbf{Alg}$ is a summability subalgebra, and hence has a natural structure of summability algebra. This gives the equalizer in $\Sigma\mathbf{Alg}$, hence in $\Sigma\mathbf{Ass}$ and $\Sigma\mathbf{Lie}$. The pointwise product of a family of summability algebras as described in Example 3.5, in the special case of the filter of all subsets, is the product in $\Sigma\mathbf{Alg}$. This construction restricts to the product in $\Sigma\mathbf{Ass}$ and $\Sigma\mathbf{Lie}$. This proves the completeness of these categories. Let us next turn to colimits.

Let $A, B$ be summability algebras and let $f, g : A \longrightarrow B$ be parallel strongly linear morphisms. Let $I(f, g)$ denote the smallest summability ideal of $B$ containing $\{f(a) - g(a) : a \in A\}$ and let $\operatorname{Coeq}(f, g) := B/I(f, g)$ be the quotient summability algebra.

**Lemma 3.18.** *The summability algebra* $\mathrm{Coeq}(f,g)$ *with the quotient map* $\mathrm{coeq}(f,g)$ *is the coequalizer of* $(f,g)$ *in* $\boldsymbol{\Sigma}\mathbf{Alg}$.

**Proof.** We have $\mathrm{coeq}(f,g)\circ f = \mathrm{coeq}(f,g)\circ g$ by definition. Let $C\in\boldsymbol{\Sigma}\mathbf{Alg}$ and $h\colon B\longrightarrow C$ an arrow with $h\circ f = h\circ g$. Then $\{f(a)-g(a):a\in A\}\subseteq\mathrm{Ker}(h)$, so $I(f,g)\subseteq\mathrm{Ker}(h)$ as $\mathrm{Ker}(h)$ is a summability ideal. Thus we have a unique morphism $\varphi\colon\mathrm{Coeq}(f,g)\longrightarrow C$ with $\varphi\circ\mathrm{coeq}(f,g)=h$, given by $\varphi(b+I(f,g))=h(b)$. It suffices in order to conclude to justify that $\varphi$ is strongly linear. If $(B_i)_{i\in I}$ is summable in $\mathrm{Coeq}(f,g)$, then by definition of the quotient summability structure, there is a summable $(b_i)_{i\in I}$ in $B$ with $b_i\in B_i$ for all $i\in I$. The family $(\varphi(B_i))_{i\in I}=(h(b_i))_{i\in I}$ is summable by strong linearity of $h$, and $\varphi(\sum_{i\in I}B_i)=\varphi(\sum_{i\in I}b_i+I(f,g))=h(\sum_{i\in I}b_i)=\sum_{i\in I}h(b_i)=\sum_{i\in I}\varphi(B_i)$. So $\varphi$ is strongly linear. □

Let $X$ be a non-empty set, and let $(A_x,\Sigma_x)_{x\in X}$ be a family of summability algebras. Consider the disjoint union $A:=\bigsqcup_{x\in X}A_x\times\{x\}$ as an alphabet for the free non-associative monoid $A^{(\star)}$. We consider the following bornology on $A^{(\star)}$. Let $\mathcal{F}_0$ be set of sets $S\times\{x\}$ where $S\subseteq A_x$ and $(s)_{s\in S}$ is summable in $(A_x,\Sigma_x)$. Let $\mathcal{F}$ be the smallest bornology on $A^{(\star)}$ that contains $\mathcal{F}_0$ and that is closed under the operation $(U,V)\mapsto(UV):=\{(uv):u\in U\wedge v\in V\}$. So $\mathcal{F}$ is the increasing union $\mathcal{F}=\bigcup_{n\in\mathbb{N}}\mathcal{F}_n$ where $\mathcal{F}_{n+1}=\mathcal{F}_n\cup\{(UV):U,V\in\mathcal{F}_n\}$. Note that no element of $\mathcal{F}$ contains the empty word on $A$.

We see $k(A^{(\star)};\mathcal{F})$ as naturally included into $k\langle\!\langle A\rangle\!\rangle^{\mathrm{na}}$ and we note by definition of $\mathcal{F}$ that it is a subalgebra of $k\langle\!\langle A\rangle\!\rangle^{\mathrm{na}}$. By definition of $\mathcal{F}$ and the product, the operation turns $(k(A^{(\star)};\mathcal{F}),\Sigma_\mathcal{F})$ into a summability algebra. Lastly, let $I$ denote the smallest summability ideal of $(k(A^{(\star)};\mathcal{F}),\Sigma_\mathcal{F})$ that contains all differences $X_{(a,x)}-\Sigma_\mathcal{F}(c_s X_{(s,x)})_{s\in S}$ where $a\in A_x$, $S\subseteq A_x$, $(\lambda_s)_{s\in S}\in k^S$ and $a=\Sigma_x(\lambda_s s)_{s\in S}$ in $(A_x,\Sigma_x)$, and all differences $X_{(a,x)}-X_{(b,x)}X_{(c,x)}$ where $a=bc$ in $A_x$. We consider the quotient summability algebra $k(A^{(\star)};\mathcal{F})/I$.

**Proposition 3.19.** *Let be a family of summability Lie algebras. The summability algebra* $k(A^{(\star)};\mathcal{F})/I$ *is the coproduct of* $(A_x,\Sigma_x)_{x\in X}$ *in* $\boldsymbol{\Sigma}\mathbf{Alg}$.

**Proof.** For $x\in X$, we have a strongly linear morphism $\iota_x\colon A_x\longrightarrow k(A^{(\star)};\mathcal{F})/I$ given by $a\mapsto\mathbb{1}_{(a,x)}+I$. Let $(B,\Sigma)$ be a summability algebra and let us be given, for each $x\in X$, a strongly linear morphism $\varphi_x\colon A_x\longrightarrow B$. We first define a strongly linear morphism $\psi\colon k(A^{(\star)};\mathcal{F})\longrightarrow B$. The support of a $P\in k(A^{(\star)};\mathcal{F})$ lies in $\mathcal{F}_n$ for some $n\in\mathbb{N}$, hence there are $m\geqslant n$, $x_1,\ldots,x_m\in X$ and families $(s_1)_{s_1\in S_1},\ldots,(s_m)_{s_m\in S_m}$ where each $S_i\subseteq A_{x_i}$ and each $(s_i)_{s_i\in S_i}$ is summable in $(A_{x_i},\Sigma_{x_i})$, such that $(X_u)_{u\in\mathrm{supp}\,P}$ is a subfamily of $(X_u)_{u\in D}$ where $D$ is a subset of the set of products of elements in $S_1\times\{x_1\}\cup\ldots\cup S_m\times\{x_m\}$ of depth $n$. Consider the family $f=(\varphi_x(s))_{(s,x)\in\bigsqcup_{i=1}^n S_i\times\{x_i\}}$, extended with zeros outside of $\bigsqcup_{i=1}^n S_i\times\{x_i\}$. Note since $\mathrm{supp}\,P$ does not contain the empty word that $u(f)$ is well defined in the magma $(B,\cdot)$ for all $u\in\mathrm{supp}\,P$. By strong linearity of each $\varphi_x$, the families $(\varphi_x(s_i))_{s_i\in S_i}$ are summable in $B$. Applying **SA** in $B$ consecutively $n$ times, we see that $(\mathrm{eu}(f))_{u\in\mathrm{supp}\,P}$ is summable in $B$. With **SM5**, we may define $\psi(P):=\sum_{u\in\mathrm{supp}\,P}P(u)\,u(f)$, and we see with Lemmas 2.11 and 3.10 that this is a strongly linear morphism of algebras $k(A^{(\star)};\mathcal{F})\longrightarrow B$. Since each $\varphi_x$ for $x\in X$ is a strongly linear morphism of algebras, the kernel of $\psi$ contains $I$. So we obtain a strongly linear morphism $\varphi\colon k(A^{(\star)};\mathcal{F})/I\longrightarrow B$, which by definition is unique to satisfy $\varphi\circ\iota_x=\varphi_x$ for all $x\in X$. This concludes the proof. □

**Remark 3.20.** For $(A,\Sigma)\in\boldsymbol{\Sigma}\mathbf{Alg}$, seeing $A$ as the coproduct of a one-point family, we see that $A$ is isomorphic to a quotient by a summability ideal of a "bornological algebra" $(k\langle\!\langle A\rangle\!\rangle;\mathcal{F})$.

To obtain coproducts for associative algebras and Lie algebras, we apply the corresponding left ajoints.

**Proposition 3.21.** *The categories of summability algebras, summability associative algebras and summability Lie algebras are complete and cocomplete.*

## 3.5 Algebra evaluations

**Definition 3.22.** *A summability algebra $(A, \Sigma)$ **has evaluations** (or is **with evaluations**) if for all summable families $a: I \longrightarrow A$, there is a strongly linear morphism $\mathrm{ev}_a : k\langle\!\langle I \rangle\!\rangle_0^{\mathrm{na}} \longrightarrow A$ with $\mathrm{ev}_a(X_i) = a(i)$ for all $i \in I$.*

Note in view of (3.4) that evaluation morphisms $\mathrm{ev}_a$ when they exist are unique. If $m \in \mathbb{N}$ and $a = (a_i)_{i \in \{0,\ldots,m\}} \in A^m$ is a family, we will sometimes write $\mathrm{ev}_{a_0,\ldots,a_m}$ or $\mathrm{ev}_{(a_0,\ldots,a_m)}$ for $\mathrm{ev}_a$.

**Lemma 3.23.** *Let $A$ be a summability algebra $A$ and let $a: I \longrightarrow A$ be a family. Then $\mathrm{ev}_a : k\langle\!\langle I \rangle\!\rangle_0^{\mathrm{na}} \longrightarrow A$ exists if and only if the family $(u(a))_{u \in I^{(\star)} \setminus \{\varnothing\}}$ is summable.*

**Proof.** If $(u(a))_{u \in I^{(\star)} \setminus \{\varnothing\}}$ is summable, then by **SM5**, for each $P \in k\langle\!\langle I \rangle\!\rangle_0^{\mathrm{na}}$, the family $(P(u) u(a))_{u \in I^{(\star)} \setminus \{\varnothing\}}$ is summable, and by Lemma 2.11, the map $P \mapsto \sum_{u \in I^{(\star)} \setminus \{\varnothing\}} P(u) u(a)$ is strongly linear. By Lemma 3.10, it is a morphism of algebras. Conversely, if $\mathrm{ev}_a : k\langle\!\langle I \rangle\!\rangle_0^{\mathrm{na}} \longrightarrow A$ exists, then as $(X_u)_{u \in I^{(\star)} \setminus \{\varnothing\}}$ is summable in $k\langle\!\langle I \rangle\!\rangle_0^{\mathrm{na}}$, so is $(\mathrm{ev}_a(X_u))_{u \in I^{(\star)} \setminus \{\varnothing\}} = (u(a))_{u \in I^{(\star)} \setminus \{\varnothing\}}$. □

**Corollary 3.24.** *Let $(A, \Sigma)$ be a non-unital summability algebra with evaluations. Then all summability subalgebras of $(A, \Sigma)$ have evaluations.*

**Proposition 3.25.** *Let $A$ be an algebra and let $\Sigma^{\min}$ be the minimal summability structure on $A$ (see Example 2.6). Then $(A, \Sigma^{\min})$ has evaluations if and only if $A$ is locally nilpotent.*

**Proof.** Recall that a family $f: I \longrightarrow A$ is summable in $(A, \Sigma^{\min})$ if and only if $f$ is finitely supported. Let $S \subseteq A$ be finite and consider the summable family $\mathrm{id}_S : S \longrightarrow A$. As $S$ is finite, the family $(u(\mathrm{id}_S))_{u \in S^{(\star)} \setminus \{\varnothing\}}$ is finitely supported if and only if for all $u \in S^{(\star)}$ of sufficiently large rank, we have $u(\mathrm{id}_S) = 0$. Thus $(u(\mathrm{id}_S))_{u \in S^{(\star)} \setminus \{\varnothing\}}$ is finitely supported if and only if the subalgebra of $A$ generated by $S$ is nilpotent. We deduce with Lemma 3.23 that $(A, \Sigma^{\min})$ has evaluations if and only if $A$ is locally nilpotent. □

**Example 3.26.** Consider a lower central complete algebra $A$ with the summability structure of Example 3.6, i.e. the inverse limit structure of $A \simeq \varprojlim (A/A_n)_{n>0}$ for the minimal summability structures $\Sigma^n$ on each quotient $A/A_n$. As each quotient $A/A_n$ is nilpotent, $(A/A_n, \Sigma^n)$ has evaluations by Proposition 3.25. It is easy to see that the inverse limit is a summability subalgebra of the corresponding product, that has evaluations (see also Section 3.8). So $A$ has evaluations.

**Proposition 3.27.** *Let $(A, \Sigma)$ be a summability algebra and let $\mathfrak{q} \subseteq A$ be a summability ideal. Then the quotient summability algebra $(A/\mathfrak{q}, \Sigma^{/\mathfrak{q}})$ has evaluations.*

**Proof.** Let $a: I \to A/\mathfrak{q}$ be summable. By definition of the summability structure $\Sigma^{/\mathfrak{q}}$ (see Example 2.8), there is a summable family $b: I \longrightarrow A$ with $a(i) = b(i) + \mathfrak{q}$ for all $i \in I$. Composing the evaluation map $\mathrm{ev}_b : k\langle\!\langle I \rangle\!\rangle^{\mathrm{na}} \longrightarrow A$ with the strongly linear projection $A \longrightarrow A/\mathfrak{q}$ yields the evaluation map $\mathrm{ev}_a : k\langle\!\langle I \rangle\!\rangle^{\mathrm{na}} \longrightarrow A/\mathfrak{q}$. □

A unital $A$ summability algebra cannot have evaluations. Indeed the one-point family $\{\bullet\} \to \{1\}$ is summable in $A$, but $(1^n)_{n \in \mathbb{N}} = (\mathrm{ev}_{\{\bullet\} \to \{1\}}(X_\bullet^n))_{n \in \mathbb{N}}$ is not summable (see Remark 2.2). This is why we adopt a distinct definition for unital algebras:

**Definition 3.28.** [7, Definition 1.39/Theorem 1.43] *Let $A = k \oplus \mathfrak{m}$ be a unital summability algebra where $\mathfrak{m}$ is a summability ideal. We say that $A$ has **unital evaluations** if $\mathfrak{m}$ has evaluations, i.e. if for all summable families $a : I \longrightarrow \mathfrak{m}$, there is a strongly linear morphism $\mathrm{ev}_a : k \langle\!\langle I \rangle\!\rangle^{\mathrm{na}} \longrightarrow A$ with $\mathrm{ev}_a(X_i) = a(i)$ for all $i \in I$.*

As in Lemma 3.23, we have:

**Lemma 3.29.** *Let $A = k \oplus \mathfrak{m}$ be a unital summability algebra where $\mathfrak{m} \subseteq A$ is a summability ideal. Let $a : I \longrightarrow A$ be a family. Then $\mathrm{ev}_a : k \langle\!\langle I \rangle\!\rangle^{\mathrm{na}} \longrightarrow A$ exists if and only if the family $(u(a))_{u \in I^{(\star)}}$ is summable.*

The most important example of summability algebra with evaluations is the associative algebra of contracting strongly linear maps on a module of Noetherian series:

**Proposition 3.30.** [7, Theorem 3.11] *Let $(X, <)$ be a non-empty ordered set and let $\mathcal{N}$ be the bornology of Noetherian subsets of $(X, <)$. Then $\mathrm{Lin}_\prec^+(k(X; \mathcal{N}))$ has evaluations.*

Let $(M, \cdot, 1, <)$ be an ordered non-associative monoid and let $\mathcal{N}$ be the bornology of Noetherian subsets of $(M, <)$. We set $k((M^{>1})) = \{f \in k((M)) : \mathrm{supp}\, f > 1\}$, seen as a summability subalgebra of $k((M))$.

**Corollary 3.31.** *The summability algebra $k((M^{>1}))$ has evaluations.*

**Proof.** We see $k((M^{>1}))$ as embedded into $\mathrm{Lin}_\prec^+(k((M)))$ via

$$\Phi : k((M^{>1})) \longrightarrow \mathrm{Lin}_\prec^+(k((M)))$$
$$f \longmapsto (g \mapsto f \cdot g).$$

A family in $k((M^{>1}))$ is summable if and only if its image under $\Phi$ is summable in $\mathrm{Lin}^+(k((M)))$. In view of Corollary 3.24, this yields the result. $\square$

**Corollary 3.32.** *Let $I$ be a set. Then $k \langle\!\langle I \rangle\!\rangle_0^{\mathrm{na}}$ and $k \langle\!\langle I \rangle\!\rangle_0$, have evaluations.*

**Proof.** Recall that $k \langle\!\langle I \rangle\!\rangle^{\mathrm{na}}$ can be realized as $k((I^{(\star)}))$ for a linear ordering on $I^{(\star)}$ for which $\varnothing$ is the minimum. Thus $k((I^{(\star) > \varnothing})) = k \langle\!\langle I \rangle\!\rangle_0^{\mathrm{na}}$. The same applies to $k \langle\!\langle I \rangle\!\rangle$. $\square$

So we can now compose evaluations by evaluating into algebras of formal power series, and then into arbitrary algebras with evaluations. The following simple but crucial fact can be interpreted as an associativity propety of evaluations.

**Proposition 3.33.** *Let $A$ be a non-unital summability algebra with evaluations, let $a : I \longrightarrow A$ be summable and let $P : J \longrightarrow k \langle\!\langle I \rangle\!\rangle_0^{\mathrm{na}}$ be summable. Then for all $Q \in k \langle\!\langle J \rangle\!\rangle_0^{\mathrm{na}}$, we have*

$$\mathrm{ev}_a(\mathrm{ev}_P(Q)) = \mathrm{ev}_{(\mathrm{ev}_a(P(j)))_{j \in J}}(Q).$$

**Proof.** As evaluation maps are strongly linear morphisms of algebras, it suffices to prove the result for $Q = X_j$ for each $j \in J$. Now for $j \in J$, we have $\mathrm{ev}_a(\mathrm{ev}_P(Q)) = \mathrm{ev}_a(P(j)) = \mathrm{ev}_{(\mathrm{ev}_a(P(j)))_{j \in J}}(X_j)$. $\square$

For the next proposition, we use the notations of Section 3.4 regarding varieties of summability algebras.

**Proposition 3.34.** *Let $\Sigma T$ be a variety of summability algebras. Then an object $A$ of $\Sigma T$ has evaluations if and only if for all summable families $a: I \longrightarrow A$, there is a strongly linear morphism $\mathrm{ev}_a^T: T(k\langle\!\langle I\rangle\!\rangle^{\mathrm{na}}) \longrightarrow A$, which is unique, with $\mathrm{ev}_a^T(X_i + I_T(A)) = a(i)$ for all $i \in I$.*

**Proof.** Let $a: I \longrightarrow A$ be summable. The unicity of evaluation maps comes from the fact that each element $f = P + I_T(k\langle\!\langle I\rangle\!\rangle^{\mathrm{na}})$ of $T(k\langle\!\langle I\rangle\!\rangle^{\mathrm{na}})$ can be writen as

$$f = \sum_{u \in I^{(\star)}} P(u)\,(X_u + I_T(k\langle\!\langle I\rangle\!\rangle^{\mathrm{na}})).$$

Let $\iota: I \longrightarrow T(k\langle\!\langle I\rangle\!\rangle^{\mathrm{na}})$ with $\iota(i) = X_i + I_T(k\langle\!\langle I\rangle\!\rangle^{\mathrm{na}})$ for all $i \in I$. Note that $\iota$ is summable. As $T(k\langle\!\langle I\rangle\!\rangle^{\mathrm{na}})$ has evaluations by Proposition 3.27, we obtain an evaluation map $\mathrm{ev}_\iota: k\langle\!\langle I\rangle\!\rangle^{\mathrm{na}} \longrightarrow T(k\langle\!\langle I\rangle\!\rangle^{\mathrm{na}})$. Suppose that $A$ has evaluations and let $\mathrm{ev}_a: k\langle\!\langle I\rangle\!\rangle^{\mathrm{na}} \longrightarrow A$ be the evaluation map. Then as $A \vDash T$, we have $\mathrm{ev}_a = 0$ on $I_T(k\langle\!\langle I\rangle\!\rangle^{\mathrm{na}})$, so $\mathrm{ev}_a$ induces the desired evaluation map $T(k\langle\!\langle I\rangle\!\rangle^{\mathrm{na}}) \longrightarrow A$. Suppose conversely that $\mathrm{ev}_a^T: T(k\langle\!\langle I\rangle\!\rangle^{\mathrm{na}}) \longrightarrow A$ exists. Then $\mathrm{ev}_a^T \circ \mathrm{ev}_\iota$ satisfies the conditions of Definition 3.22, so $A$ has evaluations. $\square$

## 3.6 Exponential and logarithm

In the associative case, the associative formal series

$$\exp(X_0) := \sum_{n \in \mathbb{N}} \frac{1}{n!} X_0^n \in 1 + k\langle\!\langle 1\rangle\!\rangle_0,$$

$$\log(1 + X_0) := \sum_{n \in \mathbb{N}} \frac{(-1)^n}{n+1} X_0^{n+1} \in k\langle\!\langle 1\rangle\!\rangle_0$$

will induce functions on summability algebras with evaluations. We first give the two following strongly linear isomorphisms for all sets $I$

$$\mathrm{Ass}(k\langle\!\langle I\rangle\!\rangle^{\mathrm{na}}) \simeq k\langle\!\langle I\rangle\!\rangle,$$
$$\mathrm{Ab}(k\langle\!\langle 2\rangle\!\rangle) \simeq k[[X_0, X_1]] \quad \text{(ring of formal series in two commuting variables)}.$$

The strong linearity of the canonical isomorphisms follows from Lemma 2.11. With

**Lemma 3.35.** *An associative summability algebra $A$ has evaluations if and only if for all summable families $a: I \longrightarrow A$, there is a strongly linear morphism $\mathrm{ev}_a^{\mathrm{Ass}}: k\langle\!\langle I\rangle\!\rangle_0 \longrightarrow A$ with $\mathrm{ev}_a^{\mathrm{Ass}}(X_i) = a(i)$.*

**Proof.** This follows from Proposition 3.34 and the isomorphism $\mathrm{Ass}(k\langle\!\langle I\rangle\!\rangle^{\mathrm{na}}) \simeq k\langle\!\langle I\rangle\!\rangle$. $\square$

**Lemma 3.36.** *A unital associative algebra $A = k \oplus \mathfrak{m}$ has unital evaluations if and only if for all summable families $a: I \longrightarrow \mathfrak{m}$, the family $(a(i_1) \cdots a(i_n))_{(i_1,\ldots,i_n) \in I^\star}$ is summable in $\mathfrak{m}$.*

**Proof.** This follows from combining Lemma 3.35 and [7, Theorem 1.43]. $\square$

We now fix a unital associative algebra $A = k \oplus \mathfrak{m}$ with unital evaluations.

**Proposition 3.37.** [7, Proposition 1.45] *The algebra $A$ is local with maximal ideal $\mathfrak{m}$. In particular, the subset $1 + \mathfrak{m} \subseteq A^\times$ is a group.*

We define, for all $\varepsilon \in \mathfrak{m}$,

$$\begin{aligned} \exp(\varepsilon) &= \mathrm{ev}_\varepsilon(\exp(X_0)) \in 1+\mathfrak{m} \\ \log(1+\varepsilon) &= \mathrm{ev}_\varepsilon(\log(1+X_0)) \in \mathfrak{m}. \end{aligned}$$

**Proposition 3.38.** *The functions* $\exp : \mathfrak{m} \longrightarrow 1+\mathfrak{m}$ *and* $\log : 1+\mathfrak{m} \longrightarrow \mathfrak{m}$ *are bijections and* $\log$ *is the functional inverse of* $\exp$. *Moreover, we have* $\exp(\varepsilon + \delta) = \exp(\varepsilon) \cdot \exp(\delta)$ *and* $\log((1+\varepsilon)(1+\delta)) = \log(1+\varepsilon) + \log(1+\delta)$ *whenever* $\varepsilon, \delta \in \mathfrak{m}$ *commute.*

**Proof.** The part of the statement about the bijections is [7, Corollary 2.5]. By Lemma 3.17 and Corollary 3.24, we may assume that $A$ is abelian, as the subalgebra of $A$ generated by $\varepsilon$ and $\delta$ is abelian. Now $\mathrm{Ab}(k\langle\!\langle 2 \rangle\!\rangle) \simeq k[[X_0, X_1]]$ where $\exp(X_0) \cdot \exp(X_1) = \exp(X_0 + X_1)$ and $\log((1+X_0) \cdot (1+X_1)) = \log(1+X_0) + \log(1+X_1)$. We conclude by applying $\mathrm{ev}_{\varepsilon,\delta}^{\mathrm{Ab}} : k[[X_0, X_1]] \longrightarrow A$ via Proposition 3.34. $\square$

For $\varepsilon \in \mathfrak{m}$ and $\lambda \in k$, we set

$$(1+\varepsilon)^\lambda := \exp(\lambda \log(1+\varepsilon)) \in 1 + \mathfrak{m}. \tag{3.5}$$

**Proposition 3.39.** *The structure* $(1+\mathfrak{m}, \cdot, 1, (a \mapsto a^\lambda)_{\lambda \in k})$ *is an exponential group.*

**Proof.** Let $\varepsilon, \delta \in \mathfrak{m}$ and $\lambda, \mu \in k$. We clearly have $(1+\varepsilon)^0 = \exp(0) = 1$ and $(1+\varepsilon)^1 = \exp(\log(1+\varepsilon)) = (1+\varepsilon)$ by Proposition 3.38, so **EG1** holds. By Proposition 3.33, we have $(1+\varepsilon)^\lambda = \mathrm{ev}_\varepsilon(\exp(\lambda \log(1+X_0)))$ for the evaluation map $\mathrm{ev}_\varepsilon : k\langle\!\langle \{0\} \rangle\!\rangle \longrightarrow A$. So $(1+\varepsilon)^\lambda$ and likewise $(1+\varepsilon)^\mu$ are power series in $\varepsilon$, and they commute by **SA**. We deduce with Proposition 3.38 that $\log((1+\varepsilon)^\mu ((1+\varepsilon)^\lambda)) = \mu \log(1+\varepsilon) + \lambda \log(1+\varepsilon)$. So $(1+\varepsilon)^{\mu+\lambda} = \exp(\mu \log(1+\varepsilon) + \lambda \log(1+\varepsilon)) = \exp(\mu \log(1+\varepsilon)) \exp(\lambda \log(1+\varepsilon))$ by Proposition 3.38, i.e. $(1+\varepsilon)^{\mu+\lambda} = (1+\varepsilon)^\mu (1+\varepsilon)^\lambda$. So **EG2** holds. It is clear by definition that **EG3** holds. We have

$$((1+\varepsilon)^\lambda)^\mu = \exp(\mu \log((1+\varepsilon)^\lambda)) = \exp(\mu \lambda \log(1+\varepsilon)) = (1+\varepsilon)^{\mu\lambda} = (1+\varepsilon)^{\lambda\mu}.$$

Suppose that $1+\varepsilon$ and $1+\delta$ commute in $1+\varepsilon$, i.e. $\varepsilon$ and $\delta$ commute in $A$. As $\lambda \log(1+\varepsilon)$ and $\lambda \log(1+\delta)$ are power series in $\varepsilon$ and $\delta$ respectively with coefficients in the commutative ring $k$, they also commute by **SA**. So $\exp(\lambda \log(1+\varepsilon) + \lambda \log(1+\delta)) = \exp(\lambda \log(1+\varepsilon)) \cdot \exp(\lambda \log(1+\delta))$ by Proposition 3.38, i.e. $((1+\varepsilon)(1+\delta))^\lambda = (1+\varepsilon)^\lambda (1+\delta)^\lambda$. Thus **EG4** holds.

Lastly, the conjugation map $a \mapsto (1+\varepsilon) \cdot a \cdot (1+\varepsilon)^{-1} : A \longrightarrow A$ is a strongly linear automorphism by **SA**, so it commutes with evaluations of power series in $k\langle\!\langle \{0\} \rangle\!\rangle$, such as the maps $\exp$ and $\log$. Therefore $((1+\varepsilon)(1+\delta)(1+\varepsilon)^{-1})^\lambda = (1+\varepsilon)(1+\delta)^\lambda (1+\varepsilon)^{-1}$, i.e. **EG5** holds. $\square$

**Proposition 3.40.** *(mainly [25, Corollary 5.9]) There is a unique homomorphism of exponential groups* $\mathrm{Free}(I) \longrightarrow 1 + k\langle\!\langle I \rangle\!\rangle$ *that sends* $i$ *to* $1+X_i$ *for all* $i \in I$.

**Proof.** This follows as in [25, Corollary 5.9], where $k$ is assumed to be a binomial domain. This is the assumption used by Warfield [47, Theorem 10.24] to define the power map on $1 + k\langle\!\langle I \rangle\!\rangle$ satisfying a stronger condition than **EG4**. However, we prove without further assumption on $k$ that $1 + k\langle\!\langle I \rangle\!\rangle$ is an exponential group in Proposition 3.39, and Jaikin-Zapirain's arguments carry over without change. $\square$

We also define

$$\varepsilon * \delta := \log(\exp(\varepsilon) \cdot \exp(\delta)) \in \mathfrak{m}.$$

for all $\varepsilon, \delta \in \mathfrak{m}$. By Proposition 3.33, the operation $*$ is a law of exponential group on $\mathfrak{m}$ for which $\exp: \mathfrak{m} \longrightarrow 1 + \mathfrak{m}$ is an isomorphism of exponential groups. The power series $X_0 * X_1 \in k\langle\!\langle 2 \rangle\!\rangle_0$ has a well-known [32, 27, 46] expansion as a Lie series

$$X_0 * X_1 = X_0 + X_1 + \frac{1}{2}[X_0, X_1] + \cdots$$

in $k\langle\!\langle 2 \rangle\!\rangle$, i.e. a series which is a linear combination of iterated Lie brackets in $X_0$ and $X_1$. For now, we only note that for all $\varepsilon, \delta \in \mathfrak{m}$, we have

$$\varepsilon * \delta = \varepsilon + \delta + \frac{1}{2}[\varepsilon, \delta] + P_1(\varepsilon, \delta), \tag{3.6}$$

where $P_1$ is a sum of iterated Lie brackets of length $\geqslant 2$, and that we have

$$\begin{aligned}\log(\llbracket \exp(\varepsilon), \exp(\delta) \rrbracket) &= (-\log \varepsilon) * (-\log \delta) * (\log \varepsilon) * (\log \delta) \\ &= [\varepsilon, \delta] + P_2(\varepsilon, \delta),\end{aligned} \tag{3.7}$$

where $P_2$ is a sum of iterated Lie brackets of length $\geqslant 2$. See [44] for more details.

**Remark 3.41.** A precise computation of the coefficients in $P_1$ (see [11]) shows that a prime number $p$ does not divide any of the denominators of its coefficients of degree $<p$. This was used by Lazard [27] to give a correspondence between non-necessarily divisible $p$-groups and Lie algebras over fields of characteristic $p$ of nilpotency class $<p$. In that case the correspondence gives a first-order bi-interpretability relation between the group and the Lie algebra, and this was used recently by D'Elbée [13, 12] to solve cases of Wilson's conjecture.

## 3.7 Free Lie algebras with evaluations

Let $\mathrm{Lie}\langle\!\langle I \rangle\!\rangle$ be the smallest summability subalgebra of $[k\langle\!\langle I \rangle\!\rangle]$ that contains $\{X_i : i \in I\}$. The summability algebra $[k\langle\!\langle I \rangle\!\rangle]$ has evaluations by Lemma 3.23. As $\mathrm{Lie}\langle\!\langle I \rangle\!\rangle$ is a summability subalgebra of $[k\langle\!\langle I \rangle\!\rangle]$, it has evaluations by Corollary 3.24. So writing $\iota(i) = X_i$ for each $i \in I$, we have a unique strongly linear morphism $\mathrm{ev}_\iota : k\langle\!\langle I \rangle\!\rangle^{\mathrm{na}} \longrightarrow \mathrm{Lie}\langle\!\langle I \rangle\!\rangle$ that sends $X_i$ to $X_i$ for all $i \in I$.

Recall that $\mathrm{Lie}(k\langle\!\langle I \rangle\!\rangle^{\mathrm{na}})$ is the quotient of $k\langle\!\langle I \rangle\!\rangle^{\mathrm{na}}$ by the smallest summability ideal $\mathfrak{l}(\!(I)\!)$ containing all Lie identities

$$P \cdot P \qquad \text{and} \qquad P \cdot (Q \cdot R) + Q \cdot (R \cdot P) + R \cdot (P \cdot Q) \tag{3.8}$$

for $P, Q, R \in k\langle\!\langle I \rangle\!\rangle^{\mathrm{na}}$.

**Proposition 3.42.** *We have $\mathrm{Ker}(\mathrm{ev}_\iota) = \mathfrak{l}(\!(I)\!)$ and $\mathrm{Im}(\mathrm{ev}_\iota) = \mathrm{Lie}\langle\!\langle I \rangle\!\rangle$. In other words, we have an isomorphism*

$$\mathrm{Lie}(k\langle\!\langle I \rangle\!\rangle^{\mathrm{na}}) \longrightarrow \mathrm{Lie}\langle\!\langle I \rangle\!\rangle$$

*which is natural in $I \in \mathbf{Set}$.*

**Proof.** That $\mathfrak{l}(\!(I)\!) \subseteq \mathrm{Ker}(\mathrm{ev}_\iota)$ follows from the fact that $\mathrm{Lie}\langle\!\langle I \rangle\!\rangle \subseteq [k\langle\!\langle I \rangle\!\rangle^{\mathrm{na}}]$ is a Lie summability subalgebra. In order to prove the converse inclusion, we claim that we may assume that $I$ is finite. Indeed suppose that the $\mathfrak{l}(\!(I)\!) \supset \mathrm{Ker}(\mathrm{ev}_\iota)$ holds whenever $I$ is a finite set. Let $J \subseteq I$ be a finite subset. Consider the projections $\pi_J^{\mathrm{na}}$ and $\pi_J$ of Remark 3.11. By unicity of evaluation maps, we have $\mathrm{ev}_\iota \circ \pi_J^{\mathrm{na}} = \pi_J \circ \mathrm{ev}_\iota$. If $P \in \mathrm{Ker}(\mathrm{ev}_\iota)$, then $\pi_J^{\mathrm{na}}(P)$ is in $\mathfrak{l}(\!(J)\!)$ for all finite subsets $J \subseteq I$ by our hypothesis. But $P$ can be written as a sum $P = \sum_{J \subseteq I \text{ finite}} \sum_{\mathrm{content}(w) = J} P(w) X_w$, so each $\sum_{\mathrm{content}(w) = J} P(w) X_w = \pi_J^{\mathrm{na}}(P)$ lies in $\mathfrak{l}(\!(J)\!)$. Now $\mathfrak{l}(\!(I)\!)$ is a summability submodule containing all $\mathfrak{l}(\!(J)\!)$ for $J \subseteq I$ finite, so $P \in \mathfrak{l}(\!(I)\!)$. Thus in the sequel we assume that $I$ is finite.

Consider the subalgebra $\operatorname{Lie}\langle I\rangle$ (resp. $k\langle I\rangle$, and $k\langle I\rangle^{\mathrm{na}}$) of $[k\langle\!\langle I\rangle\!\rangle]$ (resp. $k\langle\!\langle I\rangle\!\rangle$, and $k\langle\!\langle I\rangle\!\rangle^{\mathrm{na}}$) generated by $\{X_i : i \in I\}$ and write $\mathfrak{l}(I)$ for the ideal generated by all identities (3.8) in $k\langle I\rangle^{\mathrm{na}}$. For $k = \mathbb{Q}$, the quotient map $k\langle I\rangle^{\mathrm{na}}/\mathfrak{l}(I) \longrightarrow \operatorname{Lie}\langle I\rangle$ is injective by the Poincaré-Birkoff-Witt theorem (see [41, Corollary 0.3]). As $k$ is flat as a $\mathbb{Q}$-module, we deduce that the quotient map $k\langle I\rangle^{\mathrm{na}}/\mathfrak{l}(I) \longrightarrow \operatorname{Lie}(I)$ is injective as a morphism of modules. Clearly its image is a subalgebra of $[k\langle I\rangle]$ containing $\{X_i : i \in I\}$, so $k\langle I\rangle^{\mathrm{na}}/\mathfrak{l}(I) \simeq \operatorname{Lie}\langle I\rangle$. Taking the completions with respect to the lower central filtrations, we obtain an isomorphism $\varphi : k\langle\!\langle I\rangle\!\rangle^{\mathrm{na}}/\mathfrak{l}((I)) \longrightarrow \operatorname{Lie}\langle\!\langle I\rangle\!\rangle$. Any $P \in k\langle\!\langle I\rangle\!\rangle^{\mathrm{na}}$ is the limit of the Cauchy sequence

$$P = \lim_{n \to +\infty} \left( \sum_{w \in I^{(\star)},\, \operatorname{rank}(w) \leqslant n} P(w)\, X_w \right)_{n \in \mathbb{N}}.$$

So the map $\varphi$ is the quotient map $k\langle\!\langle I\rangle\!\rangle^{\mathrm{na}}/\mathfrak{l}((I)) \longrightarrow \operatorname{Lie}\langle\!\langle I\rangle\!\rangle$. This proves our claims. It is clear that that the isomorphism is natural in $I$. □

### 3.8 The category of summability Lie algebras with evaluations

In this section we show that the category $\mathbf{\Sigma Alg}^{\mathrm{ev}}$ of summability algebras with evaluations is complete and cocomplete. The completeness follows from the completeness of $\mathbf{\Sigma Alg}$ (Proposition 3.21) and the fact that products in $\mathbf{\Sigma Alg}$ of summability algebras with evaluations have evaluations, and that equalizers of arrows in $\mathbf{\Sigma Alg}$ whose domain has evaluations have evaluations.

It is easy to see that coequalizers of arrows in $\mathbf{\Sigma Alg}$ whose codomain has evaluations have evaluations. Thus it remains to show that coproducts of summability Lie algebras with evaluations exist.

**Proposition 3.43.** *The category of summability algebras with evaluations has coproducts.*

**Proof.** Let $(A_x, \Sigma_x)_{x \in X}$ be a family of summability algebras with evaluations and set $A := \bigsqcup_{x \in X} A_x \times \{x\}$. The construction of the coproduct is similar to that in Proposition 3.19, based on a bornology on $A^{(\star)}$. The bornology $\mathcal{F}$ that we consider on $A^{(\star)}$ is the smallest containing $\mathcal{F}_0$ which is closed under products $(U, V) \mapsto (UV)$ and star operations $U \mapsto U^{(\star)}$. So $\mathcal{F} = \bigcup_{n \in \mathbb{N}} \mathcal{F}_n$ where $\mathcal{F}_{n+1} = \mathcal{F}_n \cup \{(UV), U^{(\star)} : U, V \in \mathcal{F}_n\}$. Given a family of strongly linear morphisms $\varphi_x : A_x \longrightarrow B$ into a summability algebra $B$ with evaluations, and $P \in k(A^{(\star)}; \mathcal{F})$, there is an $n \in \mathbb{N}$ with $\operatorname{supp} P \in \mathcal{F}_n$. There are $m \in \mathbb{N}$, $x_1, \ldots, x_m \in X$ and summable families $(s_1)_{s_1 \in S_1}, \ldots, (s_m)_{s_m \in S_m}$ in $A_{x_1}, \ldots, A_{x_m}$ such that writing $f$ for the summable family $\bigsqcup_{i=1}^m S_i \times \{x_i\} \longrightarrow k(A^{(\star)}; \mathcal{F})$ sending $(s_i, x_i)$ to $X_{(s_i, x_i)}$, the family $(X_u)_{u \in \operatorname{supp} P}$ is a subfamily of $(u(f))_{u \in (\bigsqcup_{i=1}^m S_i \times \{x_i\})^{(\star)}}$. By strong linearity of each $\varphi_x$, the family $g = (\bigsqcup_{i=1}^m \varphi_x) \circ f$ is summable in $B$, and since $B$ has evaluations, the family $(P(u)\, u(g))_{u \in \operatorname{supp} P}$ is summable. Defining $\psi(P)$ to be its sum, we obtain a strongly linear morphism $k(A^{(\star)}; \mathcal{F}) \longrightarrow B$ and we conclude as in the proof of Proposition 3.19. □

With Propositions 3.16 and 3.27, we obtain:

**Theorem 3.44.** *The categories of summability algebras and summmability Lie algebras with evaluations are complete and cocomplete.*

**Remark 3.45.** In contrast, the category of locally nilpotent Lie algebras does not have infinite products, nor finite coproducts.

# 4 Multipliability exponential groups

## 4.1 Infinite linearly ordered products

Let $(\mathcal{G},\cdot,1,(g\mapsto g^\lambda)_{\lambda\in k})$ be an exponential group. Let $\Pi=(\Pi_I)_{I\in\mathbf{Los}}$ be a family of partial maps $\Pi_I:\mathcal{G}^{\underline{I}}\longrightarrow\mathcal{G}$ whose domains are denoted $\operatorname{dom}\Pi_I$. We refer to Section 1.1 for notations regarding orderings. Given $I,J\in\mathbf{Los}$, consider the following axioms:

**MG1.** The set $\operatorname{dom}\Pi_I$ is a subgroup of $\mathcal{G}^{\underline{I}}$. Moreover, for $f,g\in\operatorname{dom}\Pi_I$, writing

$$f[g]:=(\Pi_{(i^+)^*}(f^{-1}\!\upharpoonright(\underline{i^+}))\cdot g(i)\cdot\Pi_{i^+}(f\!\upharpoonright(\underline{i^+})))_{i\in\underline{I}},$$

we have $f[g]\in\operatorname{dom}\Pi_I$ and

$$\Pi_I(f\cdot g)=(\Pi_I f)\cdot(\Pi_I f[g]).$$

**MG2.** We have $\operatorname{dom}\Pi_I\supseteq\mathcal{G}^{(\underline{I})}$, and for all $g\in\mathcal{G}^{(\underline{I})}$, writing $\operatorname{supp} g=\{i_1,\ldots,i_n\}$ where $n\in\mathbb{N}$ and $i_1<\cdots<i_n$, we have

$$\Pi_I g=g(i_1)\cdots g(i_n).$$

**MG3.** If $\varphi:I\longrightarrow J$ is an isomorphism and $g\in\operatorname{dom}\Pi_J$, then $g\circ\varphi\in\operatorname{dom}\Pi_I$, and

$$\Pi_I(g\circ\varphi)=\Pi_J g.$$

**MG4.** Let $(I_j)_{j\in\underline{J}}$ be a family with $I=\coprod_J(I_j)_{j\in\underline{J}}$. Let $g\in\operatorname{dom}\Pi_I$ and for each $j\in\underline{J}$, set $g_j:=g\!\upharpoonright\underline{I_j}$. Then

    **MG4a.** $g_j\in\operatorname{dom}\Pi_{I_j}$ for each $j\in\underline{J}$.

    **MG4b.** $(\Pi_{I_j}g_j)_{j\in\underline{J}}\in\operatorname{dom}\Pi_J$.

    **MG4c.** $\Pi_J((\Pi_{I_j}g_j)_{j\in\underline{J}})=\Pi_I g$.

**MG5.** If $I=I_1\amalg I_2$ and $(g,h)\in\operatorname{dom}\Pi_{I_1}\times\operatorname{dom}\Pi_{I_2}$, then $(g\sqcup h)\in\operatorname{dom}\Pi_I$.

**MG6.** Let $g\in\operatorname{dom}\Pi_I$ and $g_0\in\mathcal{G}$. Then $g_0\cdot g\cdot g_0^{-1}=(g_0\cdot g(i)\cdot g_0^{-1})_{i\in\underline{I}}\in\operatorname{dom}\Pi_I$ and

$$\Pi_I(g_0\cdot g\cdot g_0^{-1})=g_0\cdot(\Pi_I g)\cdot g_0^{-1}.$$

**MEG7.** Let $g\in\operatorname{dom}\Pi_I$. Then $g^{-1}\in\operatorname{dom}\Pi_{I^*}$ and

$$\Pi_{I^*}g^{-1}=(\Pi_I g)^{-1}.$$

**MEG8.** For $g\in\operatorname{dom}\Sigma_I$ and a family of functions $f_i:X_i\longrightarrow k$ with finite domains $X_i$, $i\in\underline{I}$, we have $(g(i)^{f_i(x)})_{i\in\underline{I}\wedge x\in X_i}\in\operatorname{dom}\Sigma_{\{(i,x):i\in\underline{I}\wedge x\in X_i\}}$.

**MEG9.** For each linear ordering $\triangleleft$ on $\underline{I}$, we have $\operatorname{dom}\Pi_I=\operatorname{dom}\Pi_{(\underline{I},\triangleleft)}$.

**Definition 4.1.** *If $\Pi=(\Pi_I)_{I\in\mathbf{I}}$ satisfies* **MG1**–**MEG9**, *then we say that $\Pi$ is* **multipliability structure** *on $(\mathcal{G},\cdot,1)$ and that $(\mathcal{G},\Pi)$ is a* **multipliability exponential group**.

**Remark 4.2.** *If $(\mathcal{G},\cdot,1,\Pi)$ satisfies* **MG4**, *then for $I=(\underline{I},<)\in\mathbf{Los}$, $f\in\operatorname{dom}\Pi_I$ and $X\subseteq\underline{I}$ such that $(X,<)$ is a convex subset of $I$, then we have $f\!\upharpoonright X\in\operatorname{dom}\Pi_{(X,<)}$. Indeed, we can write $I=(L,<)\amalg(X,<)\amalg(R,<)$ where*

$$(L,R):=(\{i\in\underline{I}:\forall x\in X,i<X\},\{i\in\underline{I}:\forall x\in X,i>x\}).$$

Therefore, this follows from **MG4a**.

**Remark 4.3.** In **MG1**, it is implicitly assumed that $\Pi_{(i^+)^*}(f^{-1}\restriction(\underline{i^+}))$ and $\Pi_{i^+}(f\restriction(\underline{i^+}))$ are well-defined. This follows from Remark 4.2 once one knows that $(\mathcal{G},\cdot,1,\Pi)$ satisfies **MG4** and **MEG7**.

We say that a family $g:\underline{I}\longrightarrow\mathcal{G}$ is *multipliable* if $g\in\operatorname{dom}\Pi_I$. Let $(\mathcal{G},\Pi^\mathcal{G})$ and $(\mathcal{H},\Pi^\mathcal{H})$ be multipliability exponential groups. A map $\Phi:\mathcal{G}\longrightarrow\mathcal{H}$ is said *strongly multiplicative* if for all $I\in\mathbf{Los}$ and $g\in\operatorname{dom}\Pi_I^\mathcal{G}$, we have $\Phi\circ g\in\operatorname{dom}\Pi_I^\mathcal{H}$ and $\Pi_I^\mathcal{H}(\Phi\circ g)=\Phi(\Pi_I^\mathcal{G}g)$. We write **ΠGr** for the category of multipliability exponential groups, with strongly multiplicative exponential group homomorphisms as arrows.

**Remark 4.4.** The set $\operatorname{Aut}^+(\mathcal{G})$ of strongly multiplicative automorphisms of exponential group of $(\mathcal{G},\Pi)$ is a group under composition. Contrary to the case of summability Lie algebras $L$, for which the subalgebra $\operatorname{Der}^+(L)\subseteq[\operatorname{Lin}^+(L)]$ of strongly linear derivations on $L$ is a summability Lie algebra, there is no natural structure of multiliability exponential group on $(\operatorname{Aut}^+(\mathcal{G}),\circ,\operatorname{id}_\mathcal{G})$.

## 4.2 Examples of multipliability exponential groups

Throughout Section 4.2, we fix an exponential group $(\mathcal{G},\cdot,1,(g\mapsto g^\lambda)_{\lambda\in k})$.

**Example 4.5.** It is routine to check that $\Pi^{\min}$ is a full multipliability structure on $\mathcal{G}$. We call $\Pi^{\min}$ the *minimal multipliability structure* on $\mathcal{G}$.

**Definition 4.6.** Let $(\mathcal{G},\cdot,1,\Pi)$ be a multipliability group. A **multipliability exponential subgroup** of $(\mathcal{G},\Pi)$ is a subgroup $\mathcal{H}\subseteq\mathcal{G}$ such that $\Pi_I h\in\mathcal{H}$ whenever $I\in\mathbf{Los}$ and $h\in\operatorname{dom}\Pi_I\cap\mathcal{H}^{\underline{I}}$. A **multipliability ideal** of $(\mathcal{G},\Pi)$ is an ideal of $\mathcal{G}$ which is a multipliability exponential subgroup of $(\mathcal{G},\Pi)$.

**Proposition 4.7.** *Let $(\mathcal{G},\cdot,1,\Pi)$ be a multipliability exponential group and let $\mathcal{H}$ be a multipliability exponential subgroup. The function $\Pi^\mathcal{H}$ given by $\Pi_I^\mathcal{H}:=\Pi_I\restriction\mathcal{H}^{\underline{I}}$ for each $I\in\mathbf{Los}$ is a multipliability structure on $\mathcal{H}$.*

**Proof.** It is easy to see that the axioms are preserved. □

For $I\in\mathbf{Los}$, we have a formal ordered product $\operatorname{OP}_I\in 1+k\langle\!\langle I\rangle\!\rangle_0$ given by

$$\operatorname{OP}_I:=\sum_{\substack{(i_1,\dots,i_n)\in(\underline{I})^\star \\ i_1<\cdots<i_n}}X_{(i_1,\dots,i_n)}.$$

This can be thought of as the formal expansion of the product $\cdots(1+X_i)\cdots(1+X_j)\cdots$ where $i<j$ range in $\underline{I}$. Let $A=k\oplus\mathfrak{m}$ be an associative unital summability algebra with unital evaluations. Then we define a partial map $\Pi_I$ with domain $\operatorname{dom}\Pi_I=\{a+1:a\in\operatorname{dom}\Sigma_I\}$, and which sends a family $a:\underline{I}\longrightarrow 1+\mathfrak{m}$ where $a-1$ is summable, to the element

$$\Pi_I a=\operatorname{ev}_{a-1}(\operatorname{OP}_I)\in 1+\mathfrak{m}.$$

See [3, Sections 4.1 and 4.2] for more details.

**Proposition 4.8.** *The group $1+\mathfrak{m}$ with the power map of Proposition 3.39 and the family $(\Pi_I)_{I\in\mathbf{Los}}$ defined above is a multipliability exponential group.*

**Proof.** The validity of the axioms **MG1**–**MEG7** and **MEG9** is [3, Theorem 4.15], whereas **MEG8** follows from **SM5** in $k+\mathfrak{m}$. □

**Lemma 4.9.** *A family $f: I \longrightarrow \mathfrak{m}$ is summable in $\mathfrak{m}$ if and only if $\exp \circ f$ is multipliable in $1 + \mathfrak{m}$.*

**Proof.** Suppose that $f$ is summable. As $(\exp(X_i) - 1)_{i \in I}$ is summable in $k \langle\!\langle I \rangle\!\rangle_0$ and $\mathfrak{m}$ has evaluations, the family $\mathrm{ev}^{\mathrm{Ass}}_f (\exp(X_i) - 1)_{i \in I} = (\exp(f_i) - 1)_{i \in I}$ is summable in $\mathfrak{m}$, i.e. $\exp \circ f$ is multipliable in $1 + \mathfrak{m}$. Suppose conversely that $\exp \circ f$ is multipliable in $1 + \mathfrak{m}$, i.e. $\exp \circ f - 1$ is summable in $\mathfrak{m}$. As $(\log(1 + X_i))_{i \in I}$ is summable in $k \langle\!\langle I \rangle\!\rangle_0$ and $\mathfrak{m}$ has evaluations, the family $\mathrm{ev}^{\mathrm{Ass}}_{\exp \circ f - 1}(\log(1 + X_i))_{i \in I} = f$ by Proposition 3.38, is summable in $\mathfrak{m}$. $\square$

**Lemma 4.10.** *Let $G \subseteq 1 + \mathfrak{m}$ be a subset such that $\log(G)$ is a summability subalgebra of $[\mathfrak{m}]$. Then $G$ is a multipliability exponential subgroup of $1 + \mathfrak{m}$.*

**Proof.** By [3, Proposition 4.16], we only need to show that $G$ is closed under the power map (3.5), which follows from the fact that $\log(G)$ is closed under scalar products. $\square$

Recall by Proposition 3.30 that if $\mathbb{A} = k((M))$ is an algebra of Noetherian series, with its ordering $\prec$ (see Section 3.1), then the associative unital algebra $A = k \, \mathrm{id}_\mathbb{A} + \mathrm{Lin}^+_\prec(\mathbb{A})$ has unital evaluations. Thus by Proposition 4.8 the group $\mathrm{id}_\mathbb{A} + \mathrm{Lin}^+_\prec(\mathbb{A})$ is a multipliability exponential group of strongly linear maps under composition.

**Example 4.11.** The subgroup $1\text{-}\mathrm{Aut}^+(\mathbb{A}) := \mathrm{Aut}(\mathbb{A}) \cap \mathrm{id}_\mathbb{A} + \mathrm{Lin}^+_\prec(\mathbb{A}) = \log(\mathrm{Der}^+(\mathbb{A}) \cap \mathrm{Lin}^+_\prec(\mathbb{A}))$ is a multipliability exponential subgroup of $\mathrm{id}_\mathbb{A} + \mathrm{Lin}^+_\prec(\mathbb{A})$, and hence a mutipliability exponential group.

**Example 4.12.** Let $\mathbb{T}^{\mathrm{flat}}$ denote the field of flat transseries (see [3, Section 6.1] and [4, Section 3.1]), which is a subfield of an $\mathbb{R}$-algebra of Noetherian series $\mathbb{R}((M))$. The group $\mathcal{T}$ of flat transseries that are tangent to the identity, under composition, is a multipliability exponential subgroup of the opposite group of $1\text{-}\mathrm{Aut}^+(\mathbb{R}((M)))$. This follows from [3, Proposition 6.4] and the fact [3, Lemma 6.2] that it is an exponential subgroup of $1\text{-}\mathrm{Aut}^+(\mathbb{R}((M)))$.

Let $(\mathcal{G}, \cdot, 1, \Pi)$ be a multipliability exponential group and let $(\mathcal{H}, \cdot, 1)$ be an exponential group. Let $\Phi: \mathcal{G} \longrightarrow \mathcal{H}$ be a homomorphism of exponential groups and suppose that $\mathrm{Ker}(\Phi)$ is a multipliability subgroup of $(\mathcal{G}, \Pi)$. We have [3, Section 2.3] a well-defined multipliability structure $\Pi^\Phi$ on $\Phi(\mathcal{G}) \subseteq \mathcal{H}$ where for all $I \in \mathbf{Los}$, we have

$$\mathrm{dom}\, \Pi_I^\Phi = \{\Phi \circ g : g \in \mathrm{dom}\, \Pi_I\} \qquad \text{and} \qquad \forall g \in \mathrm{dom}\, \Pi_I, \Pi_I^\Phi(\Phi \circ g) := \Phi(\Pi_I g).$$

**Proposition 4.13.** *The structure $(\Phi(\mathcal{G}), \Pi^\Phi)$ is a multipliability group. It is full if $\Pi$ is full.*

**Proof.** This is the same as [3, Proposition 2.8] except that we need to check that **MEG8** holds as well, but this follows from the validity of **MEG8** in $\mathcal{G}$ and the fact that $\Phi$ is a homomorphism of exponential groups. $\square$

### 4.3 Group evaluations

As a consequence of Lemma 4.10, we have

**Lemma 4.14.** *For all sets $I$, the set $\exp(\mathrm{Lie}\langle\!\langle I \rangle\!\rangle)$ is a multipliability exponential subgroup of $1 + k \langle\!\langle I \rangle\!\rangle_0$.*

Thus $\exp(\mathrm{Lie}\langle\!\langle I\rangle\!\rangle)$ is a multipliability exponential group. We write

$$\mathrm{Gr}\langle\!\langle I\rangle\!\rangle := \exp(\mathrm{Lie}\langle\!\langle I\rangle\!\rangle).$$

**Remark 4.15.** A family $(Q_n)_{n\in\mathbb{N}}$ in $\mathrm{Gr}\langle\!\langle I\rangle\!\rangle$ where $(\mathrm{val}\, Q_n)_{n\in\mathbb{N}}$ tends to $+\infty$ in $\mathbb{N}$ is always multipliable in $\mathrm{Gr}\langle\!\langle I\rangle\!\rangle$, and the sequence $(\Pi_{m=1}^n Q_m)_{n\in\mathbb{N}}$ converges in the valuation topology in $k\langle\!\langle I\rangle\!\rangle$. This is not the case for general families, e.g. $(\exp(X_n))_{n\in\mathbb{N}}$ is multipliable in $\mathrm{Gr}((\mathbb{N}))$ but $(\Pi_{m=1}^n \exp(X_m))_{n\in\mathbb{N}}$ does not converge in the valuation topology.

**Definition 4.16.** *Let $\mathcal{G}$ be a multipliability exponential group. We say that $\mathcal{G}$ **has evaluations** if for all multipliable families $f:I\longrightarrow \mathcal{G}$, there is a strongly linear exponential group homomorphism $\mathrm{ev}_f^{\mathrm{Gr}}:\mathrm{Gr}\langle\!\langle I\rangle\!\rangle \longrightarrow \mathcal{G}$ with $\mathrm{ev}_f^{\mathrm{Gr}}(\exp(X_i))=f(i)$ for all $i\in I$.*

We write $\mathbf{\Pi Gr}^{\mathrm{ev}}$ for the full subcategory of $\mathbf{\Pi Gr}$ of multipliability exponential groups with evaluations. By definition, we have:

**Lemma 4.17.** *Any multipliability exponential subgroup of a multipliability exponential group with evaluations has evaluations.*

We fix a set $E$ and consider the family $e:E\longrightarrow \mathrm{Gr}\langle\!\langle E\rangle\!\rangle$ with $e(i)=\exp(X_i)$ for all $i\in E$.

**Lemma 4.18.** *Let $w$ be a non-associative word on $E$. Then $\mathrm{val}(w[\![e]\!]-1)\geqslant \mathrm{rank}(w)$.*

**Proof.** We prove this by induction on $\mathrm{rank}(w)$. This is clear for $w\in E\sqcup\{\varnothing\}$. Suppose that $w=(uv)$ where $u,v\neq\varnothing$ and the result holds for $u,v$. So $u[\![e]\!]=1+P_u$ and $v[\![e]\!]=1+P_v$ where $\mathrm{val}(P_u)\geqslant\mathrm{rank}(u)$ and $\mathrm{val}(P_v)\geqslant\mathrm{rank}(v)$. We have

$$w[\![e]\!] = (1+P_u)^{-1}(1+P_v)^{-1}(1+P_u)(1+P_v).$$

Recall since $\deg(P_u)>0$ that $(1+P_u)^{-1}=1-P_u+P_u^2-P_u^3+\cdots$, and likewise for $P_v$. So by expanding, we obtain

$$w[\![e]\!] = 1+P_u P_v Q_1 + P_u^2 Q_2 + P_v P_u Q_3 + P_v^2 Q_4 \tag{4.1}$$

for some $Q_1,\ldots,Q_4\in 1+k\langle\!\langle E\rangle\!\rangle_0$. We deduce since $\mathrm{val}(P_u),\mathrm{val}(P_v)\geqslant 1$ and $\mathrm{rank}(w)\leqslant \max(\mathrm{rank}(u),\mathrm{rank}(v))+1$ that $\mathrm{val}(w[\![e]\!]-1)\geqslant\mathrm{rank}(w)$. We conclude by induction. $\square$

**Lemma 4.19.** *Let $w$ be a non-associative word on $E$. For all $i\in\mathrm{content}(w)$ and all $u\in E^{(\star)}$ with $X_u\in\mathrm{supp}\, w[\![e]\!]-1$, we have $i\in\mathrm{content}(u)$.*

**Proof.** We prove this by induction on $\mathrm{rank}(w)$. This is clear for $w\in E\sqcup\{\varnothing\}$. Suppose that $w=(uv)$ where $u,v\neq\varnothing$ and where the result holds for $u,v$. Let $i\in\mathrm{content}(w)$. So $i\in\mathrm{content}(u)$ or $i\in\mathrm{content}(v)$. In view of (4.1), monomials in $w[\![g]\!]-1$ are products of monomials in $u[\![g]\!]$ and $v[\![g]\!]$, so the result follows from the induction hypothesis. $\square$

**Lemma 4.20.** *The family $(w[\![e]\!])_{w\in E^{(\star)}}$ is multipliable in $\mathrm{Gr}\langle\!\langle E\rangle\!\rangle$.*

**Proof.** We need to show that $(w[\![e]\!]-1)_{w\in E^{(\star)}}$ is summable in $k\langle\!\langle E\rangle\!\rangle$. Now for each $w\in E^{(\star)}$ there is a $P_w\in k\langle\!\langle\mathrm{content}(w)\rangle\!\rangle$ with $\mathrm{val}(P_w)\geqslant\mathrm{rank}(w)$ by Lemma 4.18, with $w[\![e]\!]-1=\mathrm{ev}_{e-1}(P_w)$. For fixed content and rank, there are only finitely many such formal power series, so the family $(P_w)_{w\in E^{(\star)}}$ is summable in $k\langle\!\langle E\rangle\!\rangle$. Thus by Corollary 3.32 the family $(\mathrm{ev}_h(P_w))_{w\in J^{(\star)}}=(w[\![h]\!]-1)_{w\in J^{(\star)}}$ is summable in $k\langle\!\langle E\rangle\!\rangle$. $\square$

**Corollary 4.21.** *Let $\mathcal{G}$ be a multipliability exponential group with evaluations and let $f:E\longrightarrow \mathcal{G}$ be multipliable. Then $(w[\![f]\!])_{f\in E^{(\star)}}$ is multipliable.*

We define
$$\mathrm{BHC}_{(E,<)} := \log(\Pi_{(E,<)} e) \in \mathrm{Lie}\langle\!\langle E \rangle\!\rangle. \tag{4.2}$$

**Lemma 4.22.** *Let $I$ be a set and let $Q = (Q_i)_{i \in E}$ be a multipliable family in $\mathrm{Gr}\langle\!\langle I \rangle\!\rangle$. We have*
$$\log(\Pi_{(E,<)} Q) = \mathrm{ev}^{\mathrm{Lie}}_{\log \circ Q}(\mathrm{BCH}_{(E,<)}).$$

**Proof.** We have
$$\begin{aligned}
\log(\Pi_{(E,<)} Q) &= \log(\mathrm{ev}_{Q-1}(\mathrm{OP}_{(E,<)})) \\
&= \mathrm{ev}_{\mathrm{ev}_{Q-1}(\mathrm{OP}_{(E,<)})-1}(\log(1+X)),
\end{aligned}$$
and
$$\begin{aligned}
\mathrm{ev}^{\mathrm{Lie}}_{\log \circ Q}(\mathrm{BCH}_J) &= \mathrm{ev}^{\mathrm{Lie}}_{\log \circ Q}(\log(\Pi_{(E,<)} g)) \\
&= \mathrm{ev}^{\mathrm{Lie}}_{\log \circ Q}(\log(\mathrm{ev}_{e-1}(\mathrm{OP}_J))) \\
&= \mathrm{ev}^{\mathrm{Lie}}_{\log \circ Q}(\mathrm{ev}_{\mathrm{ev}_{e-1}(\mathrm{OP}_{(E,<)})-1}(\log(1+X))).
\end{aligned}$$

Now
$$\begin{aligned}
\mathrm{ev}_{\log \circ Q}(\mathrm{ev}_{e-1}(\mathrm{OP}_{(E,<)}) - 1) &= \mathrm{ev}_{\log \circ Q}(\mathrm{ev}_{e-1}(\mathrm{OP}_{(E,<)})) - 1 \\
&= \mathrm{ev}_{(\mathrm{ev}_{\log \circ Q}(e(i)-1))_{i \in E}}(\mathrm{OP}_{(E,<)}) - 1 \\
&\quad \mathrm{ev}_{Q-1}(\mathrm{OP}_{(E,<)}) - 1.
\end{aligned}$$

Thus $\log(\Pi_{(E,<)} Q) = \mathrm{ev}^{\mathrm{Lie}}_{\log \circ Q}(\mathrm{BCH}_{(E,<)})$. $\square$

We fix a well-ordering $<$ on $E$ and consider the corresponding graded ordering $<_{\mathrm{gr}}$ on $E^\star$ of Section 3.3. We define a valuation $v: \mathrm{Gr}\langle\!\langle E \rangle\!\rangle \longrightarrow E^\star$ (e.g. in the sense of [42] for the reverse ordering on $E^\star$) by
$$\forall Q \neq 1, v(Q) = \min(\mathrm{supp}(Q-1), <_{\mathrm{gr}}) \quad \text{and} \quad v(1) = \varnothing.$$

**Lemma 4.23.** *For all $Q \in \mathrm{Gr}\langle\!\langle E \rangle\!\rangle \setminus \{1\}$, there is a unique $w$ Lyndon word $w$ on $E$ such that there exists a $\lambda \in k^\times$, which is also unique, such that $v(w^{()}[\![e]\!]) = v(Q)$ and $v(Q u^{()}[\![e]\!]^{-\lambda} Q_1) > v(Q)$ for all $Q_1 \in \mathrm{Gr}\langle\!\langle E \rangle\!\rangle$ with $v(P) > v(Q)$.*

**Proof.** Set $P := \log(Q) \in \mathrm{Lie}\langle\!\langle E \rangle\!\rangle$. Note that for any finite subset $S$ of $\mathrm{supp}\, P$ writing $C := \bigcup_{v \in S} \mathrm{content}(v)$, we have $\pi_C(\log(Q)) \in \mathrm{Lie}\langle\!\langle C \rangle\!\rangle$ where $\pi_C$ is the projection $k\langle\!\langle E \rangle\!\rangle \longrightarrow k\langle\!\langle C \rangle\!\rangle$ of Remark 3.11. As $Q \neq 1$, we have $P \neq 0$, so $\mathrm{supp}\, P$ has a least element $u$. By Propositions 3.42 and 3.14, this word is a Lyndon word over $C$. Proposition 3.13 gives
$$u <_{\mathrm{gr}} \mathrm{supp}(u^{()}[\log \circ e] - X_u). \tag{4.3}$$

Note that for any $R \in \mathrm{Lie}\langle\!\langle C \rangle\!\rangle$ with $u <_{\mathrm{gr}} \mathrm{supp}\, R_1$, we have $(X_u + R_1) * (-R_1) = X_u + R_1 - R_1 + R_2$ for an $R_2$ whose support consists only of words of rank $> \mathrm{rank}(u)$. So $u^{()}[\log \circ e] = X_u * R_3$ and $P = X_u * R_4$ for an $R_3, R_4$ with $u <_{\mathrm{gr}} \mathrm{supp}\, R_3 \cup \mathrm{supp}\, R_4$. Now for any $v \in E^\star \setminus \{\varnothing\}$, $P_1 \in \mathrm{Lie}\langle\!\langle E \rangle\!\rangle$ with $v <_{\mathrm{gr}} \mathrm{supp}\, P_1$, and any $\lambda \in k$, we have $\exp((\lambda X_v) * P_1) = \exp(\lambda X_v) \cdot \exp(P_1) = (1 + \lambda X_u + \cdots) \cdot (1 + P_1 + \cdots) = 1 + \lambda X_v + P_2$ for a $P_2$ with $u <_{\mathrm{gr}} \mathrm{supp}\, P_2$.

So for each $P_3$ with $u <_{\mathrm{gr}} \mathrm{supp}\, P_3$, the word $u$ lies strictly $<_{\mathrm{gr}}$-below $\mathrm{supp}(Q - (1 + Q(u) X_u))$ and $(\exp(Q(u) X_u * P_4) - (1 + \lambda X_u))$. In particular $u = v(Q)$ and $(u, Q(u))$ is unique to satisfy $v(Q (\exp(X_u))^{-Q(u)} P)$ for all $P \in \mathrm{Gr}\langle\!\langle C \rangle\!\rangle$ with $v(Q) <_{\mathrm{gr}} v(P)$.

There remains to show in order to conclude that $u <_{\mathrm{gr}} \operatorname{supp} \exp(X_u) \cdot (u^{()}\llbracket e \rrbracket)^{-1}$. In view of (4.3), it suffices to show that $u <_{\mathrm{gr}} \operatorname{supp}(\log(u^{()}\llbracket e \rrbracket) - u^{()}[\log \circ e])$. We prove this relation by induction on the rank of $u$. If $u$ is a letter $i$ in $E$, as then $\log(u^{()}\llbracket e \rrbracket) - u^{()}[\log \circ e] = \log(i\llbracket e \rrbracket) - i[\log \circ e] = \log(\exp(X_i)) - X_i = 0$ so the result holds. Suppose that $\operatorname{rank}(u) > 1$ and that the result holds for all Lyndon words of rank $<\operatorname{rank}(u)$. We can write $u = u_0 u_1$ where $u_0, u_1$ are Lyndon words and $u_1$ has maximal length. Then by (3.7), we have

$$\log(u^{()}\llbracket e \rrbracket) - u^{()}[\log \circ e] = \log\bigl(\bigl[\!\bigl[u_0^{()}\llbracket e \rrbracket, u_1^{()}\llbracket e \rrbracket\bigr]\!\bigr]\bigr) - \bigl[u_0^{()}[\log \circ e], u_1^{()}[\log \circ e]\bigr]$$
$$= \bigl[\log u_0^{()}\llbracket e \rrbracket, \log u_1^{()}\llbracket e \rrbracket\bigr] + P - \bigl[u_0^{()}[\log \circ e], u_1^{()}[\log \circ e]\bigr]$$

where $\operatorname{val}(P) > \operatorname{rank}(u_0) + \operatorname{rank}(u_1) = \operatorname{rank}(u)$, so $u <_{\mathrm{gr}} \operatorname{supp} P$. Thus it suffices to show that the support of $\bigl[u_0^{()}\llbracket e \rrbracket, u_1^{()}\llbracket e \rrbracket\bigr] - \bigl[u_0^{()}[\log \circ e], u_1^{()}[\log \circ e]\bigr]$ lies above $u$. The induction hypothesis gives $P_0, P_1 \in k\langle\!\langle E \rangle\!\rangle$ with $u_0 <_{\mathrm{gr}} \operatorname{supp} P_0$ and $u_1 <_{\mathrm{gr}} \operatorname{supp} P_1$, and

$$\bigl(\log u_0^{()}\llbracket e \rrbracket, \log u_1^{()}\llbracket e \rrbracket\bigr) = \bigl(u_0^{()}[\log \circ e] + P_0, u_1^{()}[\log \circ e] + P_1\bigr).$$

Recall since $u$ is a Lyndon word that $u <_{\mathrm{gr}} u_1$. Now for all $M_0, M_1 \in k\langle\!\langle E \rangle\!\rangle$, we have $\operatorname{supp}[M_0, M_1] \subseteq \{(m_0 m_1), (m_1 m_0) : (m_0, m_1) \in \operatorname{supp} M_0 \times \operatorname{supp} M_1\}$, so we obtain $u <_{\mathrm{gr}} \operatorname{supp}\bigl[\log u_0^{()}\llbracket e \rrbracket, P_1\bigr]$, and $u <_{\mathrm{gr}} \operatorname{supp}[P_0, P_1]$ and $u < \operatorname{supp}\bigl[P_0, \log u_1^{()}\llbracket e \rrbracket\bigr]$ by definition of $<_{\mathrm{gr}}$. This concludes the proof. $\square$

**Proposition 4.24.** *Let $Q \in \operatorname{Gr}\langle\!\langle E \rangle\!\rangle$. There are unique families $(\lambda_\beta)_{\beta<\alpha} \in (k^\times)^\alpha$ and $(w_\beta)_{\beta<\alpha} \in (\operatorname{Lyndon}(E))^\alpha$ strictly increasing in $(E^\star, <_{\mathrm{gr}})$ such that $Q = \prod_{(\alpha,<)} \bigl(w_\beta^{()}\llbracket e \rrbracket^{\lambda_\beta}\bigr)_{\beta<\alpha}$.*

**Proof.** We construct the sequence $(\lambda_\beta, w_\beta)_{\beta<\alpha}$ by induction on the class of ordinals. At each step $\alpha \in \mathbf{On}$, we write $Q_\alpha = \prod_{(\alpha,<)} \bigl(w_\beta^{()}\llbracket e \rrbracket^{\lambda_\beta}\bigr)_{\beta<\alpha}$. For $\alpha \in \mathbf{On}$, our inductive hypothesis $H_\alpha$ is that a sequence $(\lambda_\beta, w_\beta)_{\beta<\alpha} \in (k^\times \times \operatorname{Lyndon}(E))^\alpha$ is defined that is unique to satisfy the following properties. For each $\gamma < \alpha$, then we have $\lambda_\gamma \in k^\times$, and for each $R \in \operatorname{Gr}\langle\!\langle E \rangle\!\rangle$ with $\{w_\beta : \beta < \alpha\} <_{\mathrm{gr}} \operatorname{supp} R$, we have $w_\gamma <_{\mathrm{gr}} v\bigl(Q^{-1}\bigl(\bigl(\prod_{(\alpha,<)} \bigl(w_\beta^{()}\llbracket e \rrbracket^{\lambda_\beta}\bigr)_{\beta<\alpha}\bigr) \cdot R\bigr)\bigr)$. So $H_0$ holds vacuously. Let $\alpha > 0$ be an ordinal such that $H_\beta$ holds for all $\beta < \alpha$. Suppose that $\alpha$ is a limit, so the sequence $(\lambda_\beta, w_\beta)_{\beta<\alpha}$ is already defined. Let $R \in \operatorname{Gr}\langle\!\langle E \rangle\!\rangle$ with $\{w_\beta : \beta < \alpha\} <_{\mathrm{gr}} \operatorname{supp} R$ and let $\gamma < \alpha$. Note by **MG4** that $\prod_{(\alpha,<)} \bigl(w_\beta^{()}\llbracket e \rrbracket^{\lambda_\beta}\bigr)_{\beta<\alpha} = \bigl(\prod_{(\gamma,<)} \bigl(w_\beta^{()}\llbracket e \rrbracket^{\lambda_\beta}\bigr)_{\beta<\gamma+1}\bigr) R_\gamma$ where $R_\gamma = \prod_{([\gamma+1,\alpha),<)} \bigl(w_\beta^{()}\llbracket e \rrbracket^{\lambda_\beta}\bigr)_{\gamma+1\leqslant\beta<\alpha}$. Recall by Lemma 4.23 that $w_\beta = v\bigl(w_\beta^{()}\llbracket e \rrbracket^{\lambda_\beta}\bigr)$ for each $\beta < \alpha$, so $w_\gamma <_{\mathrm{gr}} w_{\gamma+1} \leqslant_{\mathrm{gr}} \operatorname{supp} R_\gamma R$. Therefore $H_{\gamma+1}$ applied at $R R_\gamma$ gives

$$w_\gamma < v\Biggl(Q^{-1}\Biggl(\biggl(\prod_{(\gamma+1,<)} (w_\beta^{()}\llbracket e \rrbracket^{\lambda_\beta})_{\beta<\gamma+1}\biggr) R_\gamma R\Biggr)\Biggr)$$
$$= v\Biggl(Q^{-1}\Biggl(\biggl(\prod_{(\alpha,<)} (w_\beta^{()}\llbracket e \rrbracket^{\lambda_\beta})_{\beta<\alpha}\biggr) R\Biggr)\Biggr).$$

The sequence $(\lambda_\beta, w_\beta)_{\beta<\alpha}$ is unique by unicity in $H_\beta$ for all $\beta < \alpha$, so $H_\alpha$ holds. Suppose now that $\alpha = \gamma + 1$ is a successor. If the element $Q_\gamma = 1$, then we stop and obtain the desired result. Suppose that $Q_\gamma \neq 1$. Then by Lemma 4.23 applied at $Q_\gamma$, which by $H_\gamma$ satisfies $v(Q') > w_\gamma$, we can extend the sequence uniquely. As $(E^\star, <_{\mathrm{gr}})$ is a well-ordered set and $(w_\beta)_{\beta<\alpha}$ is always a strictly increasing sequence, the process must stop at an ordinal $\alpha$ below or equal to the order type of $(E^\star, <_{\mathrm{gr}})$, and we obtain $Q = \prod_{(\alpha,<)} \bigl(w_\beta^{()}\llbracket e \rrbracket^{\lambda_\beta}\bigr)_{\beta<\alpha}$. $\square$

**Proposition 4.25.** *Let $A = k \oplus \mathfrak{m}$ be an associative unital summability algebra with unital evaluations. Then the multipliability exponential group $1 + \mathfrak{m}$ has evaluations. Moreover, we have*

$$\mathrm{ev}_f^{\mathrm{Gr}} := \exp \circ \, \mathrm{ev}_{\log \circ f}^{\mathrm{Lie}} \circ \log \tag{4.4}$$

*for all multipliable families $f$ in $1 + \mathfrak{m}$.*

**Proof.** Let $f : I \longrightarrow 1 + \mathfrak{m}$ be multipliable. So $f - 1$ is summable in $\mathfrak{m}$, whence $\log \circ f$ is summable in $\mathfrak{m}$ (as $(\log(1 + X_i))_{i \in I}$ is summable in $k \langle\!\langle I \rangle\!\rangle$), hence in $[\mathfrak{m}]$. We set

$$\mathrm{ev}_f^{\mathrm{Gr}} := \exp \circ \, \mathrm{ev}_{\log \circ f}^{\mathrm{Lie}} \circ \log,$$

where $\mathrm{ev}_f^{\mathrm{Lie}} : \mathrm{Lie}\langle\!\langle I \rangle\!\rangle \longrightarrow [\mathfrak{m}]$ is the evaluation map. So $\mathrm{ev}_f^{\mathrm{Gr}}(\exp(X_i)) = \exp(\mathrm{ev}_{\log \circ f}^{\mathrm{Lie}}(X_i)) = \exp(\log(f(i))) = f(i)$ for all $i \in I$. As $\mathrm{ev}_{\log \circ f}^{\mathrm{Lie}}$ is linear, this map commutes with the power operations. Let us show that it is strongly multiplicative.

Let $J = (\underline{J}, <) \in \mathrm{Los}$ and let $Q = (Q_j)_{j \in \underline{J}}$ be a multipliable family in $\mathrm{Gr}\langle\!\langle E \rangle\!\rangle$. We have $\log(\Pi_J Q) = \mathrm{ev}_{\log \circ Q}^{\mathrm{Lie}}(\mathrm{BCH}_J)$ by Lemma 4.22 for $(E, <) = J$. We deduce that

$$\begin{aligned}
\log(\mathrm{ev}_f^{\mathrm{Gr}}(\Pi_J Q)) &= \mathrm{ev}_{\log \circ f}^{\mathrm{Lie}}(\log(\Pi_J Q)) \\
&= \mathrm{ev}_{\log \circ f}^{\mathrm{Lie}}(\mathrm{ev}_{\log \circ Q}^{\mathrm{Lie}}(\mathrm{BCH}_J)) \\
&= \mathrm{ev}_{\mathrm{ev}_{\log \circ f}^{\mathrm{Lie}}(\log(Q(j)))_{j \in \underline{J}}}^{\mathrm{Lie}}(\mathrm{BCH}_J)
\end{aligned}$$

Likewise, writing $h$ for the family $(\mathrm{ev}_f^{\mathrm{Gr}}(Q_j))_{j \in \underline{J}}$, Lemma 4.22 gives

$$\log(\Pi_J (\mathrm{ev}_f^{\mathrm{Gr}}(Q_j))_{j \in \underline{J}}) = \mathrm{ev}_{\log \circ h}(\mathrm{BCH}_J).$$

Now $\log(h(j)) = \mathrm{ev}_{\log \circ f}^{\mathrm{Lie}}(\log(Q(j)))$ for all $j \in \underline{J}$ by definition of $\mathrm{ev}_f^{\mathrm{Gr}}$. So $\log(\mathrm{ev}_f^{\mathrm{Gr}}(\Pi_J Q)) = \log(\Pi_J (\mathrm{ev}_f^{\mathrm{Gr}}(Q_j))_{j \in \underline{J}})$, whence by injectivity of $\log$ (see Proposition 3.38), we get

$$\mathrm{ev}_f^{\mathrm{Gr}}(\Pi_J Q) = \Pi_J (\mathrm{ev}_f^{\mathrm{Gr}}(Q_j))_{j \in \underline{J}},$$

i.e. $\mathrm{ev}_f^{\mathrm{Gr}}$ is strongly multiplicative. $\square$

**Corollary 4.26.** *The multipliability exponential group $\mathrm{Gr}\langle\!\langle E \rangle\!\rangle$ has evaluations.*

**Proof.** This follows from Lemmas 4.14 and 4.17. $\square$

**Corollary 4.27.** *The multipliability exponential group $\mathcal{T}$ of flat transseries of Example 4.12 has evaluations.*

**Proof.** This follows from Proposition 4.25, Lemma 4.17 and the fact that $\mathcal{T}$ is a multipliability subgroup of a group of the form $1 + \mathfrak{m}$ where $\mathfrak{m}$ has evaluations. $\square$

**Proposition 4.28.** *Let $\mathcal{G}$ be a multipliability exponential group with evaluations. Let $I = (\underline{I}, <) \in \mathrm{Los}$, let $f : E \longrightarrow \mathcal{G}$ and $P : \underline{I} \longrightarrow \mathrm{Gr}\langle\!\langle E \rangle\!\rangle$ be multipliable families. Then for all $Q \in \mathrm{Gr}\langle\!\langle \underline{I} \rangle\!\rangle$, we have*

$$\mathrm{ev}_f^{\mathrm{Gr}}(\mathrm{ev}_P^{\mathrm{Gr}}(Q)) = \mathrm{ev}_{(\mathrm{ev}_f^{\mathrm{Gr}}(P(i)))_{i \in \underline{I}}}^{\mathrm{Gr}}(Q).$$

**Proof.** As $\mathrm{ev}_f^{\mathrm{Gr}} \circ \mathrm{ev}_P^{\mathrm{Gr}}$ and $\mathrm{ev}_{(\mathrm{ev}_f^{\mathrm{Gr}}(P(i)))_{i \in \underline{I}}}^{\mathrm{Gr}}$ are both strongly multiplicative homomorphisms of exponential groups, and in view of Proposition 4.24, it suffices to prove the result for $Q = \exp(X_i)$ for each $i \in \underline{I}$. Let $i \in \underline{I}$. We have $\mathrm{ev}_f^{\mathrm{Gr}}(\mathrm{ev}_P^{\mathrm{Gr}}(\exp(X_i))) = \mathrm{ev}_f^{\mathrm{Gr}}(P(i)) = \mathrm{ev}_{(\mathrm{ev}_f^{\mathrm{Gr}}(P(i)))_{i \in \underline{I}}}^{\mathrm{Gr}}(\exp(X_i))$, hence the result. $\square$

**Remark 4.29.** Our conventions for the definitions of $\mathrm{Lie}\langle\!\langle I\rangle\!\rangle$ and $\mathrm{Gr}\langle\!\langle I\rangle\!\rangle$ correspond implicitly to the choice of Hopf algebra structure, and hence of primitite elements, on $k\langle\!\langle I\rangle\!\rangle$. We do not introduce summability Hopf algebras in the paper because algebras of Noetherian series and associative algebras of strongly linear maps on such algebras do not have natural structures of Hopf algebras.

**Remark 4.30.** One of the difficulties of working with nilpotency for an exponential group $\mathcal{G}$ is that elements $g$ in the lower central series may not be product of powers $w[\![\mathrm{id}_S]\!]^{\lambda_w}$ of pure group commutators for $S \subseteq \mathcal{G}$ finite as is the case in $\mathrm{Gr}\langle\!\langle I\rangle\!\rangle$ (see Proposition 4.24). Indeed the arguments of the commutators themselves may involve powers. For instance, expressions of the form

$$[\![[\![s_0, s_1]\!]^{\lambda_0}, [\![s_2, [\![s_3, s_4]\!]^{\lambda_1}]\!]^{\lambda_2}]\!] \tag{4.5}$$

are in $\mathcal{G}_3$. The Hall-Petresco identities [47, Chapter 6] (see Section 6.4), in the nilpotent or lower central complete cases, eliminate this distinction, and as we will see (Proposition 6.9), the existence of evaluations $\mathrm{Gr}\langle\!\langle I\rangle\!\rangle \longrightarrow \mathcal{G}$ entails these identities in the locally nilpotent case. This suggests that contrary to the case of summability algebras (see Lemma 3.23), the definition of having evaluations cannot be replaced by the multipliability criterion of Corollary 4.21 or even a stronger one involving expressions like (4.5), as this should not suffice to entail the Hall-Petresco identities.

# 5 The formal Lie correspondence

## 5.1 The Baker-Campbell-Hausdorff operations

Let $(L, \Sigma)$ be a summability Lie algebra with evaluations. We write $\mathrm{Gr}(L)$ for the set $L$ together with the following structure. For each $I = (\underline{I}, <) \in \mathbf{Los}$, we set $\mathrm{dom}\,\Pi_I = \mathrm{dom}\,\Sigma_{\underline{I}}$ and

$$\mathrm{BCH}_I := \log(\Pi_I\,(\exp(X_i))_{i \in \underline{I}}) \in \mathrm{Lie}\langle\!\langle I\rangle\!\rangle.$$

For $a \in \mathrm{dom}\,\Pi_I$, we define $\Pi_I a$ to be $\mathrm{ev}_a^{\mathrm{Lie}}(\mathrm{BCH}_I)$. For $a, b \in \mathrm{Gr}(L)$, we set

$$a * b := \Pi_{(2, <)}(a, b) = a + b + \frac{1}{2}[a, b] + \frac{1}{12}([a, [a, b]] - [b, [a, b]]) + \cdots \tag{5.1}$$

For $\lambda \in k$ and $a \in \mathrm{Gr}(L)$, we set $a^\lambda := \lambda\,a$.

**Lemma 5.1.** [46, Theorem 2.4.1] *If $\mathfrak{g}$ is a locally nilpotent Lie algebra over $\mathbb{Q}$, then $(\mathfrak{g}, *, 0, (a \mapsto a^q)_{q \in \mathbb{Q}})$ is a uniquely divisible group.*

For $n > 0$, we define

$$\mathrm{Lie}_n\langle\!\langle I\rangle\!\rangle := \{P \in \mathrm{Lie}\langle\!\langle I\rangle\!\rangle : \mathrm{val}(P) \geqslant n - 1\},$$

so $\mathrm{Lie}_1\langle\!\langle I\rangle\!\rangle = \mathrm{Lie}\langle\!\langle I\rangle\!\rangle$. For all $u, v \in I^\star$, we have $[X_u, X_v] \subseteq \{uv, vu\}$. It follows that for all $P, Q \in \mathrm{Lie}((k, I))$, we have

$$\mathrm{val}([P, Q]) \geqslant \max(\mathrm{val}(P), \mathrm{val}(Q)) + 1. \tag{5.2}$$

We deduce that $\mathrm{Lie}_n\langle\!\langle I\rangle\!\rangle$ is clearly a summability ideal of $\mathrm{Lie}\langle\!\langle I\rangle\!\rangle$. The quotient summability Lie algebra

$$\mathrm{Lie}^{/n}\langle\!\langle I\rangle\!\rangle := \mathrm{Lie}\langle\!\langle I\rangle\!\rangle / \mathrm{Lie}_n\langle\!\langle I\rangle\!\rangle,$$

is $n$-nilpotent by (5.2). Note that $\mathrm{Lie}((k, I))$ is the inverse limit

$$\mathrm{Lie}\langle\!\langle I\rangle\!\rangle = \varprojlim (\mathrm{Lie}^{/n}((I)))_{n>0}$$

where the morphisms $\mathrm{Lie}^{/n}\langle\!\langle I\rangle\!\rangle \longrightarrow \mathrm{Lie}^{/m}\langle\!\langle I\rangle\!\rangle$ for $n \geqslant m$ are the quotient maps.

**Lemma 5.2.** *If $I$ is a set, then $(\mathrm{Lie}\langle\!\langle I\rangle\!\rangle, 0, *, (g \mapsto g^\lambda)_{\lambda \in k})$ is an exponential group over $k$.*

**Proof.** Note that the axioms of exponential group over $k$ are preserved under inverse limits, so it suffices to show that for each $n > 0$, the structure $(\mathrm{Lie}^{/n}\langle\!\langle I\rangle\!\rangle, 0, *, (g \mapsto g^\lambda)_{\lambda \in k})$ is an exponential group over $k$. Fix an $n > 0$. We know by Lemma 5.1 that $(\mathrm{Lie}^{/n}\langle\!\langle I\rangle\!\rangle, 0, *, (g \mapsto g^q)_{q \in \mathbb{Q}})$ is a uniquely divisible group, so it suffices to check that the other axioms for the power map $k \times \mathrm{Lie}^{/n}\langle\!\langle I\rangle\!\rangle \longrightarrow \mathrm{Lie}^{/n}\langle\!\langle I\rangle\!\rangle$ hold for non-rational constants in this structure. Now the operation induced on $\exp(\mathrm{Lie}^{/n}\langle\!\langle I\rangle\!\rangle)$ by the scalar law on $\mathrm{Lie}_n\langle\!\langle I\rangle\!\rangle$ is the power map on the exponential subgroup $\exp(\mathrm{Lie}^{/n}\langle\!\langle I\rangle\!\rangle) \subseteq 1 + k\langle\!\langle I\rangle\!\rangle_0 / (1 + k\langle\!\langle I\rangle\!\rangle_0^n)$. As $1 + k\langle\!\langle I\rangle\!\rangle_0$ is an exponential group, we deduce that $\mathrm{Lie}^{/n}\langle\!\langle I\rangle\!\rangle$ is an exponential group. $\square$

**Proposition 5.3.** *The structure $(\mathrm{Gr}(L), *, 0, (a \mapsto a^\lambda)_{\lambda \in k}, \Pi)$ is a multipliability exponential group with evaluations.*

**Proof.** That $(\mathrm{Gr}(L), *, 0, (a \mapsto a^c)_{c \in k})$ is an exponential group follows from Lemma 5.2, using evaluation morphisms to recover the axioms of exponential groups. Let us show that $\Pi$ is a multipliability structure on $\mathrm{Gr}(L)$. By definition, it satisfies **MEG9**. As $(L, \Sigma)$ is a summability algebra and multipliability on $(\mathrm{Gr}(L), \Pi)$ coincides with summability $(L, \Sigma)$, all axioms regarding multipliability of families hold, i.e. **MG5**, **MEG8** and the clauses of multipliability in **MG2**, **MG3**, **MG4**, **MG6** and **MEG7** hold.

We now prove the other parts of the axioms, starting with **MG4c**. Let $g : I \longrightarrow \mathrm{Gr}(L)$ be multipliable and let $<$ be a linear ordering on $I$. Let $(I_j)_{j \in J}$ be a partition of $I$ into convex subsets $I_j, j \in J$ and let $<_J$ be the ordering $i <_J j \Longleftrightarrow I_i < I_j$ on $J$. We want to prove that

$$\Pi_{(J,<_J)}(\Pi_{(I_j,<)} g)_{j \in J} = \Pi_{(I,<)} g.$$

By definition, we have

$$\Pi_{(I,<)} g = \mathrm{ev}_g^{\mathrm{Lie}}(\mathrm{BCH}_{(I,<)})$$

and

$$\begin{aligned}\Pi_{(J,<_J)}(\Pi_{(I_j,<)} g)_{j \in J} &= \mathrm{ev}_{\left(\mathrm{ev}_{g\upharpoonright I_j}^{\mathrm{Lie}}(\mathrm{BCH}_{(I_j,<)})\right)}^{\mathrm{Lie}}(\mathrm{BCH}_{(J,<_J)}) \\ &= \mathrm{ev}_g^{\mathrm{Lie}}\left(\mathrm{ev}_{(\mathrm{BCH}_{(I_j,<)})_{j \in J}}^{\mathrm{Lie}}(\mathrm{BCH}_{(J,<_J)})\right),\end{aligned}$$

so it suffices to show that $\mathrm{BCH}_{(I,<)} = \mathrm{ev}_{(\mathrm{BCH}_{(I_j,<)})_{j \in J}}^{\mathrm{Lie}}(\mathrm{BCH}_{(J,<_J)})$. Composing with the exponential $\exp : \mathrm{Lie}\langle\!\langle I\rangle\!\rangle \longrightarrow \mathrm{Gr}\langle\!\langle I\rangle\!\rangle$ and using Lemma 4.22, this reduces to proving

$$\mathrm{OP}_{(I,<)} = \mathrm{ev}_{(\mathrm{OP}_{(I_j,<)})_{j \in J}}(\mathrm{OP}_{(J,<_J)}),$$

which follows by Lemma 4.14 from the validity of **MG4c** in $1 + k\langle\!\langle I\rangle\!\rangle_0$. So **MG4** holds. Likewise **MG2**, **MG3**, **MG4**, **MG6** and **MEG7** hold in $\mathrm{Gr}(L)$ as they holds in $\mathrm{Gr}\langle\!\langle I\rangle\!\rangle$ for all sets $I$.

For **MG1**, we need to show that given another multipliable $f: I \longrightarrow \mathrm{Gr}(L)$, the family $f[g]$ is multipliable in $\mathrm{Gr}(L)$, i.e. summable in $L$. For $i \in I$ consider the series

$$Q(i) := (\Pi_{i^*}(\exp(X_j))_{j>i}) \exp(X_{(i,\varsigma)}) (\Pi_{i^*}(\exp(X_j))_{j>i})^{-1}$$

in $1 + k\langle\!\langle I \sqcup I \times \{\varsigma\}\rangle\!\rangle_0$ where $I \times \{\varsigma\}$ is a disjoint copy of $I$. This lies in $\mathrm{Gr}((I \sqcup I \times \{\varsigma\}))$ so $\log(Q(i)) \in \mathrm{Lie}((\{i,(i,\varsigma)\}))$. The family $(\log(Q(i)))_{i \in I}$ is summable in $\mathrm{Lie}((I \sqcup I \times \{\varsigma\}))$. Now $f \sqcup g_\varsigma$ is summable in $L$ by **SM4** and **SM2**, where $g_\varsigma$ is the copy $g_\varsigma = \{((i,\varsigma), g(i)) : i \in I\}$ of $g$. Thus the family $f[g] = (\mathrm{ev}^{\mathrm{Lie}}_{f \sqcup g_\varsigma}(\log(Q_i)))_{i \in I}$ is summable in $L$. Again we obtain the identity $\Pi_{(I,<)}(f_i g_i)_{i \in I} = (\Pi_{(I,<)} f) \cdot (\Pi_{(I,<)} f[g])$ from Lemma 4.22 and the validity of the corresponding identity in $\mathrm{Gr}\langle\!\langle I\rangle\!\rangle$. □

**Proposition 5.4.** *Let $L_1, L_2$ be summability Lie algebras with evaluations, and let $\Phi: L_1 \longrightarrow L_2$ be a strongly linear morphism of algebras. Then $\Phi: \mathrm{Gr}(L_1) \longrightarrow \mathrm{Gr}(L_2)$ is a strongly multiplicative homomorphism of exponential groups.*

**Proof.** Let $a: I \longrightarrow \mathrm{Gr}(L_1)$ be multipliable. Then $a: I \longrightarrow L_1$ is summable, so $\Phi(a)$ is summable, hence multipliable in $\mathrm{Gr}(L_1)$. Let $<$ be a linear ordering on $I$. Since $\Phi$ is a strongly linear morphism of Lie algebras, the set of series $Q \in \mathrm{Lie}\langle\!\langle I\rangle\!\rangle$ such that $\mathrm{ev}^{\mathrm{Lie}}_{\Phi(a)}(Q) = \Phi(\mathrm{ev}^{\mathrm{Lie}}_a(Q))$ is a summablity subalgebra of $\mathrm{Lie}\langle\!\langle I\rangle\!\rangle$ that contains $\{X_i : i \in I\}$, hence it is $\mathrm{Lie}((k, I))$ itself. The series $\mathrm{BCH}_{(I,<)}$ lies in $\mathrm{Lie}\langle\!\langle I\rangle\!\rangle$, so $\Pi_{(I,<)} \Phi(a) = \mathrm{ev}^{\mathrm{Lie}}_{\Phi(a)}(\mathrm{BCH}_{(I,<)}) = \Phi(\mathrm{ev}^{\mathrm{Lie}}_a(\mathrm{BCH}_{(I,<)})) = \Phi(\Pi_{(I,<)} a)$. □

**Corollary 5.5.** *Summability ideals of $L$ are multipliability ideals of $\mathrm{Gr}(L)$.*

**Proposition 5.6.** *The multipliability exponential group $\mathrm{Gr}(L)$ has evaluations with $\mathrm{ev}^{\mathrm{Gr}}_g := \mathrm{ev}^{\mathrm{Lie}}_g \circ \log$ for all multipliable families $g$ in $\mathrm{Gr}(L)$.*

**Proof.** Let $g: I \longrightarrow \mathrm{Gr}(L)$ be multipliable. Then $\mathrm{ev}^{\mathrm{Lie}}_g: \mathrm{Gr}(\mathrm{Lie}\langle\!\langle I\rangle\!\rangle) \longrightarrow \mathrm{Gr}(L)$ is a strongly multiplicative exponential group homorphism by Proposition 5.4. Now $\exp: \mathrm{Gr}(\mathrm{Lie}\langle\!\langle I\rangle\!\rangle) \longrightarrow \mathrm{Gr}\langle\!\langle I\rangle\!\rangle$ is an isomorphism of multipliability exponential groups by Lemma 4.22, so $\mathrm{ev}^{\mathrm{Gr}}_g := \mathrm{ev}^{\mathrm{Lie}}_g \circ \log$ is the desired evaluation map. □

## 5.2 The Mal'cev-Lazard-Stewart operations

Let $I$ be a set. By definition of the functor $\mathrm{Gr}$, the map $\exp: \mathrm{Gr}(\mathrm{Lie}\langle\!\langle I\rangle\!\rangle) \longrightarrow \mathrm{Gr}\langle\!\langle I\rangle\!\rangle$ is an isomorphism of multipliability exponential groups. Thus by Corollary 5.5, for each $n > 0$, the set

$$\mathrm{Gr}_n\langle\!\langle I\rangle\!\rangle := \exp(\mathrm{Lie}_n\langle\!\langle I\rangle\!\rangle)$$

is a multipliability ideal of $\mathrm{Gr}\langle\!\langle I\rangle\!\rangle$. Furthermore, as $\exp$ is injective and $\bigcap_{n>0} \mathrm{Lie}_n\langle\!\langle I\rangle\!\rangle = \{0\}$, we have

$$\bigcap_{n>0} \mathrm{Gr}_n\langle\!\langle I\rangle\!\rangle = \{0\}. \tag{5.3}$$

**Lemma 5.7.** *For all $n > 0$, we have $\mathrm{Gr}\langle\!\langle I\rangle\!\rangle_n \subseteq \mathrm{Gr}_n\langle\!\langle I\rangle\!\rangle$.*

**Proof.** Recall that $\mathrm{Gr}\langle\!\langle I\rangle\!\rangle_n$ denotes the $n$-th term of the lower central series of the exponential group $\mathrm{Gr}\langle\!\langle I\rangle\!\rangle$ over $k$. We have $\mathrm{Gr}\langle\!\langle I\rangle\!\rangle_1 = \mathrm{Gr}\langle\!\langle I\rangle\!\rangle = \mathrm{Gr}_1\langle\!\langle I\rangle\!\rangle$. Let $n > 0$ such that $\mathrm{Gr}\langle\!\langle I\rangle\!\rangle_n \subseteq \mathrm{Gr}_n\langle\!\langle I\rangle\!\rangle$. In order to show that $\mathrm{Gr}\langle\!\langle I\rangle\!\rangle_{n+1} \subseteq \mathrm{Gr}_{n+1}\langle\!\langle I\rangle\!\rangle$, it suffices to prove that $\mathrm{Gr}_{n+1}\langle\!\langle I\rangle\!\rangle$ is an ideal of $\mathrm{Gr}\langle\!\langle I\rangle\!\rangle$ that contains $[\![\mathrm{Gr}\langle\!\langle I\rangle\!\rangle, \mathrm{Gr}\langle\!\langle I\rangle\!\rangle_n]\!]$. We know that it is an ideal of $\mathrm{Gr}((k, I))$, and by the inductive hypothesis it suffices to show that $[\![\mathrm{Gr}\langle\!\langle I\rangle\!\rangle, \mathrm{Gr}_n\langle\!\langle I\rangle\!\rangle]\!] \subseteq \mathrm{Gr}_{n+1}\langle\!\langle I\rangle\!\rangle$. Let $(f, g) \in \mathrm{Gr}\langle\!\langle I\rangle\!\rangle \times \mathrm{Gr}_n\langle\!\langle I\rangle\!\rangle$. We want to show that $\log(f^{-1}g^{-1}fg) \in \mathrm{Lie}_{n+1}\langle\!\langle I\rangle\!\rangle$. Now

$$P := \log(\exp(X_0)^{-1}\exp(X_1)^{-1}\exp(X_0)\exp(X_1)) = (-X_0) * (-X_1) * X_0 * X_1$$

is a Lie series in $\mathrm{Lie}\langle\!\langle 2\rangle\!\rangle \subseteq k\langle\!\langle 2\rangle\!\rangle_0$. By (3.6), we have $P = [X_0, X_1] + R$ where $R \in \mathrm{Lie}_2\langle\!\langle 2\rangle\!\rangle$. It follows that the single letter words 0 and 1 do not lie in $\operatorname{supp} P_\lambda$. with $1 \in \operatorname{content}(u)$ for all $u \in \operatorname{supp}(-X_0) * (-X_1) * X_0 * X_1$. In particular, as $\mathrm{Lie}_{n+1}\langle\!\langle I\rangle\!\rangle$ is an ideal containing all $[Q, \log(g)]$ for all $Q \in \mathrm{Lie}\langle\!\langle I\rangle\!\rangle$, we obtain $\log(f^{-1}g^{-1}fg) \in \mathrm{Lie}_{n+1}((k, I))$. We conclude by induction. $\square$

**Corollary 5.8.** *The exponential group $\mathrm{Gr}\langle\!\langle I\rangle\!\rangle$ is residually nilpotent.*

Let $\mathcal{G}$ be a multipliability exponential group with evaluations. We write $\mathrm{Lie}(\mathcal{G})$ for the set $\mathcal{G}$ together with the following structure. Let $b: J \longrightarrow \mathrm{Lie}(\mathcal{G})$ be a family. We say that $b$ is summable if it is multipliable. In that case, writing

$$\mathrm{MLS}_J := \exp\left(\sum_{j \in J} X_j\right) \in \mathrm{Gr}\langle\!\langle I\rangle\!\rangle,$$

we define $\Sigma_J b$ to be $\mathrm{ev}_b^{\mathrm{Gr}}(\mathrm{MLS}_J)$. We Also define

$$\begin{aligned}
a + b &:= \mathrm{ev}_{a,b}^{\mathrm{Gr}}(\exp(X_0 + X_1)), \\
[a, b] &:= \mathrm{ev}_{a,b}^{\mathrm{Gr}}(\exp([X_0, X_1])), \quad \text{and} \\
\lambda a &:= a^\lambda
\end{aligned}$$

for all $a, b \in \mathrm{Lie}(\mathcal{G})$ and $\lambda \in k$.

**Remark 5.9.** Set $(\xi_0, \xi_1) := (\exp(X_0), \exp(X_1))$. The series $\exp(X_0 + X_1)$ and $\exp([X_0, X_1])$ are infinite ordered products of families of powers of iterated commutators $w_n[\![(\xi_0, \xi_1)]\!]$, $w_n \in \{0, 1\}^{(\star)}$ of nondecreasing length. We have

$$\exp(X_0 + X_1) = \xi_0 \xi_1 [\![\xi_0, \xi_1]\!]^{-1/2} \cdots \quad \text{and} \tag{5.4}$$

$$\exp([X_0, X_1]) = [\![\xi_0, \xi_1]\!] [\![\xi_0, [\![\xi_0, \xi_1]\!]]\!]^{1/2} [\![\xi_1, [\![\xi_0, \xi_1]\!]]\!]^{1/2} \cdots. \tag{5.5}$$

As $\{0, 1\}$ is finite, the condition of Remark 4.15 is satisfied, and thus these products are ordered products in $\mathrm{Gr}\langle\!\langle 2\rangle\!\rangle$. See [27, (2.8)] and [46, Section 2.3] for more details and [11] for computations of the exponents.

**Lemma 5.10.** [46, Theorem 2.4.1] *If $G$ is a locally nilpotent uniquely divisible group, then $(G, +, 0, [\cdot, \cdot], (g \mapsto qg)_{q \in \mathbb{Q}})$ is a Lie algebra over $\mathbb{Q}$.*

**Lemma 5.11.** *If $I$ is a set, then $(\mathrm{Gr}((I)), +, 0, [\cdot, \cdot], (g \mapsto \lambda g)_{\lambda \in k})$ is a Lie algebra over $k$.*

**Proof.** This follows as in the proof of Lemma 5.2, considering $\mathrm{Gr}_n\langle\!\langle I\rangle\!\rangle$ instead of $\mathrm{Lie}_n\langle\!\langle I\rangle\!\rangle$. This time, the law induced on $\mathrm{Lie}_n\langle\!\langle I\rangle\!\rangle = \log\langle\!\langle I\rangle\!\rangle$ by the power law on $\mathrm{Gr}_n\langle\!\langle I\rangle\!\rangle$ is the scalar Lie algebra law, thus $\mathrm{Gr}_n\langle\!\langle I\rangle\!\rangle$ is a Lie algebra over $k$. $\square$

**Proposition 5.12.** *The structure $(\mathrm{Lie}(\mathcal{G}), +, [\cdot, \cdot], 1, (a \mapsto \lambda a)_{\lambda \in k}, \Sigma)$ is a summability Lie algebra.*

**Proof.** That $(\mathrm{Lie}(\mathcal{G}), +, [\cdot,\cdot], 1, (a \mapsto \lambda a)_{\lambda \in k})$ is a Lie algebra follows from Lemma 5.11, using evaluation morphisms to recover the axioms of Lie algebras. Let us show that $\Sigma$ is a summability structure on $\mathrm{Lie}(\mathcal{G})$. Let $f: I \longrightarrow \mathcal{G}$ be multipliable. Let $J$ be a set and let $\varphi: J \longrightarrow I$ be a bijection. Let $<$ be a linear ordering on $<$ and let $<_\varphi$ be its pullback on $J$, so $\varphi: (J, <_\varphi) \longrightarrow (I, <)$ is an isomorphism. We have $f \in \mathrm{dom}\,\Pi_{(I,<)}$ by **MEG9**, $f \circ \varphi \in \mathrm{dom}\,\Pi_{(J,<)}$ by **MG3**, so again by **MEG9** the family $f \circ \varphi$ is multipliable. We have $\Sigma_J f \circ \varphi = \mathrm{ev}^{\mathrm{Gr}}_{f \circ \varphi}(\exp(\sum_{j \in J} X_j))$. Now we have a unique isomorphism $\mathrm{ev}^{\mathrm{Gr}}_{\exp(X_\varphi)}: \mathrm{Gr}\langle\!\langle J \rangle\!\rangle \longrightarrow \mathrm{Gr}\langle\!\langle I \rangle\!\rangle$ sending each $\exp(X_j)$ to $\exp(X_{\varphi(j)})$ (that it is an isomorphism follows from the fact that we can give its inverse $\mathrm{ev}^{\mathrm{Gr}}_{\exp(X_{\varphi^{\mathrm{inv}}})}$) so by unicity of evaluation maps, we must have $\mathrm{ev}^{\mathrm{Gr}}_{f \circ \varphi} = \mathrm{ev}^{\mathrm{Gr}}_f \circ \mathrm{ev}^{\mathrm{Gr}}_{\exp(X_\varphi)}$. So $\Sigma_J f \circ \varphi = \mathrm{ev}^{\mathrm{Gr}}_f(\exp(\sum_{i \in I} X_i)) = \Sigma_I f$. Thus **SM2** holds in $\mathrm{Lie}(\mathcal{G})$.

The axiom **SM1** holds because because clearly $\mathrm{ev}^{\mathrm{Lie}}_f(\exp(\sum_{i \in I} X_i)) = \mathrm{ev}^{\mathrm{Lie}}_{f \restriction \mathrm{supp}\,f}(\exp(\sum_{i \in \mathrm{supp}\,f} X_i))$. The axioms **SM4**, **SM5** and **SM3a** hold by definition of the summability structure $\Sigma$ on $\mathrm{Lie}(\mathcal{G})$. There remains to prove that **SM3b** and **SM3c** hold. Let $J$ be a set and consider a partition $I =: \bigsqcup_{j \in J} I_j$ (allowing certain elements of the partition to be empty). Since $(\Sigma_{I_j}(X_i)_{i \in I_j})_{j \in J}$ is summable in $\mathrm{Lie}\langle\!\langle I \rangle\!\rangle$, $(\exp(\Sigma_{I_j}(X_i)_{i \in I_j}))_{j \in J}$ is multipliable in $\mathrm{Gr}\langle\!\langle I \rangle\!\rangle$ by Lemma 4.9. So the family $S := (\mathrm{ev}^{\mathrm{Gr}}_f(\exp(\Sigma_{I_j}(X_i)_{i \in I_j})))_{j \in J} = (\Sigma_{I_j} f \restriction I_j)_{j \in J}$ is multipliable in $\mathcal{G}$, hence summable in $\mathrm{Lie}(\mathcal{G})$, so **SM3b** holds. We have

$$\sum_{j \in J} \sum_{i \in I_j} f_i = \mathrm{ev}^{\mathrm{Gr}}_{(\mathrm{ev}^{\mathrm{Gr}}_f(\exp(\Sigma_{I_j}(X_i)_{i \in I_j})))_{j \in J}}\left(\exp\left(\sum_{j \in J} X_j\right)\right)$$

$$= \mathrm{ev}^{\mathrm{Gr}}_f\left(\mathrm{ev}^{\mathrm{Gr}}_{(\exp(\Sigma_{I_j}(X_i)_{i \in I_j}))_{j \in J}}\left(\exp\left(\sum_{j \in J} X_j\right)\right)\right) \quad \text{(by Proposition 4.28)}$$

$$= \mathrm{ev}^{\mathrm{Gr}}_f\left(\exp\left(\mathrm{ev}^{\mathrm{Lie}}_{(\Sigma_{I_j}(X_i)_{i \in I_j})_{j \in J}}\left(\sum_{j \in J} X_j\right)\right)\right) \quad \text{(by (4.4))}$$

$$= \mathrm{ev}^{\mathrm{Gr}}_f\left(\exp\left(\sum_{i \in I} X_j\right)\right) \quad \text{(by \textbf{SM3c} in } \mathrm{Lie}\langle\!\langle I \rangle\!\rangle)$$

$$= \sum_{i \in I} f_i.$$

Therefore **SM3c** holds in $\mathrm{Lie}(\mathcal{G})$. This concludes the proof. $\square$

**Proposition 5.13.** *Let $\mathcal{G}_1, \mathcal{G}_2$ be multipliability exponential groups with evaluations, and let $\Psi: \mathcal{G}_1 \longrightarrow \mathcal{G}_2$ be a strongly multiplicative homomorphism of exponential groups. Then $\Psi: \mathrm{Lie}(\mathcal{G}_1) \longrightarrow \mathrm{Lie}(\mathcal{G}_2)$ is a strongly linear morphism of algebras.*

**Proof.** The proof is the same as Proposition 5.4, relying on Proposition 4.24 to justify that $\Psi(\mathrm{ev}_a(Q)) = \mathrm{ev}_{\Psi(a)}(Q)$ for all multipliable $a: I \longrightarrow \mathcal{G}_1$ and all $Q \in \mathrm{Gr}\langle\!\langle I \rangle\!\rangle$. Indeed this holds for $Q \in \{\exp(X_i) : i \in I\}$ by definition of group evaluations. $\square$

**Proposition 5.14.** *The structure $(\mathrm{Lie}(\mathcal{G}), +, [\cdot,\cdot], 1, (a \mapsto \lambda a)_{\lambda \in k}, \Sigma)$ has evaluations and we have $\mathrm{ev}^{\mathrm{Lie}}_f = \mathrm{ev}^{\mathrm{Gr}}_f \circ \exp$ for all summable families $f$ in $(\mathrm{Lie}(\mathcal{G}), \Sigma)$.*

**Proof.** Let $f: I \longrightarrow \mathrm{Lie}(\mathcal{G})$ be multipliable. Then $\mathrm{ev}^{\mathrm{Gr}}_f: \mathrm{Lie}(\mathrm{Gr}\langle\!\langle I \rangle\!\rangle) \longrightarrow \mathrm{Lie}(\mathcal{G})$ is a strongly linear algebra morphism by Proposition 5.13. The map $\log: \mathrm{Lie}(\mathrm{Gr}\langle\!\langle I \rangle\!\rangle) \longrightarrow \mathrm{Lie}\langle\!\langle I \rangle\!\rangle$ is an isomorphism of summability Lie algebras by Lemma 4.22, so $\mathrm{ev}^{\mathrm{Lie}}_f := \mathrm{ev}^{\mathrm{Gr}}_f \circ \exp$ is the desired evaluation map. $\square$

## 5.3 The formal Lie correspondence

**Proposition 5.15.** *Let $L$ be a summability Lie algebra with evaluations and $\mathcal{G}$ be a multipliability exponential group with evaluations. Then $\mathrm{Lie}(\mathrm{Gr}(L)) = L$ and $\mathrm{Gr}(\mathrm{Lie}(\mathcal{G})) = \mathcal{G}$.*

**Proof.** We first prove the identity $\mathrm{Lie}(\mathrm{Gr}(L)) = L$. In view of the definitions, it suffices to show that the sum, Lie bracket and infinite sums defined on $\mathrm{Lie}(\mathrm{Gr}(L))$ are exactly that on $L$. By Proposition 5.6, we have $\mathrm{ev}_f^{\mathrm{Gr}} = \mathrm{ev}_f^{\mathrm{Lie}} \circ \log$ where $\mathrm{ev}_f^{\mathrm{Lie}} \colon \mathrm{Lie}\langle\!\langle J \rangle\!\rangle \longrightarrow L$ and $\mathrm{ev}_f^{\mathrm{Gr}} \colon \mathrm{Gr}\langle\!\langle J \rangle\!\rangle \longrightarrow \mathrm{Gr}(L)$ are the evaluation maps. So for $a, b \in L$, we have

$$\mathrm{ev}_{a,b}^{\mathrm{Gr}}(\exp(X_0 + X_1)) = a + b \quad \text{and} \quad \mathrm{ev}_{a,b}^{\mathrm{Gr}}(\exp([X_0, X_1])) = [a, b].$$

The summable families in $\mathrm{Lie}(\mathrm{Gr}(L))$ are exactly the multipliable families in $\mathrm{Gr}(L)$ which are exactly the summable families in $L$. Let $f \colon J \longrightarrow L$ be summable. We have $\Sigma_J f = \mathrm{ev}_f^{\mathrm{Gr}}(\mathrm{MLS}_J) = \mathrm{ev}_f^{\mathrm{Lie}}(\log(\mathrm{MLS}_J)) = \mathrm{ev}_f^{\mathrm{Lie}}(\sum_{j \in J} X_j) = \sum_J f$. So $\mathrm{Lie}(\mathrm{Gr}(L)) = L$. The proof of $\mathrm{Gr}(\mathrm{Lie}(\mathcal{G})) = \mathcal{G}$ is symmetric. $\square$

We have thus proved:

**Theorem 5.16.** *The correspondence*

$$L \mapsto \mathrm{Gr}(L) \quad \mathrm{Lie}(\mathcal{G}) \mapsfrom \mathcal{G}$$
$$(\Phi \colon L_1 \to L_2) \mapsto (\Phi \colon \mathrm{Gr}(L_1) \to \mathrm{Gr}(L_2)) \quad (\Phi \colon \mathrm{Lie}(\mathcal{G}_1) \to \mathrm{Lie}(\mathcal{G}_2)) \mapsfrom (\Phi \colon \mathcal{G}_1 \to \mathcal{G}_2)$$

*is an isomorphism of categories between the category $\mathbf{\Sigma Lie}^{\mathrm{ev}}$ of summability Lie algebras with evaluations and the category $\mathbf{\Pi Gr}^{\mathrm{ev}}$ multipliability exponential groups with evaluations.*

With Theorem 3.44, we deduce:

**Corollary 5.17.** *The category $\mathbf{\Pi Gr}^{\mathrm{ev}}$ is complete and cocomplete.*

**Proposition 5.18.** *Let $L$ be a summability Lie algebra with evaluations and let $\mathcal{G}$ be a multipliability exponential group with evaluations. Then summability ideals of $L$ are exactly multipliability ideals of $\mathrm{Gr}(L)$, and multipliability ideals of $\mathcal{G}$ are exactly summability ideals of $\mathrm{Lie}(\mathcal{G})$.*

**Proof.** Summability ideals of $L$ are exactly kernels of strongly linear morphisms defined on $L$, which are, by Propositions 5.4, 5.13 and 5.15, exactly kernels of strongly multiplicative homomorphisms of exponential groups defined on $\mathrm{Gr}(L)$, i.e. multipliability ideals of $\mathrm{Gr}(L)$. The other equivalence follows from Proposition 5.15. $\square$

**Lemma 5.19.** *Let $\mathcal{G}$ be a multipliability exponential group with evaluations and let $f, g \in \mathcal{G}$. Then a summability ideal of $\mathcal{G}$ contains $[\![f, g]\!]$ if and only if it contains $[f, g]$.*

**Proof.** We have $[\![f, g]\!] = [f, g] + \mathrm{ev}_{f,g}(R)$ where $R \in \mathrm{Lie}_2\langle\!\langle 2 \rangle\!\rangle$ by (3.6), whence $\mathrm{ev}_{f,g}(R)$ lies in all summability/multipliability ideals of $\mathrm{Lie}(\mathcal{G})$ containing $[f, g]$. On the other hand, we have $[f, g] = [\![f, g]\!] \mathrm{P}$ where by Remark 5.9, P is a product of powers of iterated commutators $w[\![(f, g)]\!]$ in $(f, g)$. Thus P lies in all multipliability/summability ideals of $\mathcal{G}$ containing $[\![f, g]\!]$. $\square$

**Lemma 5.20.** *Let $L$ be a summability Lie algebra with evaluations. Then the summability subalgebras of $L$ are exactly the multipliability exponential subgroups of $\mathrm{Gr}(L)$.*

**Proof.** Any summability subalgebra of $L$ is in particular closed under (5.1) and power maps, and hence is a multipliability exponential subgroup of $\mathrm{Gr}(L)$. Conversely, a multipliability exponential subgroup of $\mathrm{Gr}(L)$ is closed under (5.4) and (5.5) and scalar products, hence it is a summability subalgebra of $L$. □

**Lemma 5.21.** *Let $L$ be a summability Lie algebra with evaluations and let $a \in L$. Then*

  a) *the centraliser of $a$ in $L$ is the centraliser of $a$ in $\mathrm{Gr}(L)$*

  b) *the center of $L$ is the center of $\mathrm{Gr}(L)$.*

**Proof.** One notices that these are summability subalgebras and summability ideals of $L$ respectively, and in view of (5.1), (5.4) and (5.5), for $b \in L$, we have $[a,b] = 0$ if and only if $[\![a,b]\!] = 0$. □

## 5.4 Some specializations

We now show that the formal Lie correspondence extend the ones mentioned in the introduction.

**Proposition 5.22.** *Let $\mathfrak{g}$ be a locally nilpotent Lie algebra and let $n > 0$. Then $\mathfrak{g}_n = \mathrm{Gr}(\mathfrak{g})_n$. In particular $\mathfrak{g}$ is n-nilpotent if and only if $\mathrm{Gr}(\mathfrak{g})$ is n-nilpotent.*

**Proof.** By Proposition 3.25, $\mathfrak{g}$ has evaluations for the minimal summatility structure $\Sigma^{\min}$. By Proposition 5.18, summability ideals of $\mathfrak{g}$ and multipliability ideals of $\mathrm{Gr}(\mathfrak{g})$ are the same thing. By Lemma 5.19, for such an ideal, containing $[\mathfrak{g}, S]$ for a subset $S$ is equivalent to containing $[\![\mathrm{Gr}(\mathfrak{g}), S]\!]$. Thus by induction we have $\mathfrak{g}_n = \mathrm{Gr}(\mathfrak{g})_n$ for all $n > 0$. □

Note that this does not entail that for a nilpotent exponential group $\mathcal{G}$, the underlying Lie algebra of $\mathrm{Lie}(\mathcal{G})$ is nilpotent. Indeed to conclude that we would need to know that $\mathcal{G}$ with its minimal multipliability structure has evaluations.

**Proposition 5.23.** *Let $\mathcal{H}$ be an exponential group together with the minimal multipliability structure $\Pi^{\min}$. Then $(\mathcal{H}, \Pi^{\min})$ has evaluations if and only if $\mathcal{H}$ is locally nilpotent.*

**Proof.** Suppose that $(\mathcal{H}, \Pi^{\min})$ has evaluations, let $S \subseteq H$ be a finite subset. Fix $n > 0$. Let $\langle S \rangle$ denote the exponential subgroup of $\mathcal{H}$ generated by $S$. Let $e: S \longrightarrow \mathrm{Gr}\langle\!\langle S \rangle\!\rangle$ with $e(s) = \exp(X_s)$ for all $s \in S$. This is a multipliabiity subgroup of $\mathcal{H}$, hence by Lemma 4.17 it has evaluations. Consider the multipliable family $\iota = \mathrm{id}_S$ and the corresponding evaluation map $\mathrm{ev}_\iota^{\mathrm{Gr}}: \mathrm{Gr}\langle\!\langle S \rangle\!\rangle \longrightarrow \langle S \rangle$. As $S$ is finite, this is a surjective exponential group homomorphism. It follows that we have inclusions $\mathrm{Gr}\langle\!\langle S \rangle\!\rangle_n \supseteq \langle S \rangle_n$ of the terms of the lower central series. We have $\mathrm{Gr}\langle\!\langle S \rangle\!\rangle_n \subseteq \mathrm{Gr}_n\langle\!\langle S \rangle\!\rangle$ by Lemma 5.7. Consider the decomposition by Proposition 4.24 of a $Q \in \mathrm{Gr}_n\langle\!\langle S \rangle\!\rangle = \exp(\mathrm{Lie}_n\langle\!\langle S \rangle\!\rangle)$. Any word $w \in S^{(\star)}$ appearing with non-zero exponent in that decomposition has rank $\geqslant n$, so $w[\![e]\!]$ is an iterated commutator of length $n$ in elements of $\langle S \rangle$, hence it lies in $\langle S \rangle_n$. We deduce as $\langle S \rangle_n$ is a multipliability subgroup of $\mathcal{H}$ that $\mathrm{ev}_\iota^{\mathrm{Gr}}(Q) \in \langle S \rangle_n$. This shows that $\mathrm{Gr}\langle\!\langle S \rangle\!\rangle_n = \langle S \rangle_n$. Assume for contradiction that $\langle S \rangle$ is not nilpotent. So there is a family $(g_n)_{n>0} \in \prod_{n>0} \langle S \rangle_{n+1} \setminus \langle S \rangle_n$, and we find a family $(Q_n)_{n>0} \in \Pi_{n>0} \mathrm{Gr}_{n+1}\langle\!\langle S \rangle\!\rangle \setminus \mathrm{Gr}_n\langle\!\langle S \rangle\!\rangle$ with $g_n = \mathrm{ev}_\iota^{\mathrm{Gr}}(Q_n)$ for all $n > 0$. The family $(Q_n)_{n>0}$ is multipliable by Lemma 4.20 and Corollary 4.26, so the infinitely supported family $(g_n)_{n>0}$ is multipliable in $\mathcal{H}$: a contradiction. We deduce that $\langle S \rangle$ is nilpotent. Therefore $\mathcal{H}$ is locally nilpotent.

Conversely, suppose that $\mathcal{H}$ is a locally nilpotent exponential group and let $S \subseteq \mathcal{H}$ be a finite subset. We want to prove that there is a strongly multiplicative exponential group homomorphism $\operatorname{ev}_S^{\operatorname{Gr}} : \operatorname{Gr}\langle\!\langle S \rangle\!\rangle \longrightarrow \mathcal{H}$ with $\operatorname{ev}_S^{\operatorname{Gr}}(\exp(X_s)) = s$ for all $s \in S$. Let $n$ be the nilpotency class of the exponential subgroup $\langle S \rangle$ generated by $S$. Let $\operatorname{Free}(S)$ denote the free exponential group over $S$, and let $\varphi : \operatorname{Free}(S) \longrightarrow \langle S \rangle$ and $\psi : \operatorname{Free}(S) \longrightarrow \operatorname{Gr}\langle\!\langle S \rangle\!\rangle$ be the morphisms given by the universal property of the free exponential group with $\varphi(s) = s$ and $\psi(s) = \exp(X_s)$ for all $s \in S$. By Lemma 1.8, we can pass to the lower central completions and obtain exponential group homomorphisms

$$\tilde{\psi} : \widetilde{\operatorname{Free}(S)} \longrightarrow \widetilde{\operatorname{Gr}\langle\!\langle S \rangle\!\rangle} \qquad \text{and} \qquad \tilde{\varphi} : \widetilde{\operatorname{Free}(S)} \longrightarrow \widetilde{\langle S \rangle}.$$

Since $\langle S \rangle$ is nilpotent, it is isomorphic to its lower central completion. Let us show that $\tilde{\psi}$ is an isomorphism of exponential groups. We claim that $\psi$ is injective. Indeed, we have a unique exponential group homomorphism $\eta : \operatorname{Free}(S) \hookrightarrow 1 + k\langle\!\langle S \rangle\!\rangle$ with $\eta(s) = 1 + X_s$ for all $s \in S$.

Since $1 + k\langle\!\langle S \rangle\!\rangle$ has evaluations, we also have an exponential group homomorphism $\varepsilon : \operatorname{Gr}\langle\!\langle S \rangle\!\rangle \longrightarrow 1 + k\langle\!\langle S \rangle\!\rangle$ with $\varepsilon(\exp(X_s)) = 1 + X_s$ for all $s \in S$. Now $\varepsilon \circ \psi$ and $\eta$ must coincide. By Proposition 3.40, the map $\eta$ is injective, so $\psi$ is injective. It follows that $\tilde{\psi}$ is injective. Let $Q \in \operatorname{Gr}\langle\!\langle S \rangle\!\rangle$ and let $Q = \prod_{w \in S^{(\star)}} w[\![e]\!]^{c_w}$ be its decomposition as per Proposition 4.24. For all $w \in S^{(\star)}$, the element $w[\![e]\!] \in \operatorname{Gr}\langle\!\langle S \rangle\!\rangle$ is in the image of $\psi$. As $S$ is finite, sets of non-associative words over $S$ with bounded rank are finite. Thus all but finitely elements $w$ of $S^{(\star)}$ satisfy $w[\![e]\!] \notin \operatorname{Gr}\langle\!\langle S \rangle\!\rangle_n$. Thus $Q$ lies in the adherence of $\psi(\operatorname{Free}(S))$ in the lower central series topology on $\operatorname{Gr}\langle\!\langle S \rangle\!\rangle$. Therefore $\psi(\operatorname{Free}(S))$ is dense in $\widetilde{\operatorname{Gr}\langle\!\langle S \rangle\!\rangle}$, so $\tilde{\psi}(\widetilde{\operatorname{Free}(S)}) = \widetilde{\operatorname{Gr}\langle\!\langle S \rangle\!\rangle}$.

We obtain $\operatorname{ev}_S^{\operatorname{Gr}} : \operatorname{Gr}\langle\!\langle S \rangle\!\rangle \longrightarrow \mathcal{H}$ by composing the natural exponential group homomorphism $\operatorname{Gr}\langle\!\langle S \rangle\!\rangle \longrightarrow \widetilde{\operatorname{Gr}\langle\!\langle S \rangle\!\rangle}$, the inverse of $\tilde{\psi}$, the map $\tilde{\varphi}$, and the natural isomorphism $\widetilde{\langle S \rangle} \longrightarrow \langle S \rangle$. If $(Q_i)_{i \in I}$ is multipliable in $\operatorname{Gr}\langle\!\langle S \rangle\!\rangle$, then as $S$ is finite, we must have $\log(Q_i) \in \operatorname{Lie}_n\langle\!\langle S \rangle\!\rangle$, whence $Q_i \in \operatorname{Gr}_n\langle\!\langle S \rangle\!\rangle$ for all but finitely many $i \in I$. It follows since $\operatorname{ev}_S^{\operatorname{Gr}}$ ranges in $\langle S \rangle$ that $\operatorname{ev}_S^{\operatorname{Gr}}(Q_i) = 1$ for all but finitely many $i \in I$, i.e. $(\operatorname{ev}_S^{\operatorname{Gr}}(Q_i))_{i \in I}$ is finitely supported. Thus this family is multipliable in $(\mathcal{H}, \Pi^{\min})$ and $\operatorname{ev}_S^{\operatorname{Gr}}$ is trivially strongly multiplicative. $\square$

The proof also shows that for all finite sets, $\operatorname{Free}(S)$ embeds into $\operatorname{Gr}\langle\!\langle S \rangle\!\rangle$, and $\widetilde{\operatorname{Free}(S)} \simeq \widetilde{\operatorname{Gr}\langle\!\langle S \rangle\!\rangle}$. If $S$ is an infinite set, then $\operatorname{Free}(S)$ is the direct limit of the system $(\operatorname{Free}(S_0))_{S_0 \subseteq S \text{ is finite}}$. As $\operatorname{Gr}\langle\!\langle S \rangle\!\rangle$ contains the direct limit of $(\operatorname{Gr}\langle\!\langle S_0 \rangle\!\rangle)_{S_0 \subseteq S \text{ is finite}}$ (by applying the universal property of direct limits to the evaluation morphisms $\operatorname{Gr}\langle\!\langle S_0 \rangle\!\rangle \longrightarrow \operatorname{Gr}\langle\!\langle S \rangle\!\rangle$), we have an embedding $\operatorname{Free}(S) \longrightarrow \operatorname{Gr}\langle\!\langle S \rangle\!\rangle$ for all sets $S$. With Corollary 5.8, we deduce that:

**Proposition 5.24.** *For each set $I$, the free exponential group $\operatorname{Free}(I)$ is residually nilpotent.*

**Remark 5.25.** We do not have $\widetilde{\operatorname{Free}(I)} \simeq \widetilde{\operatorname{Gr}\langle\!\langle I \rangle\!\rangle}$ for infinite $I$. Indeed, the product $\Pi_{(\mathbb{N},<)} (\exp(X_n))_{n \in \mathbb{N}}$ exists in $\operatorname{Gr}\langle\!\langle \mathbb{N} \rangle\!\rangle$ whereas $(\Pi_{m=0}^n m)_{n \in \mathbb{N}}$ is not a Cauchy sequence in the lower central uniform structure on $\operatorname{Free}(\mathbb{N})$.

**Theorem 5.26.** *For $m > 0$, the formal Lie correspondence specializes to isomorphisms of categories between the categories of:*

*a) $m$-nilpotent Lie algebras and $m$-nilpotent exponential groups,*

b) *locally nilpotent Lie algebras and locally nilpotent exponential groups*

c) *lower central complete Lie algebras lower central complete exponential groups.*

**Proof.** We first prove a. If $\mathfrak{g}$ is an $m$-nilpotent Lie algebra, with its minimal summability structure then by Proposition 5.22, the underlying exponential group $\mathrm{Gr}(\mathfrak{g})$ is $m$-nilpotent. If $\mathcal{G}$ is an $m$-nilpotent exponential group, then with its minimal multipliability structure, it has evaluations by Proposition 5.23. Thus $\mathrm{Lie}(\mathcal{G})$ is well-defined, and as its summability structure is the minimal one by definition of the functor Lie, it must by Proposition 3.25 be locally nilpotent. So by Proposition 5.22, as the exponential group $\mathcal{G} = \mathrm{Gr}(\mathrm{Lie}(\mathcal{G}))$ is $m$-nilpotent, the underlying Lie algebra of $\mathrm{Lie}(\mathcal{G})$ is $m$-nilpotent. We obtain b as a consequence of Propositions 3.25 and 5.23. Lastly follows c from a by taking inverse limits. □

**Remark 5.27.** The correspondence between formal groups seen as locally conilpotent cocommutative Hopf algebras (see [19, Section 7.2]) and Lie algebras is orthogonal to the formal Lie correspondence. Indeed, although its objects are colimits of nilpotent objects, and the points of a formal groups form residually nilpotent groups, the corresponding Lie algebras are arbitrary Lie algebras whose structure does not come from the inversion of a Baker-Campbell-Hausdorff operation.

Some examples of summability Lie algebras and multipliability exponential groups with evaluations that we have seen so far are residually nilpotent. Here we give an example which is perfect (and non-trivial), i.e. with $\mathcal{G}_n = \mathcal{G}$ for all $n > 0$.

**Example 5.28.** Consider the Hahn series field $\mathbb{K} = \mathbb{R}((\mathbb{R}))$ where $\mathbb{R}$ is seen as the reverse-ordered additive ordered group $(\mathbb{R}, +, 0, >)$ of real numbers. We write $x^\lambda = \mathbb{1}_\lambda \in \$$ for each $\lambda \in \mathbb{R}$, and set $x := x^1$. This is a real closed field [20, 31]. We have a unique strongly linear derivation $\partial : \mathbb{K} \longrightarrow \mathbb{K}$ with $\partial(\mathbb{1}_\lambda) = \lambda \mathbb{1}_{\lambda-1}$ for all $\lambda \in \mathbb{R}$, so $\partial(x) = 1$. Consider the relation $\prec$ of (3.2), and the summability Lie algebra over $\mathbb{R}$

$$L = \{a \in \mathbb{K} : a \prec x\} = \{a \in \mathbb{K} : \mathrm{supp}\, f < 1\},$$

under the operation $[a, b] := a\, \partial(b) - \partial(a)\, b$. This algebra has evaluations by [4, after Definition 1.1]. For $\lambda, \mu < 1$ with $\lambda \neq \mu$, we have $[x^\lambda, x^\mu] = (\mu - \lambda)\, x^{\lambda+\mu-1}$. Thus $[L, x^\mu] = \{a \in \mathbb{K} : a \prec x^\mu\}$. We deduce that

$$[L, b] = \{a \in \mathbb{K} : a \prec b\} \tag{5.6}$$

for all $b \in L \setminus \{0\}$, and thus that $[L, L] = L$. So $L$ is perfect.

The exponential map of [7] sends $L$ onto a subgroup $\mathcal{P}$ of the opposite group $\mathcal{T}^{\mathrm{op}}$ of the group $\mathcal{T}$ of flat transseries of Example 4.12 (see [4, Proposition 3.1]), which consists in series of the form $x + \delta$ where $\mathrm{supp}\, \delta \subseteq \mathbb{R}$ and $\delta \prec x$. More precisely, we have $\mathcal{P} \simeq \mathrm{Gr}(L)$. Now (5.6) and Lemma 5.19 entail that the smallest multipliability ideal containing $[\![\mathcal{P}, \mathcal{P}]\!]$ is $\mathcal{P}$, so $\mathcal{P}$ is perfect as well.

**Remark 5.29.** As $\mathcal{P}$ lies in the perfect hull of $\mathcal{T}$, this exponential group is not hypocentral. Our arguments can be used to show that $\mathcal{T}$ is hypoabelian of length $\omega$ (i.e. its derived series has trivial intersection). It is also a direct limit $\mathcal{T} = \varinjlim (\mathcal{T}_{(n)})_{n \in \mathbb{N}}$ of lower central complete multipliability exponential groups $\mathcal{T}_n = \{f \in \mathcal{T} : f - x \preccurlyeq x^{1-2^{-n}}\}$. Replacing $\mathbb{R}$ with a non-archimedean field in Example 5.28, one also obtains an example which is hypoabelian of length $\geqslant \omega_1$ and is not a colimit of lower central complete multipliability exponential groups with evaluations. However we do now know of examples of multipliability exponential groups that are not hypoabelian of possibly large transfinite length.



# 6 Applications

## 6.1 Mixed structures

Let $\mathcal{L}_{k,0}$ denote the one-sorted first-order language

$$\mathcal{L}_{k,0} = \langle +, 0, [\cdot, \cdot], *, (\lambda \,.\, )_{\lambda \in k} \rangle.$$

For each $m \in \mathbb{N}$, we write $\text{Term}_{k,0}(m)$ for the set of $\mathcal{L}_{k,0}$-terms in $m$ variables. We interpret any summability Lie algebra with evaluations $L$ and any multipliability exponential group with evaluations $\mathcal{G}$ as $\mathcal{L}_{k,0}$-structures by interpreting the alien function symbols in $\text{Gr}(L)$ and $\text{Lie}(\mathcal{G})$ respectively. We shall call $\mathcal{L}_{k,0}$-structures that are both Lie algebras and exponential groups *mixed structures* (over $k$). For instance $\mathcal{K} = (k, +, 0, 0, +, (\lambda \cdot )_{\lambda \in k})$ is a mixed structure for the minimal summability structure on $k$. We can see any mixed structure $L$ as an $\mathcal{L}_k$-structure where

$$\mathcal{L}_k = \langle (f)_{f \in \text{Lie}\langle\!\langle m \rangle\!\rangle, m \in \mathbb{N}} \rangle,$$

and each $f \in \text{Lie}\langle\!\langle m \rangle\!\rangle$ is interpreted as the function $L^m \longrightarrow L : \bar{b} \mapsto \text{ev}_{\bar{b}}^{\text{Lie}}(f)$. Given $m \in \mathbb{N}$, we write $\text{Term}_k(m)$ for the set of $\mathcal{L}_k$-terms in $m$ variables.

**Lemma 6.1.** *There are unique maps $\text{Term}_{k,0}(m) \longrightarrow \text{Lie}\langle\!\langle m \rangle\!\rangle ; t \mapsto \hat{t}$ and $\text{Term}_k(m) \longrightarrow \text{Lie}\langle\!\langle m \rangle\!\rangle ; t \mapsto \hat{t}$ for each $m \in \mathbb{N}$ such that for all composable terms $t \in \text{Term}_{k,0}(m)$, $t_1, \ldots, t_m \in \text{Term}_{k,0}(l)$, or $t \in \text{Term}_k(m)$, $t_1, \ldots, t_m \in \text{Term}_k(l)$, we have*

$$\widehat{t(t_1, \ldots, t_m)} = \hat{t}(\hat{t_1}, \ldots, \hat{t_m}).$$

*Moreover, for all $m \in \mathbb{N}$, $t \in \text{Term}_{k,0}(m), \text{Term}_k(m)$, all mixed structures $L$ and all $\bar{a} \in L^m$, we have*

$$t(\bar{a}) = \hat{t}(\bar{a}). \tag{6.1}$$

**Proof.** For $\mathcal{L}_{k,0}$-terms, we simply define the maps by induction on the complexity of terms, starting for function symbols in $\mathcal{L}_{k,0}$ with the correspondence

$$(+, 0, [\cdot, \cdot], *, (\lambda \,.\, )_{\lambda \in k}) = \left( X_0 + X_1, 0, [X_0, X_1], X_0 + X_1 + \frac{1}{2}[X_0, X_1] + \cdots, (\lambda)_{\lambda \in k} \right).$$

For $\mathcal{L}_k$-terms, we proceed similarly, defining $\widehat{t(t_1, \ldots, t_m)} := \text{ev}_{(t_1, \ldots, t_m)}^{\text{Lie}}(t)$ for all tuples of composable terms. The identity (6.1) follows by induction from Proposition 3.33. $\square$

Thus terms in $\mathcal{L}_{k,0}$ and $\mathcal{L}_k$ have canonical representations as function symbols in $\mathcal{L}_k$.

**Lemma 6.2.** *For all $m \in \mathbb{N}$, $c > 0$ and all $t \in \text{Term}(m)$, there are an $\mathcal{L}_{k\text{-Lie}}$-term $t_{c,\text{Lie}}(\bar{x})$ and an $\mathcal{L}_{k\text{-Gr}}$-term $t_{c,\text{Gr}}(\bar{x})$ a such that for all $c$-nilpotent mixed structures $L$, and all $\bar{a} \in L^m$, we have $t(\bar{a}) = t_{c,\text{Lie}}(\bar{a}) = t_{c,\text{Gr}}(\bar{a})$. Moreover, one can choose the families $(t_{c,\text{Lie}}(\bar{x}))_{c>0}$ such that $(\widehat{t_{c+1,\text{Lie}}}(\bar{x}) - \widehat{t_{c,\text{Lie}}}(\bar{x}))_{c>0}$ is summable in $\text{Term}(m)$ with*

$$\sum_{i>0} \widehat{t_{c+i+1,\text{Lie}}}(\bar{x}) - \widehat{t_{c+i,\text{Lie}}}(\bar{x}) = 0 \tag{6.2}$$

*for all $c > 0$.*



**Proof.** Recall that $\mathrm{Lie}\langle\!\langle m\rangle\!\rangle \simeq \mathrm{Lie}(k\langle\!\langle m\rangle\!\rangle^{\mathrm{na}})$ by Proposition 3.42, so we can write $t(\bar{x}) = \mathrm{ev}_{(X_0,\ldots,X_{m-1})}(Q)$ for a $Q \in k\langle\!\langle m\rangle\!\rangle^{\mathrm{na}}$. Now write $Q = Q_{<c} + Q_{\geqslant c}$ where

$$Q_{\leqslant c} = \sum_{u \in m^{(\star)}, \mathrm{rank}(u) \leqslant c} Q(u)\, X_u \quad \text{and}$$

$$Q_{>c} = \sum_{u \in m^{(\star)}, \mathrm{rank}(u) > c} Q(u)\, X_u,$$

and set $P_{\leqslant c} = \mathrm{ev}_{(X_0,\ldots,X_{m-1})}(Q_{\leqslant c})$ and $P_{>c} = \mathrm{ev}_{(X_0,\ldots,X_{m-1})}(Q_{>c})$ so $t(\bar{x}) = P_{\leqslant c} + P_{>c}$.

Given $L$, $m$ and $\bar{a}$ as in the statement, we have $\mathrm{ev}_{\bar{a}}^{\mathrm{Lie}}(P_{>c}) \in L_c = \{0\}$ so $t(\bar{a}) = \mathrm{ev}_{\bar{a}}^{\mathrm{Lie}}(P_{<c})$. Thus we can set $t_{c,\mathrm{Lie}}(\bar{x})$ to be any term corresponding to $P_{\leqslant c}$ in all Lie algebras. We see that $(\widehat{t_{c+1,\mathrm{Lie}}(\bar{x})} - \widehat{t_{c+1,\mathrm{Lie}}(\bar{x})})_{c>0} = (P_{\leqslant c+1} - P_{\leqslant c})_{c>0}$ is summable as its terms have piecewise distinct supports. We clearly have $P_{\leqslant c} = \sum_{i>0} P_{\leqslant c+1} - P_{\leqslant c}$

Now $P_{<c} \in \mathrm{Lie}\langle m\rangle$ so by Proposition 3.14, we can write $P$ as a finite linear combination of terms of the form $u[(X_0,\ldots,X_{m-1})]$ for Lyndon words $u$. Replacing each Lie bracket and then its sum iteratively by the finite products of powers of Remark 5.9 truncated at iterator depth $c-1$, we obtain a finite product of powers, i.e. a term in $t_{c,\mathrm{Gr}}(\bar{x})$ in $\mathcal{L}_{k\text{-Gr}}$, that evaluates as $t(\bar{x})$ in all $c$-nilpotent exponential groups, hence in all $c$-nilpotent mixed structures. $\square$

## 6.2 Transfer theorems

In this section, we assume that $k$ is a field of characteristic zero. We define the dimension of a mixed structure as the dimension its underlying vector space. We say that a mixed structure is nilpotent if the underlying Lie algebra / exponential group is nilpotent, and then define its nilpotency class as that of the latter.

**Lemma 6.3.** *Let $m \geqslant 2$ and $n \geqslant 1$. Then $\mathrm{Lie}\langle\!\langle m\rangle\!\rangle / \mathrm{Lie}_n\langle\!\langle m\rangle\!\rangle$ is a finite dimensional nilpotent Lie algebra of nilpotency class $n$.*

**Proof.** By definition of $\mathrm{Lie}_n\langle\!\langle m\rangle\!\rangle$, any product of $n$ elements of $\mathrm{Lie}\langle\!\langle m\rangle\!\rangle$ is zero, so $\mathrm{Lie}\langle\!\langle m\rangle\!\rangle / \mathrm{Lie}_n\langle\!\langle m\rangle\!\rangle$ is $n$-nilpotent. If $n > 1$, then is not $(n-1)$-nilpotent $[X_0, [X_0, \ldots [X_0, X_1] \ldots]] \in \mathrm{Lie}\langle\!\langle m\rangle\!\rangle \setminus \mathrm{Lie}_n\langle\!\langle m\rangle\!\rangle$. The family $(u[X_v])_{u \in m^{(\star)}, v \in m^{\star}}$ is clearly a finite generating family of $\mathrm{Lie}\langle\!\langle m\rangle\!\rangle / \mathrm{Lie}_n\langle\!\langle m\rangle\!\rangle$. $\square$

An $\mathcal{L}_k$-formula is said positive if the quantifier-free part of its prenex normal form contains no negation. We write $\mathcal{L}_k^{+,1}$ for the set of $\mathcal{L}_k$-formulas of the form $\forall \bar{x} \exists! \bar{y}(\theta(\bar{x}, \bar{y}))$ where $\theta(\bar{x}, \bar{y})$ is positive and quantifier-free.

**Proposition 6.4.** *Let $\varphi = \forall \bar{x} \exists! \bar{y}(\theta(\bar{x}, \bar{y}))$ be a sentence in $\mathcal{L}_k^{+,1}$ that holds in all finite dimensional nilpotent mixed structures. Then $\forall \bar{x} \exists \bar{y}(\theta(\bar{x}, \bar{y}))$ holds in all mixed structures.*

**Proof.** Write $\varphi = \forall \bar{x} \exists \bar{y}(\theta(\bar{x}, \bar{y}))$ as above. Write $l, m$ for the respective sizes of the variable $\bar{x}$ and $\bar{y}$ in $\theta$. Let $n > 0$. By Lemma 6.3, the sentence $\varphi$ holds in $\mathrm{Lie}^{/n}\langle\!\langle m\rangle\!\rangle := \mathrm{Lie}\langle\!\langle m\rangle\!\rangle / \mathrm{Lie}_n\langle\!\langle m\rangle\!\rangle$. Now let $\overline{P} \subseteq \mathrm{Lie}\langle\!\langle m\rangle\!\rangle$ and write $\overline{P_n}$ for the tuple of reductions modulo $\mathrm{Lie}_n\langle\!\langle m\rangle\!\rangle$. The set $\theta(\overline{P_n}, \mathrm{Lie}^{/n}\langle\!\langle m\rangle\!\rangle)$ is a singleton $\{\overline{a_n}\}$. Now for each tuple $\overline{Q_{n+1}} \in \theta(\overline{P_{n+1}}, \mathrm{Lie}^{/n+1}\langle\!\langle m\rangle\!\rangle)$, since $\theta$ is a positive quantifier-free formula and the projection $\varphi_n: \mathrm{Lie}^{/n+1}\langle\!\langle m\rangle\!\rangle \longrightarrow \mathrm{Lie}^{/n}\langle\!\langle m\rangle\!\rangle$ is an $\mathcal{L}_k$-homomorphism, we have $\overline{Q_n} \in \theta(\overline{P_n}, \mathrm{Lie}^{/n}\langle\!\langle m\rangle\!\rangle)$, so $\varphi_n(\overline{a_{n+1}}) = \overline{a_n}$. We deduce that we have we have $\overline{a_n} = \overline{Q} + \mathrm{Lie}_n\langle\!\langle m\rangle\!\rangle$ for a tuple $\overline{Q} \in \mathrm{Lie}\langle\!\langle m\rangle\!\rangle^l$. As any $\mathcal{L}_k$-term $t(\bar{x}, \bar{y})$ occurring in $\theta(\bar{x}, \bar{y})$ satisfies

$$t(\overline{P}, \overline{Q}) = t(\overline{P}, \overline{Q}) - t(\overline{P_n}, \overline{Q} + \mathrm{Lie}_n\langle\!\langle m\rangle\!\rangle) \in \mathrm{Lie}_n\langle\!\langle m\rangle\!\rangle$$



for all $m \in \mathbb{N}$, we have $t(\overline{P}, \overline{Q}) = 0$. Let $L$ be a mixed structure. Let $\overline{a} \subseteq L^l$. Let there is a $\overline{Q} \in \theta((X_0, \ldots, X_{l-1}), \mathrm{Lie}\langle\!\langle l \rangle\!\rangle)$. The relations $\theta(\overline{a}, \mathrm{ev}_{\overline{a}}(\overline{Q}))$ are satisfied as $\theta$ is quantifier-free and positive and $\mathrm{ev}_{\overline{a}}$ is an $\mathcal{L}_k$-homomorphism. This shows that $\forall \overline{x} \exists \overline{y}(\theta(\overline{x}, \overline{y}))$ holds in $L$. □

**Proposition 6.5.** *Let $\varphi \in \mathcal{L}_k$ be a universal sentence that holds in all nilpotent finite dimensional mixed structures. Then $\varphi$ holds in all mixed structures.*

**Proof.** We may assume that $\varphi = \forall \overline{x}(\bigvee_{i=1}^r \bigwedge_{j=1}^s \theta_{i,j}(\overline{x}))$ is in disjunctive-conjunctive prenex normal form, where each $\theta_{i,j}(\overline{x})$ is atomic. Let $l \in \mathbb{N}$ be the length of $\overline{x}$. Each $\theta_{i,j}(\overline{x})$ has the form $P_{i,j}(\overline{x}) = 0$ or $P_{i,j}(\overline{x}) \neq 0$ for a $P_{i,j} \in k\langle\!\langle l \rangle\!\rangle$. Let $L$ be a mixed structure and let $\overline{a} \in L^l$. Consider the kernel $\mathfrak{m}$ of the evaluation map $\mathrm{ev}_{\overline{a}}^{\mathrm{Lie}} : \mathrm{Lie}\langle\!\langle l \rangle\!\rangle \longrightarrow L$. Let $n > 0$, and write $\mathfrak{m}_n$ for the summability ideal of $\mathrm{Lie}\langle\!\langle l \rangle\!\rangle$ generated by $\mathfrak{m}$ and $\mathrm{Lie}_n\langle\!\langle l \rangle\!\rangle$, and $\mathfrak{q}_n$ for the summability ideal of $\mathrm{Lie}\langle\!\langle l \rangle\!\rangle / \mathfrak{m}$ generated by $\mathfrak{m}/(\mathfrak{m} \cap \mathrm{Lie}_n\langle\!\langle l \rangle\!\rangle)$. By Lemma 6.3, the Lie algebra $\mathrm{Lie}\langle\!\langle l \rangle\!\rangle / \mathfrak{m}_n$ is nilpotent and finite dimensional. Write $\pi_n$ for the projection $\mathrm{Lie}\langle\!\langle l \rangle\!\rangle / \mathfrak{m} \longrightarrow \mathrm{Lie}\langle\!\langle l \rangle\!\rangle / \mathfrak{m}_n$.

Let $\overline{X} = (X_0 + \mathfrak{m}, \ldots, X_{l-1} + \mathfrak{m}) \in (\mathrm{Lie}\langle\!\langle l \rangle\!\rangle / \mathfrak{m})^l$. There is an $i_n \in \{1, \ldots, r\}$, which we choose minimal, such that all $\theta_{i_n,j}(\pi_n(\overline{X}))$ holds in $\mathrm{Lie}\langle\!\langle l \rangle\!\rangle / \mathfrak{m}_n$ for each $j \in \{1, \ldots, s\}$. We may pick a strictly increasing sequence $\psi : \mathbb{N} \longrightarrow \mathbb{N}$ such that $(i_{\psi(n)})_{n \in \mathbb{N}}$ is constant with constant value $i \in \{1, \ldots, r\}$. Let $j \in \{1, \ldots, s\}$ and $n \in \mathbb{N}$. Write $R_{i,j,n} = P_{i,j} - t_{\psi(n),\mathrm{Lie}}$ where $t_{\psi(n),\mathrm{Lie}}$ is as in Lemma 6.2 for $P_{i,j}$. We have $P_{i,j}(\pi_{\psi(n)}(\overline{X})) = t_{\psi(n),\mathrm{Lie}}(\pi_{\psi(n)}(\overline{X}))$. Now

$$\begin{aligned} P_{i,j}(\overline{X}) &= t_{\psi(n),\mathrm{Lie}}(\overline{X}) + R_{i,j,n}(\overline{X}) \\ &\in P_{i,j}(\pi_{\psi(n)}(\overline{X})) + \mathfrak{q}_{\psi(n)}. \qquad \text{(by construction of } t_{\psi(n),\mathrm{Lie}} \text{ in Lemma 6.2)} \end{aligned}$$

If $\theta_{i,j}$ is $P_{i,j}(\overline{x}) = 0$, then $t_{\psi(n),\mathrm{Lie}}(\pi_{\psi(n)}(\overline{X})) \in \mathfrak{q}_n$, so $P_{i,j}(\overline{X}) \in \mathfrak{q}_{\psi(n)}$. As this holds for all $n > 0$, it follows that $P_{i,j}(\overline{X}) = 0$ in $\mathrm{Lie}\langle\!\langle l \rangle\!\rangle / \mathfrak{m}$. If $\theta_{i,j}$ is $\neg(P_{i,j}(\overline{x}) = 0)$, then $P_{i,j}(\overline{X}) \notin \mathfrak{q}_{\psi(1)}$ so $P_{i,j}(\overline{X}) \neq 0$ in $\mathrm{Lie}\langle\!\langle l \rangle\!\rangle / \mathfrak{m}$. Thus $\theta_{i,j}(\overline{X})$ holds in $\mathrm{Lie}\langle\!\langle l \rangle\!\rangle / \mathfrak{m}$ for all $j \in \{1, \ldots, s\}$.

The map $\mathrm{ev}_{\overline{a}}^{\mathrm{Lie}} : \mathrm{Lie}\langle\!\langle l \rangle\!\rangle / \mathfrak{m} \longrightarrow L$ is an injective morphism, so $\theta_{i,j}(\mathrm{ev}_{\overline{a}}^{\mathrm{Lie}}(\overline{X}))$ holds in $L$ for all $j \in \{1, \ldots, s\}$, i.e. $\bigwedge_{j=1}^s \theta_{i,j}(\overline{a})$ holds. This shows that $L \vDash \varphi$. □

**Corollary 6.6.** *Let $\varphi \in \mathcal{L}_k$ be a sentence of the form $\forall \overline{x} \exists! \overline{y}(\theta(\overline{x}, \overline{y}))$, where $\theta$ is positive and quantifier-free, that holds in all nilpotent finite dimensional mixed structures. Then $\varphi$ holds in all mixed structures.*

**Proof.** By Proposition 6.4, the sentence $\forall \overline{x} \exists \overline{y}(\theta(\overline{x}, \overline{y}))$ holds in all mixed structures. Note that the universal sentence $\phi : \forall \overline{x} \forall y \forall z (t(\overline{x}, y) = t(\overline{x}, z) = 0 \to y = z)$ holds in $\mathcal{G}$, hence in all mixed structures by Proposition 6.5. Now $\phi \wedge \psi$ is equivalent to $\varphi$ in all mixed structures, hence the result. □

## 6.3 Non-singular equations

Let $l, m \in \mathbb{N}$. For $f \in \mathrm{Lie}\langle\!\langle l+1 \rangle\!\rangle$ and $g \in \mathrm{Lie}\langle\!\langle m+1 \rangle\!\rangle$, we write $f \circ g = \mathrm{ev}_{(X_0, \ldots, X_{l-1}, g)}^{\mathrm{Lie}}(f) \in \mathrm{Lie}\langle\!\langle (l+1) \cup (m+1) \rangle\!\rangle$. By Proposition 3.33, for $l = m$, this defines an associative law on $\mathrm{Lie}\langle\!\langle l+1 \rangle\!\rangle$ with identity $X_l$. The same thing applies to the composition law $f \circ g = \mathrm{ev}_{(X_0, \ldots, X_{l-1}, g)}^{\mathrm{Ass}}(f)$ on $k\langle\!\langle l+1 \rangle\!\rangle$.

By the universal property of coproducts, there is a unique strongly linear homomorphism $\mathcal{E} : \mathrm{Lie}\langle\!\langle l+1 \rangle\!\rangle \longrightarrow k$ that is zero on $\mathcal{K}_0, \ldots, \mathcal{K}_{l-1}$ and the identity on $\mathcal{K}_l$. Given $g \in \mathrm{Lie}\langle\!\langle l+1 \rangle\!\rangle$, the functions $f \mapsto \mathcal{E}(g) \mathcal{E}(f)$ and $f \mapsto \mathcal{E}(f \circ g)$ satisfy the same universal propery, so they coincide. In other words $\mathcal{E} : (\mathrm{Lie}\langle\!\langle l+1 \rangle\!\rangle, \circ, X_l) \longrightarrow (k, \cdot, 1)$ is a morphism of monoids. By Lemma 6.2, for $t \in \mathrm{Lie}\langle\!\langle l+1 \rangle\!\rangle$, there is a $c > 0$ such that $\mathcal{E}(t) = \mathcal{E}(\widehat{t_{n,\mathrm{Lie}}})$ for all $n \geqslant c$. Indeed, as $\mathcal{E}$ is strongly linear for the minimal summability structure on $k$, otherwise by (6.2) the family $(\mathcal{E}(\widehat{t_{c+1,\mathrm{Lie}}} - \widehat{t_{c,\mathrm{Lie}}}))_{c>0}$ would not be finitely supported.



**Theorem 6.7.** *Let $L$ be a mixed structure, let $l \in \mathbb{N}$, $\bar{a} \in L^l$ and $t \in \mathrm{Term}(l,1)$ with $\mathcal{E}(t) \neq 0$. Then there is a unique element $b \in L$ with $t(\bar{a}, b) = 0$.*

**Proof.** Let $c > 0$, let $\mathcal{G}$ be a nilpotent mixed structure seen as a multipliability exponential group, and let $\bar{a} \in \mathcal{G}^l$. Note since $k$ is a field that $\mathcal{G}$ is residually torsion-free nilpotent in the sense of [2, Section 7]. We have $\mathcal{E}(t) = \mathcal{E}(\widehat{t_{n,\mathrm{Gr}}})$ for a certain $n \geqslant c$ by Lemma 6.2 and the arguments that precede this corollary. The morphism $\mathcal{E}$ restricted to the coproduct in **Gr** of $\mathcal{K}_1, \ldots, \mathcal{K}_l$ and $\mathcal{K}_{l+1}$ is the morphism $\alpha$ defined in [2, Section 7.1] (as they satisfy the same universal property), i.e. $\alpha(t_{n,\mathrm{Gr}}) \in k^\times$. So by [2, Theorem 7.10], there is a unique $b \in \mathcal{G}$ with $t_{n,\mathrm{Gr}}(\bar{a}, b) = 1$, i.e. $\mathcal{G} \vDash \forall \bar{x} \exists y (t(\bar{x}, y) = 0 \wedge t)$ where $\varphi = \forall \bar{x} \exists y (t(\bar{x}, y) = 0)$ is in $\mathcal{L}_k^{+,1}$. We conclude with Corollary 6.6. $\square$

**Corollary 6.8.** *Let $\mathcal{G}$ be a multipliaility exponential group with evaluations, let $m \in \mathbb{N}$ and $g_1, \ldots, g_n \in \mathcal{G}$ and $\lambda_1, \ldots, \lambda_m \in k$ with $\lambda_1 + \cdots + \lambda_m \neq 0$. Then there is a unique $f \in \mathcal{G}$ with*

$$g_1 f^{\lambda_1} \cdots g_n f^{\lambda_n} = 1.$$

## 6.4 The Hall-Petresco identities

Let $n \geqslant 2$ and write $\mathrm{Free}_\mathbb{Z}(n)$ for the free pure group over a set of $n$ elements. The Hall-Petresco identities for groups, consequences of Hall's commutator collection formula [21, Theorem 3.1] were proved by Petresco [39, 5.3] for $n = 2$, and by Warfield [47, Theorem 6.1] for arbitary $n$. They give an expansion of the product $g_0^m \cdots g_{n-1}^m$ of powers where $m \geqslant 1$ and $g = (g_0, \ldots, g_{n-1}) \in \mathcal{G}^n$ as

$$g_0^m \cdots g_{n-1}^m = (g_0 \cdots g_{n-1})^m \tau_2(g)^{\binom{m}{2}} \cdots \tau_{m-1}(g)^{\binom{m}{m-1}} \tau_m(g)$$

where $\tau_2(x_0, \ldots, x_{n-1}) \in \mathrm{Free}_\mathbb{Z}(n)_2, \ldots, \tau_m(x_0, \ldots, x_{n-1}) \in \mathrm{Free}_\mathbb{Z}(n)_m$ do not depend on $\mathcal{G}$.

We now assume that $k$ is a domain and we fix $c \geqslant 1$. The Hall-Petresco identities for an exponential group $\mathcal{G}$ over $k$ that is $c$-nilpotent as a pure group are the sentences

$$\forall g_0, \ldots, g_{n-1} \left( g_0^\lambda \cdots g_{n-1}^\lambda = (g_0 \cdots g_{n-1})^\lambda \tau_2(g)^{\binom{\lambda}{2}} \cdots \tau_{c-1}(g)^{\binom{\lambda}{c-1}} \right),$$

where $\lambda \in k$. They are taken as part of the definition of $k$-exponential groups by Warfield [47, Definition 10.4].

**Proposition 6.9.** *If $k$ is a domain, then any nilpotent exponential group satisfies the Hall-Petresco identities. In particular, nilpotent exponential groups are exactly nilpotent $k$-groups as per [47, Definition 10.19].*

**Proof.** Let $\mathcal{G}$ be a nilpotent exponential group of nilpotenty class $n \geqslant 1$, let $g : n \longrightarrow \mathcal{G}$ be a finite family and let $\lambda \in k$. The exponential group $1 + k \langle\!\langle n \rangle\!\rangle$ is [47, Theorem 10.23] a complete filtered $k$-group in the sense of [47, Definition 10.19], meaning that for all $Q_0, \ldots, Q_{n-1} \in 1 + k \langle\!\langle n \rangle\!\rangle$, we have

$$Q_0^\lambda \cdots Q_{n-1}^\lambda = (g_0 \cdots g_{n-1})^\lambda \tau_2(g)^{\binom{\lambda}{2}} \cdots \tau_l(g)^{\binom{\lambda}{l}} \cdots \tag{6.3}$$

(see Remark 4.15). Applying this to $(Q_0, \ldots, Q_{n-1}) = (\exp(X_0), \ldots, \exp(X_{n-1})) \in \mathrm{Gr} \langle\!\langle n \rangle\!\rangle$ and evaluating at $g$ with the evaluation map $\mathrm{ev}_g^{\mathrm{Gr}}$, we obtain the corresponding Hall-Petresco formula. $\square$



Let $(\tilde{\mathbb{L}}^{>\mathbb{R}}, \circ, x, <)$ be the linearly ordered group under composition of finitely nested hyperseries of [5]. It contains the group $\mathcal{T}$ of Corollary 4.27 as a subgroup. For any positive $f \in \tilde{\mathbb{L}}^{>\mathbb{R}}$, there is [5, Proposition 4] a unique ordered group isomorphism

$$(\mathbb{R}, +, 0, <) \longrightarrow (\mathcal{C}(f), \circ, x, <); \lambda \mapsto f^{[\lambda]}$$

onto the centraliser $\mathcal{C}(f)$ of $f$ that sends 1 to $f$. Therefore $(\mathbb{R}, +, 0, <)$ is interpretable in $(\tilde{\mathbb{L}}^{>\mathbb{R}}, \circ, x, <)$. Moreover $(\lambda, f) \mapsto f^{[\lambda]}$ is a power map on $\tilde{\mathbb{L}}^{>\mathbb{R}}$ for which it is an exponential group (see [5, Section 4.5]).

**Proposition 6.10.** *The power map* $\mathbb{R} \times \tilde{\mathbb{L}}^{>\mathbb{R}} \longrightarrow \tilde{\mathbb{L}}^{>\mathbb{R}}$ *is interpretable in* $(\tilde{\mathbb{L}}^{>\mathbb{R}}, \circ, x)$.

**Proof.** Any two positive elements of $\tilde{\mathbb{L}}^{>\mathbb{R}}$ are conjugate [5, Theorem 3]. The ordering $<$ is definable with one positive parameter $f \in \tilde{\mathbb{L}}^{>\mathbb{R}}$, as positive elements are exactly conjugates of $f$. Now $(\tilde{\mathbb{L}}^{>\mathbb{R}}, \circ, x, <)$ is [5, Proposition 5] a growth order group as per [6, Definition 2.20], so it has [2, Section 4.4] definable c-dominance relations $\preccurlyeq$ and $\prec$ as per [2, Definition 2.1].

If $f \in \mathcal{T}$, then $f^{[\lambda]}$ coincides with the power $f^\lambda$ in the exponential group $\mathcal{T}$. Indeed consider map sending $\lambda \in \mathbb{R}$ to the unique $\rho \in \mathbb{R}$ with $f^\lambda = f^{[\rho]}$. That $f^\lambda \in \mathcal{C}(f)$ is a general property [33, Property 2] of exponential groups. By **EG3** in $\mathcal{T}$ and [5, Propositions 4.27 and 4.28] in $\tilde{\mathbb{L}}^{>\mathbb{R}}$, this map is a ring endomorphism of $\mathbb{R}$, so it must be the identity. Pick non-commuting $f, g \in \mathcal{T}$ with $f > \mathcal{C}(g)$. We have $[\![f, g]\!] > 1$ by the property [6, **GOG2**] of growth order groups. For $\lambda, \mu \in \mathbb{R}$, the infinite Hall-Petresco identities (6.3) evaluated in $\mathcal{T}$ give

$$[f^{[\lambda]}, g^{[\mu]}] = [f, g]^{[\lambda \mu]} h$$

where $h$ is an infinite product of powers of commutators in $f$ and $g$ of length $\geqslant 3$. Form [2, (6.3–6.5)], we see that $[s, t] \prec s, t$ whenever $s, t \in \mathcal{T}$ are non-trivial. Indeed either $\operatorname{res}(s^{-1} t^{-1}) \prec \operatorname{res}(s), \operatorname{res}(t)$ and $\operatorname{res}(s\, t) \prec \operatorname{res}(f), \operatorname{res}(g)$, or $\operatorname{res}(s^{-1} t^{-1}) = -(\operatorname{res}(s) + \operatorname{res}(t))$ and $\operatorname{res}(s\, t) = \operatorname{res}(s) + \operatorname{res}(t)$, or $\operatorname{res}(s^{-1} t^{-1}) = -\operatorname{res}(s\, t) \in \{\operatorname{res}(s), \operatorname{res}(t)\}$. It follows that $[f, g] \prec f, g$ by [2, (6.3–6.4)] in the first case, and by [2, (6.5)] in the second and third case. We deduce by [3, Lemma 6.10] that $h \prec [f, g] \preccurlyeq [f, g]^{[\lambda \mu]}$. So $[f, g]^{[\lambda \mu]}$ is the unique element of $\mathcal{C}([f, g])$ such that $[f^{[\mu]}, g^{[\lambda]}] ([f, g]^{[\lambda \mu]})^{-1} \prec [f^{[\mu]}, g^{[\lambda]}]$. Now fixing $h_f, h_g \in \tilde{\mathbb{L}}^{>\mathbb{R}}$ with $h_f [f, g] h_f^{-1} = f$ and $h_g [f, g] h_g^{-1} = g$, the corresponding conjugation maps are definable automorphisms of $\tilde{\mathbb{L}}^{>\mathbb{R}}$ that send $[f, g]^{[\lambda]}$ to $f^{[\lambda]}$ and $g^{[\lambda]}$ respectively for each $\lambda \in \mathbb{R}$. Thus the operation $([f, g]^{[\lambda]}, [f, g]^{[\mu]}) \longmapsto [f^{[\lambda]}, g^{[\mu]}]$ is definable, and so is $\star: ([f, g]^{[\lambda]}, [f, g]^{[\mu]}) \longmapsto [f, g]^{[\lambda \mu]}$. This shows that $(\mathbb{R}, +, \cdot, 0, 1)$ is isomorphic to the definable set $(\mathcal{C}([f, g]), \circ, \star, x, [f, g])$. Now for $[f, g]^{[\lambda]} \in \mathcal{C}([f, g])$ and positive $h \in \tilde{\mathbb{L}}^{>\mathbb{R}}$, the element $h^{[\lambda]}$ is $\varphi^{-1} \circ ([f, g]^{[\lambda]}) \circ \varphi$ where $\varphi$ is any element that conjugates $h$ to $[f, g]$. For $h$ negative, we define $h^{[\lambda]}$ as $(h^{\mathrm{inv}})^{[-\lambda]}$, and we set $x^{[\lambda]} = x$. Thus, *via* the definable map

$$\begin{aligned}\mathcal{C}([f, g]) \times \tilde{\mathbb{L}}^{>\mathbb{R}} &\longrightarrow \tilde{\mathbb{L}}^{>\mathbb{R}} \\ ([f, g]^{[\lambda]}, h) &\longmapsto h^{[\lambda]},\end{aligned}$$

the power map is interpretable. $\square$

## 6.5 Central formula

Let $(L, +, 0, \cdot, (a \mapsto \lambda a)_{\lambda \in k}, \Sigma)$ be a summability Lie algebra with evaluations. For $a \in L$, we write $\mathrm{ad}_a$ for the derivation $L \longrightarrow L; b \mapsto a \cdot b$ on $L$.



**Proposition 6.11.** *For all $a, b \in L$ we have $\exp(\mathrm{ad}_a)(b) = a * b * (-a)$.*

**Proof.** As $L$ has evaluations, the map $\exp(\mathrm{ad}_a) = \sum_{i \in \mathbb{N}} \frac{\mathrm{ad}_a^{\circ i}}{i!}$ is well-defined and strongly linear. The identity

$$\forall b \left( \sum_{i \in \mathbb{N}} \frac{\overset{i \text{ times}}{a \cdot (a \cdot (\cdots (a \cdot b)))}}{i!} = b * a * (-a) \right) \tag{6.4}$$

is a universal sentence in $\mathcal{L}_k$ that holds [9, Chapter 3, § 6 Corollaire 3] in all nilpotent finite dimensional structures for $k = \mathbb{C}$. Thus it holds in $\mathbb{C}\langle\!\langle 2 \rangle\!\rangle$ by Proposition 6.5, whence in its summability subalgebra $\mathbb{Q}\langle\!\langle 2 \rangle\!\rangle$. The reduct of $L$ in $\mathbb{Q}$ has evaluations (see Remark 3.7), so evaluating from $\mathbb{Q}\langle\!\langle 2 \rangle\!\rangle$ into $L$ yields (6.4) in $L$. $\square$

## 6.6 Engel conditions

Let $\mathcal{G}$ be a group and let $L$ be a Lie algebra. The group $\mathcal{G}$ is said $n$-Engel if for all $f, g \in \mathcal{G}$, we have

$$[\![f, [\![f, [\![f, [\![\ldots [\![f, g]\!] \ldots ]\!]]\!]]\!]]\!] = 1.$$
$$\underset{n \text{ occurrences of } f}{}$$

The Lie algebra $L$ is said $n$-Engel if for all $a, b \in L$ we have $\mathrm{ad}_a^{\circ n}(b) = 0$. Clearly, any nilpotent Lie algebra of nilpotency class $\leqslant n$ is $n$-Engel. Zel'manov proved [48] that a converse holds:

**Zel'manov's theorem.** [48, Theorem] *For $n > 0$, all $n$-Engel lie algebras over a field of characteristic $0$ are nilpotent.*

**Proposition 6.12.** *Let $\mathcal{G}$ be a multipliability exponential group with evaluations and let $n > 0$. If $\mathcal{G}$ is $n$-Engel, then it is nilpotent.*

**Proof.** We work in the $\mathbb{Q}$-reduct $\mathcal{G}_\mathbb{Q}$ of $\mathcal{G}$ which, by Remark 3.7 and the formal Lie correspondence, is a multipliability exponential group with evaluations.

For $i > 0$, consider the non-associative word of rank $i$ defined by $u_i = (0(0(\ldots (01)\ldots))) \in 2^{(\star)}$. For $m > 0$, let $B_m(x_0, y_1) := u_m[(X_0, X_1)] \in \mathrm{Lie}\langle\!\langle 2 \rangle\!\rangle$ and $C_m := u_m[\![(X_0, X_1)]\!] \in \mathrm{Lie}\langle\!\langle 2 \rangle\!\rangle$. We claim that the universal $\mathcal{L}_k$-sentence $\psi_m : \forall x_0 \forall x_1 (C_m(x_0, x_1) = 0 \Longrightarrow B_m(x_0, x_1) = 0)$ holds in all finite dimensional nilpotent mixed structures over $\mathbb{Q}$. Indeed this holds for $m = 1$ by (5.5). Let $m > 0$ such that $\psi_m$ always holds. Let $\mathcal{H}$ be a nilpotent mixed structure of finite dimension over $\mathbb{Q}$ and $f, g \in \mathcal{H}$ with $C_{m+1}(f, g) = 0$, we have $B_{m+1}(f, g) = B_m(f, B_1(f, g))$. Now one sees by induction that for all $h \in \mathcal{H}$, the element $B_m(f, h)$ is a finite product of powers of elements $w[\![(f, h)]\!]$ in which $u_m$ is a subword of $w$. Therefore $B_{m+1}(f, g)$ is a product of powers of elements $w[\![(f, g)]\!]$ in which $u_{m+1}$ is a subword of $w$. As $u_{m+1}[\![(f, g)]\!] = C_{m+1}(f, g) = 0$, we deduce that $B_{m+1}(f, g) = 0$. So $\psi_{m+1}$ holds in $\mathcal{H}$. We deduce with Proposition 6.5 that each $\psi_m$ for $m > 0$ holds in $\mathcal{G}_\mathbb{Q}$, hence in $\mathcal{G}$. This shows that if $\mathcal{G}$ is $n$-Engel, then $\mathrm{Lie}(\mathcal{G}_\mathbb{Q})$ is $n$-Engel, hence nilpotent by Zelmanov's theorem. By Theorem 5.26(a), the underlying exponential group of $\mathcal{G}$ is nilpotent. $\square$

**Acknowledgments.** We thank Christian D'Elbée and Adrien Deloro for answering our questions. We are also greatful to Isabel Müller who first told us about the Mal'cev and Lazard's correspondences.